\documentclass{article}
\usepackage[T1]{fontenc}

\usepackage{geometry}
\usepackage{authblk}

\usepackage{subfig}

\usepackage{comment}
\usepackage[english]{babel}
\usepackage{xcolor}

\usepackage{graphicx}
\usepackage{latexsym}

\usepackage{amsthm}
\usepackage{mathtools}

\newtheorem{remark}{Remark}
\newtheorem{theorem}{Theorem}
\newtheorem{lemma}{Lemma}
\newtheorem{definition}{Definition}
\usepackage{amssymb}

\renewcommand{\appendix}{
\setcounter{section}{0}
\renewcommand{\thesection}{\Roman{section}}
\vspace{0.5cm}
{\Large{\bf APPENDIX}}}

\newcommand{\T}{\mathcal{T}}

\newcommand{\Vv}{\mathcal{V}}

\newcommand{\B}{\mathcal{B}}
\newcommand{\C}{\mathcal{C}}
\newcommand{\D}{\mathcal{D}}

\newcommand{\RR}{\mathcal{R}}

\newcommand{\Rbb}{\mathbb{R}}

\newcommand{\nn}{{\bf n}}
\newcommand{\mm}{{\bf m}}
\newcommand{\uu}{{\bf u}}
\newcommand{\xb}{{\bf x}}

\newcommand{\ww}{{\bf w}}
\newcommand{\vv}{{\bf v}}

\newcommand{\fb}{{\bf f}}
\newcommand{\gb}{{\bf g}}
\newcommand{\hb}{{\bf h}}

\newcommand{\Fb}{{\bf F}}

\newcommand{\Ss}{{\mathcal{ S}}}

\newcommand{\pf}{{\mathfrak{p}}}

\newcommand{\Gf}{{\mathfrak{G}}}
\newcommand{\Ff}{{\mathfrak{F}}}
\newcommand{\Hf}{{\mathfrak{H}}}
\newcommand{\Df}{{\mathfrak{D}}}

\newcommand{\Ec}{{\mathcal{E}}}

\usepackage{mathtools}

\title{Analysis of an alternative Navier-Stokes system: Weak entropy solutions and a convergent numerical scheme}
\author{Magnus Sv{\"a}rd}

  \affil{Dept. of Mathematics, University of Bergen, P.O. Box 7803, 5020 Bergen, Norway, \\ Magnus.Svard@uib.no}

\date{\today}
  
\begin{document}

\maketitle



\begin{abstract}

We consider an alternative Navier-Stokes model for compressible viscous ideal gases, originally proposed in \cite{Svard18}. We derive a priori estimates that are sufficiently strong to support a weak entropy solution of the system.  Guided by these estimates, we propose a finite volume scheme, derive the analogous estimates and demonstrate grid convergence towards a weak entropy solution. Furthermore, this existence proof is valid for ``large'' initial data and no a priori assumptions on the solution are needed.
\end{abstract}

\section{Introduction}

The importance of modelling viscous and heat conducting compressible fluids in physics and engineering can not be understated. Such fluids are typically modelled by the well-known Navier-Stokes-Fourier (NSF) equations.  The complexity of these equations often necessitates numerical methods to generate approximate solutions. However, the design of numerical schemes is severely hampered by the lack of well-posedness results for the NSF system. This invariably leads to uncertainties pertaining to the validity of numerical results obtained with state-of-the-art codes. (See e.g. \cite{GassnerSvard21}.) 

 The NSF equations model the dynamics of (Lagrangian) mass elements by Newton's 2nd Law, which by definition precludes mass diffusion. (See \cite{Svard18,Stokes1845}.) The effect of the random motions at the molecular level enters the NSF system via the stress tensor and a model of viscosity that exerts a force on each mass element. Furthermore, heat diffusion is modelled by Fourier's law in the energy equation. A consequence of the Lagrangian perspective is that the continuity equation becomes hyperbolic, which poses a significant mathematical challenge when trying to develop a well-posedness theory.

In recent years, many authors have, for physical reasons, proposed alternative models. See e.g. \cite{Reddy_etal19} where a class of mass-diffusive continuum models is derived and demonstrated to be more accurate than the NSF system in many cases. Another well-known mass-diffusive  modification of the NSF system was given in 
\cite{Brenner05_1,Brenner05_2}. For a constant mass-diffusive coefficient, this system was proven to have  weak solutions in \cite{FeireislVasseur10}.

Herein, we consider the alternative Navier-Stokes model proposed in \cite{Svard18}. Let $\xb=(x,y,z)$ be the Cartesian coordinates and $\Omega$ a bounded spatial domain. Then the model is given by
\begin{align}
\partial_t \rho + \nabla\cdot(\rho \vv )&= \nabla\cdot  (\nu \nabla \rho),\nonumber \\
\partial_t (\rho \vv) + \nabla\cdot(\rho \vv \otimes \vv) + \nabla p&=\nabla\cdot 
(\nu \nabla \rho \vv)), \qquad t\in [0,\T],\,\xb\in\Omega\subset \Rbb^3 \label{eulerian_ideal} \\
\partial_t (E) + \nabla\cdot(E \vv +p\vv)  &= \nabla\cdot( \nu \nabla E), \nonumber \\ 
p&=\rho R T, \quad \textrm{(ideal gas law).}\nonumber 
\end{align}
 $\rho,\rho \vv, E$ are the conserved variables density, momentum and total energy.  $\vv=(v_1,v_2,v_3)^T=(u,v,w)$ is the velocity vector. $p,T$ are pressure and temperature. Furthermore, the total energy is given by $E=\frac{p}{\gamma-1}+\frac{\rho |\vv|^2}{2}$ where $\gamma=c_p/c_v$ and $c_{p,v}$ are the heat capacities at constant pressure or volume. The model is (at least) valid for ideal gases with $1<\gamma\leq 5/3$. The conservative variables are collected in the vector $\uu=(\rho,\mm^T,E)^T$ where $\mm=(m_1,m_2,m_3)^T=\rho \vv$ is the momentum vector. $R$ is the gas constant and $\nu=\nu(\rho,T)$ is the diffusion coefficient.

An ideal gas is one where the molecules bounce elastically against each other and have no rotational energy. Furthermore, the mean free path is large relatively to their sizes. In such a gas, the range of the diffusion is proportional to the mean free path  which scales as the inverse of density. Hence, it was proposed in \cite{Svard18} that, 
\begin{align}
\nu =\frac{\mu_0}{\rho},\label{ideal_nu}
\end{align}
where $\mu_0$ is the dynamic viscosity.  This value results in viscous terms that resemble those of the Navier-Stokes equations, and Fourier's law appears as a diffusive term in the energy equation of (\ref{eulerian_ideal}). In \cite{DolejsiSvard21} and \cite{SayyariDalcin21}, it was shown, in a suite of problems, that the new system (\ref{eulerian_ideal})-(\ref{ideal_nu}) produces solutions that are next to indistinguishable from those of the NSF system. (The differences are much less than what can be measured in experiments.)

\begin{remark}
  With $\mu_0$ being the dynamic viscosity, the value of $\nu$ is identical to the kinematic viscosity. However, we emphasise that physically $\nu$ is a diffusion coefficient, and not a viscosity coefficient.
\end{remark}

Unlike the NSF system, (\ref{eulerian_ideal}) is derived in an Eulerian\footnote{The system (\ref{eulerian_ideal}), was termed ``Eulerian model'' in \cite{Svard18}. Here, we call it an ``alternative Navier-Stokes model'' since the original name has caused some confusion that it is a (numerical) regularisation of the Euler equations rather than an alternative to the standard Navier-Stokes model. } frame by directly appealing to the principles of conservation of mass, momentum, and energy. In the Eulerian perspective, the random motions of molecules are naturally modelled as diffusion in all three conservative variables. (Indeed, (\ref{eulerian_ideal}) differs from the NSF system only in the right-hand side.)

Our goal is to demonstrate that the system (\ref{eulerian_ideal}) has weak entropy solutions by proving that solutions to a finite volume scheme converge as the grid is refined. In this effort, we have benefited greatly from \cite{FeireislVasseur10} and have used many similar techniques and arguments.

\subsection{Adding more physics to the model}

Although (\ref{eulerian_ideal}) is accurate in the regime it is intended to model, considerations for reduced models (isothermal and isentropic) reveal that (\ref{eulerian_ideal}) appears to lack mechanisms that prevent concentration of mass. Hence, we propose some minor modifications to (\ref{eulerian_ideal}) below. These modifications do \emph{not} extend validity range of the model and can thus alternatively be viewed as \emph{technical assumptions}.

Nevertheless, we continue with a short physical discussion: The choice $\nu=\mu_0/\rho$ models diffusion whereby molecules of the ideal gas travel relatively freely in space.   Their random crossings of (Eulerian) control-volume boundaries, is modelled as diffusion on the global scale and the more rarefied the gas is, the more significant is the diffusion relatively to the convection. This is consistent with the diffusion coefficient being proportional to the mean free path, i.e., $\rho^{-1}$. However, this mechanism is not prevalent when the gas becomes denser since molecules will bounce into other molecules before they have travelled any significant distance. The latter process is conductive rather than diffusive and not accounted for with $\nu=\mu_0/\rho$. Hence, the physics suggests that the model should be augmented with collisional/frictional diffusion.

There are several options to account for collisional diffusion: The dynamic viscosity coefficient appearing in the Navier-Stokes equations is often assumed to depend on temperature as $T^r$, for some $0<r<1$. Although $\nu$ has a very different physical interpretation it might be tempting to try $\nu\sim T^r/\rho$. However, some preliminary (unpublished) numerical investigations suggest that a temperature dependent diffusion is less accurate for \eqref{eulerian_ideal} than $\nu=\mu_0/\rho$ where $\mu_0=constant$. Hence, we will not consider models with $\nu\sim T^r/\rho$.

In a dense gas, the molecules will constantly experience the repelling forces of neighbouring molecules, which will randomly deviate their paths, which in turn is a diffusive process at the global scale. Moreover, the more densely packed the gas is, the more prominent will the internal friction be in relation to the convective process. Hence, it seems plausible that there is some $\rho^s$, $s>0$ dependence of $\nu$. We take $s=1$. Since  $\nu=\mu_0/\rho$, gives accurate results under normal conditions, we assume that the diffusion coefficient is given by,
\begin{align}
\nu = \frac{\mu_0}{\rho} + \mu_1\rho, \quad \mu_0>>\mu_1>0\label{nu}.
\end{align}

Finally, to control pressure in the full system, a sufficiently strong estimate of temperature is also needed. Such estimate is unattainable for the base model and we add a Fourier-type heat flux that models radiation. That is, $\kappa\nabla T= 4\kappa_rT^3 \nabla T$, where $\kappa_r$ is a material dependent constant. Under normal conditions (modest temperatures),  $4\kappa_rT^3$ should be orders of magnitude less than $\mu_0$ such that the radiation has a negligible effect.
\begin{remark}
The extra radation term does not increase the physical validity range. For that, more aspects of radiation physics would have to be included in the model. In fact, if the radiation diffusion is significant, the solution is already outside the validity range. In this sense, one can view the extra temperature diffusion as a technical modification. 
\end{remark}

In the remainder of this text, we consider the following alternative Navier-Stokes system:
\begin{subequations}\label{eulerian} \\
\begin{align}
\partial_t \rho + \nabla\cdot(\rho \vv )&= \nabla\cdot  (\nu \nabla \rho),\label{continuity} \\
\partial_t (\rho \vv) + \nabla\cdot(\rho \vv \otimes \vv) + \nabla p&=\nabla\cdot 
(\nu \nabla \rho \vv)), \qquad t\in [0,\T] \label{momentum} \\
\partial_t (E) + \nabla\cdot(E \vv +p\vv)  &= \nabla\cdot ( \nu \nabla E) + \nabla\cdot(\kappa_r \nabla T^4), \label{energy} \\ 
p&=\rho R T, \quad \textrm{ideal gas law,}\label{gaslaw}
\end{align}
\end{subequations}
with $\nu$ given by \eqref{nu}. Finally, we take $\Omega\subset \Rbb^3$ to be bounded, and its boundary, $\partial \Omega$, models an adiabatic wall, which is given by the no-slip  boundary condition
\begin{align}
  \vv=0,\label{no-slip}
\end{align}
along with
\begin{align}
\partial_n T=\partial_n\rho=0,\label{adiabatic_bc}
\end{align}
where $\partial_n=\nn\cdot\nabla$ and $\nn$ is the wall normal. (See \cite{Svard18} for a discussion on boundary conditions.)

\subsection{Outline}

\begin{itemize}
\item In Section \ref{sec:balance}, auxiliary balance laws are derived.
\item Using the equations and the auxiliary balance laws, we derive a priori estimates in Section \ref{sec:apriori_eqns}. As a consequence of the estimates, we are able to demonstrate positivity of $\rho$ and $T$ in Section \ref{sec:spec_vol}.
\item In Section \ref{sec:weak_seq}, we formally demonstrate that the a priori estimates supports a weak entropy solution.

\item In Section \ref{sec:weak_sol}, we  define weak entropy solutions. 
  
\item Next, we propose a numerical scheme in Section \ref{sec:scheme} and give the main result of this paper: Existence of a weak entropy solution to (\ref{eulerian}). The remainder of the paper proves this result.
  \item Some auxiliary (numerical) balance laws are derived in Section \ref{sec:add_laws}.
\item We derive a priori estimates for the numerical scheme in Section \ref{sec:apriori_scheme}.
\item In Section \ref{sec:ex}, we use the global estimates to infer existence of weak solutions to the full system by proving weak convergence of the numerical scheme.
  \item Section \ref{sec:final}, contains some concluding remarks.
  \end{itemize}

\subsection{Notation}

The following notation will be used throughout the article.
\begin{itemize}
\item The space $L^p(0,\T;L^q(\Omega))$ is denoted $L^p_q$ for short, when there is no risk of confusion.

\item Similarly, we write $S_1(S_2)=S_1(0,\T;S_2(\Omega))$ where $S_{1,2}$ are two generic function spaces.

\item $W^{r,p}(\Omega)$ is the Sobolev space where all derivatives of order $0,1,...,r$ are bounded in $L^p(\Omega)$. The space $W^{1,2}$ is denoted $H^1$. Furthermore, the norm of $W^{r,p}(\Omega)$ is denoted $\|\cdot\|_{r,p}$.

\item Throughout, we let $\C$ denote a generic a priori bounded constant obtained from the initial data. 
Likewise, $\epsilon$ and $\delta$ are used as small positive constants. Note that all three constants may change their actual values throughout a calculation.

  \item We denote the $i$th component of a vector ${\bf a}$ as ${\bf a}_i$.

\item In the remainder of the paper, we will use a number of standard results which have been collected in  Appendix \ref{app1} for the reader's convenience.
\end{itemize}

\section{Additional balance laws}\label{sec:balance}

\subsection{Entropy balance}

The first relation we derive is an entropy equation, which is an expression of the Second Law of Thermodynamics. The \emph{specific entropy} is $S=\log(p/\rho^{\gamma})$ and here we take the standard approach and derive an equation for the entropy function $U=-\rho S$. Associated with the \emph{entropy function} $U$ are the \emph{entropy fluxes} $\Fb_i=-m_iS$, $i=1,2,3$, and the \emph{entropy variables},
$$ 
U_\uu=\ww^T=(-(S-\gamma)-\frac{|\vv|^2}{2c_vT}, \frac{\vv}{c_v T}, -\frac{1}{c_vT}).
$$
Contracting \eqref{eulerian} with the entropy variables gives,
\begin{align}
  U_t + \nabla \cdot \Fb=\sum_{i=1}^5\ww_i\nabla\cdot (\nu\nabla \uu_i)+\ww_5\nabla \cdot(\kappa_r\nabla T^4),\nonumber
\end{align}
which can be recast as,
\begin{align}
  U_t + \nabla \cdot \Fb=\sum_{i=1}^5\left(\nabla\cdot (\ww_i \nu\nabla \uu_i)-(\nabla \ww_i^T)\nu (\nabla \uu_i)\right)+\nabla\cdot( \ww_5\kappa_r \nabla T^4)
-\nabla\ww_5\cdot( \kappa_r \nabla T^4),\nonumber
\end{align}
and finally,
\begin{align}
  (-\rho S)_t + \nabla \cdot (-{\mm} S)&=\nonumber\\
  \nabla\cdot \left(-(S-\gamma)\nu\nabla\rho+\frac{\nu}{c_vT}\nabla\left(\frac{\rho|\vv|^2}{2}\right)-\frac{1}{c_vT}\nu\nabla e -\frac{1}{c_vT}\kappa_r \nabla T^4\right)\label{ent_balance}\\
  -\nu(\gamma-1)\rho|\nabla\ln \rho|^2-\nu\rho|\nabla \ln T|^2
  -\frac{\nu\rho}{c_vT}\sum_{i=1}^3|\nabla v_i|^2 - 4\kappa_rT^3 \frac{|\nabla T|^2}{T^2}\nonumber ,
\end{align}
which is the entropy balance.

\subsection{Kinetic Energy balance}

To obtain the kinetic energy equation for the model (\ref{eulerian}), we multiply the momentum equation by $\vv$. For the velocity terms on the left-hand side, we obtain
\begin{align}
   T1=\vv\cdot (\partial_t (\rho \vv) + \nabla\cdot(\rho \vv \otimes \vv))&=\nonumber \\
   (\partial_t\rho+\nabla\cdot(\rho\vv))|\vv|^2+\rho(\partial_t\vv+\vv\cdot\nabla\vv)\cdot \vv&=\nonumber \\
    (\partial_t\rho+\nabla\cdot(\rho\vv))|\vv|^2+\rho\partial_t\frac{|\vv|^2}{2}+\rho\vv\cdot\nabla\frac{|\vv|^2}{2}&=\nonumber\\
    \partial_t\left(\frac{\rho|\vv|^2}{2}\right)+\nabla\cdot\left( \frac{\rho|\vv|^2}{2}\vv \right)+(\partial_t\rho + \nabla\cdot(\rho\vv))\frac{|\vv|^2}{2}\nonumber.
\end{align}
Using the continuity equation (\ref{continuity}) in the last line, we obtain
\begin{align}
  T1=\partial_t\left(\frac{\rho|\vv|^2}{2}\right)+\nabla\cdot\left( \frac{\rho|\vv|^2}{2}\vv+\nu\frac{|\vv|^2}{2}\nabla\rho\right)-\nu\nabla\rho\cdot \nabla \frac{|\vv|^2}{2}\label{lhs}.
\end{align}
Next, we turn to the remaining terms of the momentum equation:
\begin{align}
T2=\vv\cdot(- \nabla p+\nabla\cdot
(\nu \nabla \rho \vv)))&=\nonumber \\
- \nabla\cdot( p\vv)+p\nabla\cdot \vv+\vv\cdot (\nabla\cdot 
(\nu \nabla \rho \vv))&= \nonumber \\
- \nabla\cdot( p\vv)+p\nabla\cdot(\vv)+ \nabla\cdot 
\left(\nu\nabla\rho|\vv|^2+\nu\rho\nabla \frac{ |\vv|^2}{2}\right)-\nu\rho|\nabla \vv|^2-\nu\nabla\rho\nabla\frac{|\vv|^2}{2}&\nonumber.
\end{align}
Combining $T1$ and $T2$, yields the kinetic energy balance,
\begin{align}
    \partial_t\left(\frac{\rho|\vv|^2}{2}\right)+\nabla\cdot\left( \frac{\rho|\vv|^2}{2}\vv+\nu\frac{|\vv|^2}{2}\nabla\rho\right)&=\nonumber\\
    - \nabla\cdot( p\vv)+p\nabla\cdot(\vv)+ 
    \nabla\cdot 
\left(\nu(\nabla\rho)|\vv|^2+\nu\rho\nabla \frac{ |\vv|^2}{2}\right)-\nu\rho|\nabla \vv|^2,\nonumber
\end{align}
that simplifies to
\begin{align}
    \partial_t\left(\frac{\rho|\vv|^2}{2}\right)+\nabla\cdot\left( \frac{\rho|\vv|^2}{2}\vv\right)&=\nonumber\\
    - \nabla\cdot( p\vv)+p\nabla\cdot(\vv)+ 
    \nabla\cdot 
    \left(\nu\nabla\left(\rho\frac{|\vv|^2}{2}\right)\right)-\nu\rho|\nabla \vv|^2.\label{KE_balance}
\end{align}

\subsection{Specific volume balance}

Multiply the continuity equation \eqref{continuity} by $-\rho^{-2}$, to obtain
\begin{align}
\frac{1}{2}\partial_t \rho^{-1}+\nabla\cdot (\rho^{-1} \vv) + \nabla(\rho^{-2})\cdot \rho\vv  =-\nabla(\rho^{-2}\nabla (\nu\rho))-\frac{2}{\rho^3}\nabla \rho\cdot (\nu\nabla \rho), \nonumber
\end{align}
or,
\begin{align}
\frac{1}{2}\partial_t \rho^{-1}+\nabla\cdot (\rho^{-1} \vv) +2 \nabla(\rho^{-1})\cdot \vv  =-\nabla(\rho^{-2}\nabla (\nu\rho))-2 \nu\rho |\nabla \rho^{-1}|^2.\label{svolume_balance} 
\end{align}

\subsection{Renormalised internal energy balance}

Subtracting the kinetic energy balance (\ref{KE_balance}), from the total energy equation (\ref{energy}), gives the internal energy balance
\begin{align}
  \partial_t (c_v\rho T) + \nabla\cdot(c_v\rho T \vv) &=     - p\nabla\cdot(\vv)+ \nabla\cdot ( \nu \nabla (c_v\rho T)) +\nu\rho|\nabla \vv|^2+
  \nabla\cdot(\kappa_r\nabla T^4). \label{IE_eq}
\end{align}
We follow \cite{FeireislVasseur10} and introduce $H(T)=(1+T)^{1-\omega}$, $\omega>0.$ Then we renormalise the equation by multiplying \eqref{IE_eq} by $H'(T)$. We carry out the procedure step by step and begin with the left-hand side. First,
\begin{align}
  H'(T) (c_v\rho_t T+c_v\rho T_t)&=\nonumber \\
  c_v (H'(T) \rho_t T+\rho H_t) + \rho_t c_v H -\rho_t c_v H&=\nonumber\\
  c_v (\rho H)_t + c_v (H'(T) T -H)\rho_t. \nonumber
\end{align}
In the same way,
\begin{align}
  H'(T) (c_v\nabla \cdot (\rho\vv T ) )&=\nonumber
  c_v \nabla\cdot (\rho \vv H) + c_v (H'(T) T -H)\nabla\cdot (\rho \vv). \nonumber
\end{align}
Using the continuity equation, the left-hand side of the renormalised equation becomes,
\begin{align}
c_v(\rho H)_t + c_v \nabla \cdot (\rho \vv T)) + c_v(H'(T)T-T)\nabla\cdot (\nu\nabla \rho).\nonumber
\end{align}

Turning to the right-hand side, we begin with 
\begin{align}
  H'\nabla\cdot (\nu \nabla(c_v\rho T))=c_vH'\nabla\cdot (\nu T \nabla\rho +\nu \rho \nabla T))&= \nonumber \\
  c_vH'T\nabla\cdot (\nu\nabla\rho) +  c_vH'\nu\nabla\rho\cdot \nabla T
  +c_v\nabla\cdot(H' \nu \rho \nabla T))-c_v\rho\nu H''|\nabla T|^2&= \nonumber\\
    c_vH'T\nabla\cdot (\nu\nabla\rho) +  c_v\nu\nabla\rho\cdot \nabla H
    +c_v\nabla\cdot ( \nu \rho \nabla H))-c_v\rho\nu H''|\nabla T|^2&= \nonumber\\
    c_vH'T\nabla\cdot (\nu\nabla\rho) +  c_v\nabla\cdot(\nu H \nabla\rho)-c_v H \nabla\cdot(\nu \nabla\rho)
    +c_v\nabla\cdot ( \nu \rho \nabla H))-c_v\rho\nu H''|\nabla T|^2&= \nonumber\\
    c_v(H'T-H)\nabla\cdot (\nu\nabla\rho) +  c_v\nabla\cdot(\nu H \nabla\rho)
    +c_v\nabla\cdot ( \nu \rho \nabla H))-c_v\rho\nu H''|\nabla T|^2&. \nonumber 
\end{align}
Collecting the terms, we obtain the renormalised internal energy balance,
\begin{align}
  c_v(\rho H)_t + c_v \nabla \cdot (\rho \vv H)) &=\nonumber \\
      c_v\nabla\cdot(\nu H \nabla\rho)
    +c_v\nabla\cdot ( \nu \rho \nabla H))-c_v\rho\nu H''|\nabla T|^2& \label{RIE}\\ 
    +H'p \nabla \cdot \vv +H'\nu\rho |\nabla\vv|^2 + H'\nabla(\kappa_r \nabla T^3)\nonumber.
\end{align}

\section{A priori estimates for the full system}\label{sec:apriori_eqns}

We begin with the standard conservation properties. That is, we assume that $p,\rho>0$ and (separately) integrate \eqref{continuity} and \eqref{energy}. Using the boundary conditions (\ref{no-slip}) and (\ref{adiabatic_bc}) we obtain the bounds,
\begin{align}
  \rho\in L^{\infty}(0,\T;L^1(\Omega)),& \quad\quad
  E\in L^{\infty}(0,\T;L^1(\Omega)), \label{con_bounds}\\
    \sqrt{\rho} |\vv|\in L^{\infty}(0,\T;L^2(\Omega)),&\quad\quad
    p\in L^{\infty}(0,\T;L^1(\Omega)). \nonumber
\end{align}

\subsection{Entropy estimate}

Integrating (\ref{ent_balance}) in space and time,  using the boundary conditions (\ref{no-slip}) and (\ref{adiabatic_bc}),  leads to
\begin{align}
  \int_\Omega(-\rho S)|_{t=\T}-(-\rho S)|_{t=0} \,d\xb &=\label{ent_int}\\
  -\int_0^\T\int_{\Omega}\nu(\gamma-1)\rho|\nabla\ln \rho|^2+\nu\rho|\nabla \ln T|^2+
 \frac{\nu\rho}{c_vT}\sum_{i=1}^3|\nabla v_i|^2+  4\kappa_3T^3 \frac{|\nabla T|^2}{T^2}\, d\xb\,dt.\nonumber
\end{align}
The second term is bounded by initial data. Provided that we can control the first, we obtain the bounds
\begin{align}
  \int_0^\T\int_{\Omega}\nu(\gamma-1)\rho|\nabla\ln \rho|^2+\nu\rho|\nabla \ln T|^2+
 \frac{\nu\rho}{c_vT}\sum_{i=1}^3|\nabla v_i|^2+  4\kappa_rT |\nabla T|^2\, d\xb\,dt\leq \C.\label{ent_est}
\end{align}

Hence, we proceed to analyse the first integral of \eqref{ent_int}. (This was shown in \cite{FeireislVasseur10} and we repeat their arguments here.) We introduce $[z]^+=\max(z,0)$ and $[z]^-=\min(z,0)$ and carry out the following manipulations:
\begin{align}
  \int_\Omega(-\rho S) \,d\xb =
  \int_\Omega -\rho\log(RT)+(\gamma-1)\rho\log(\rho)\,d\xb=\nonumber \\
  \int_\Omega -\rho([\log(RT)]^++[\log(RT)]^-)+(\gamma-1)\rho([\log(\rho)]^++[\log(\rho)^-]\,d\xb\label{U_est}.
\end{align}
We need a bound on the negative terms. $\int_\Omega-\rho[\log(RT)]^+\,d\xb$ is bounded by $p\in L^\infty_1$ given by (\ref{con_bounds}). Since  $|\rho [\log(\rho)]^-|\leq 1$  and the domain is bounded, the last term is also bounded.  Hence, we have deduced that
\begin{align}
  \rho \log(\rho)&\in L^\infty(0,\T,L^1(\Omega)),  \label{rlogr} \\
  \rho \log(T)&\in L^\infty(0,\T,L^1(\Omega)),\nonumber 
\end{align}
and, with $\nu\sim \rho^{-1}+\rho$, we obtain from the diffusive terms \eqref{ent_est}, 
\begin{align}
  \nabla (\log(T))&\in L^2(0,\T;L^2)),\nonumber \\
  \rho\nabla (\log(T))&\in L^2(0,\T;L^2)),\nonumber \\
  \nabla (\log(\rho))&\in L^2(0,\T;L^2)),\label{ent_diff} \\
  \nabla \rho&\in L^2(0,\T;L^2)),\nonumber\\
  \frac{\rho}{\sqrt{T}}|\nabla \vv|&\in  L^2(0,\T;L^2)),\nonumber\\
  \frac{1}{\sqrt{T}}|\nabla \vv|&\in  L^2(0,\T;L^2)).\nonumber
\end{align}
From the radiation dissipation (the last term of (\ref{ent_est})), we obtain 
\begin{align}
  \nabla T^{3/2}&\in L^2(0,\T;L^2)).\label{rad_bound}
\end{align}

From the $\nabla\rho$ estimate in (\ref{ent_diff}), the $\rho\in L^{\infty}_1$ bound given in (\ref{con_bounds}) and Poincare's inequality, we also have,
\begin{align}
\rho \in L^2_2.\nonumber
\end{align}

From the current estimates, it is possible to obtain a better estimate on temperature. We begin by showing that $\rho>\delta>0$ on a subset of finite measure. (Again, we repeat the arguments laid out in \cite{FeireislVasseur10}.)

That is, we wish to show that
\begin{align}
|\B|=|\{x\in\Omega:\,| \rho(\xb,t)>\delta\}|>m>0.\label{B}
\end{align}
First, let $M_0$ be the total mass, which is constant thanks to the boundary conditions. Then consider the bound (\ref{rlogr}). Since $|\rho\log(\rho)|>\rho$ for $0<\rho<<1$ and $\rho>>1$, there must exist an $\alpha>0$ such that,
\begin{align}
\int_{\{\rho(\cdot,t)>\alpha\}}\rho(\cdot,t)\,d\xb\leq \frac{M_0}{3}.\nonumber
\end{align}
We write,
\begin{align}
  M_0=\int_\Omega\rho(\cdot,t)\,d\xb&=\nonumber \\
  \int_{\{\rho(\cdot,t)\leq \delta\}}\rho(\cdot,t)\,d\xb+
  \int_{\{\delta<\rho(\cdot,t)\leq\alpha\}}\rho(\cdot,t)\,d\xb+
  \int_{\{\rho(\cdot,t)>\alpha\}}\rho(\cdot,t)\,d\xb& \leq \nonumber \\
  \delta|\Omega|+\alpha m +\frac{M_0}{3}.\nonumber
\end{align}
We choose $\delta=M_0/(3|\Omega|)$ and deduce that $m\geq M_0/(3\alpha)>0$.

Next, we use $p\in L^{\infty}_1$ from (\ref{con_bounds}). Thanks to $\B$ having non-zero measure, we conclude that $T\in L^{\infty}(0,\T;L^1(\B))$. Hence, by a generalised Poincare inequality (Theorem \ref{theo:poincare}, Eqn. (\ref{poincare}) in Appendix \ref{app1}) and (\ref{rad_bound}), we have
\begin{align}
  T^{3/2}&\in L^2(H^1) \quad \textrm{and by Sobolev embedding},\nonumber \\
  T&\in L^3(L^9).\label{temp_bound}
\end{align}

\subsection{Kinetic energy estimate}

Integrating the balance \eqref{KE_balance},  using the boundary conditions (\ref{no-slip}) and (\ref{adiabatic_bc}), yields,
\begin{align}
   \int_\Omega \partial_t\left(\frac{\rho|\vv|^2}{2}\right)\,d\xb&=\int_{\Omega} \left( p\nabla\cdot(\vv) -\nu\rho|\nabla \vv|^2\right)\,d\xb\label{KE_est}.
\end{align}
We need to control $p \nabla\cdot(\vv)$, which is accomplished by
\begin{align}
\int_0^T \int_{\Omega} p\nabla\cdot(\vv)\,d\xb dt \leq 
\int_0^T \int_{\Omega} \eta^{-1}T^2 + \eta \rho^2|\nabla\cdot(\vv)|^2\,d\xb dt\label{pdivv}.
\end{align}
The first term on the right-hand side is controlled by (\ref{temp_bound}). Since $\nu \sim \mu_1\rho$, we can control the last term by choosing $\eta<\mu_1/2$ and from (\ref{KE_est}) we thus obtain the estimates,
\begin{align}
  \int_0^T\|\rho \nabla \vv\|^2_2\,dt \leq \C,\label{v_est} \\
  \int_0^T\| \nabla \vv\|^2_2\,dt \leq \C.\nonumber
\end{align}
Since $\sqrt{\rho} \vv\in L^2_2$ by (\ref{con_bounds}), and $\rho>\delta$ on $\B$, we have that $\vv\in L^2(0,\T;L^2(\B))$. Using the generalised Poincare inequality (\ref{poincare}), we obtain 
\begin{align}
  \vv &\in L^2(H^1), \label{v_H1} \quad \textrm{and by Sobolev embedding}\\
  \vv &\in L^2_6. \nonumber
\end{align}

\subsection{Estimates of the specific volume}\label{sec:spec_vol}

Integrating \eqref{svolume_balance} in space and time, using the boundary conditions (\ref{no-slip}) and (\ref{adiabatic_bc}),  gives,
\begin{align}
\int_{0,\Omega}^\T \frac{1}{2}\partial_t \rho^{-1}\,d\xb\, dt+2 \int_{0,\Omega}^\T\nabla(\rho^{-1})\cdot \vv \,d\xb\,dt  =-\int_{0,\Omega}^\T2 \nu\rho |\nabla \rho^{-1}|^2\,d\xb\,dt,
\end{align}
and using Young's inequality, gives,
\begin{align}
\int_{0,\Omega}^\T\frac{1}{2}\partial_t \rho^{-1} d\xb\,dt\leq 2 \int_{0,\Omega}^\T\eta|\nabla(\rho^{-1})|^2 +\eta^{-1}|\vv|^2d\xb\,dt  -2\int_{0,\Omega}^\T \nu\rho |\nabla \rho^{-1}|^2d\xb\,dt.\label{svolume_est} 
\end{align}
Using (\ref{v_H1}), choosing $\eta\leq \mu_0/2$ and noting that $\rho>0$, we obtain the bounds,
\begin{align}
  \sup_{t\in [0,\T]}\|\rho^{-1}\|_1\leq \C,\label{svolume_bounds} \\
  \int_0^T \|\nabla(\rho^{-1})\|^2_2\,dt\leq \C.\nonumber
\end{align}
By Poincare's inequality (\ref{poincare}) and Sobolev embedding, we also have
\begin{align}
  \rho^{-1}\in L^2(H^1), \label{svolume_b2}\\
  \rho^{-1}\in L^2(L^6). \nonumber
\end{align}
Consequently, we conclude that $\rho>0$ a.e.
Furthermore, from (\ref{rlogr}) we have a bound on $\rho\log(T)$. Hence, since $\rho>0$ a.e., we have that $T>0$ a.e.

\begin{remark}
Being able to establish the absence of large vacuum regions is a remarkable feature of (\ref{eulerian}) and it is a key property for the existence proof below. Note that it is a consequence of the mass diffusion.
\end{remark}

We can also obtain a useful estimate on temperature. By using the gas law, $p=\rho R T$, we have $T^{1/2}\sim p^{1/2}\rho^{-1/2}$. Since both $p^{1/2}$ and $\rho^{-1/2}$ are bounded in $L^{\infty}_2$ by \eqref{con_bounds} and \eqref{svolume_bounds}, we have
\begin{align}
\sqrt{T}\in L^{\infty}_1.\label{T_infty}
\end{align}

\subsection{Estimates from the continuity equation}

Multiplying \eqref{continuity} by $\rho$ and integrating in space, using conditions (\ref{no-slip}) and (\ref{adiabatic_bc}), yield
\begin{align}
    \frac{1}{2}d_t\|\rho\|_2^2-\int_{\Omega} \sqrt{\rho}\nabla \rho\cdot\sqrt{\rho} \vv \,d\xb =- \int_{\Omega} \nu(\nabla \rho)^T(\nabla \rho)\, d\xb.\nonumber
\end{align}
Using Young's inequality, we obtain,
\begin{align}
    \frac{1}{2}d_t\|\rho\|_2^2\leq -  \int_{\Omega} \left(\sqrt{\nu}|\nabla \rho|\right)^2\, d\xb + \int_{\Omega} \sqrt{\eta}(\sqrt{\rho}|\nabla \rho|)^2 +\eta^{-1/2}(\sqrt{\rho} |\vv|)^2 \, d\xb.\nonumber
\end{align}
Choose, $\eta\leq \mu_1/2$ and use $\sqrt{\rho}\vv\in L^{\infty}_2$ from \eqref{con_bounds} to obtain the bounds,
\begin{align}
  \sup_{t\in[0,\T]}\|\rho\|_2^2\leq \C.\label{rho_bound1}\\
\int_0^T\int_{\Omega} \left(\sqrt{\nu}|\nabla \rho|\right)^2\, d\xb<\C.\nonumber 
\end{align}
By \eqref{rho_bound1}, Poincare, Sobolev embedding, and a standard interpolation inequality (see Appendix \ref{app1}, eqn (\ref{interp1})), we have,
\begin{align}
  \rho\in L^{\infty}_2,\quad\quad
  \rho\in L^2(H^1),\quad\quad
  \rho\in L^2_6,\quad\quad
  \rho \in L^{10/3}_{10/3}.\label{rho_bound2} 
\end{align}

With (\ref{rho_bound2}) at hand, even stronger estimates on density are obtained by testing the continuity equation against $\rho^2$  ,
\begin{align}
  \frac{1}{3}d_t \|\rho\|_3^3 +\int_{\Omega} \rho^2(\nabla\cdot( \rho \vv))\,d\xb = \int_{\Omega} \rho^2 \,\nabla\cdot (\nu \, \nabla \rho)\,d\xb.\nonumber
\end{align}
Integrating by parts with the use  of the boundary conditions (\ref{no-slip}) and (\ref{adiabatic_bc}), results in
\begin{align}
\frac{1}{3}d_t \|\rho\|_3^3 -\int_{\Omega} \nabla \rho^2 \cdot (\rho \vv) \,d\xb = -\int_{\Omega} \nabla \rho^2 \cdot (\nu \, \nabla \rho)\,d\xb,\nonumber
\end{align}
and 
\begin{align}
\frac{1}{3}d_t \|\rho\|_3^3 -\int_{\Omega} (\frac{2}{3}(\nabla \rho^3)\cdot \vv)\,d\xb =  -\int_{\Omega}(2\rho \nabla \rho) \cdot (\nu \, \nabla \rho)\,d\xb.\nonumber
\end{align}
Partial integration and no-slip, lead to
\begin{align}
\frac{1}{3}d_t \|\rho\|_3^3 +\int_{\Omega} (\frac{2}{3} \rho^3 \nabla\cdot \vv )\,d\xb = -\int_{\Omega}(2\rho \nabla \rho) \cdot (\nu \, \nabla \rho)\,d\xb.\nonumber
\end{align}
Using Young's inequality, we split the term, $ \rho^3 \nabla\cdot \vv<\eta \rho^4 +\eta^{-1}(\rho \nabla\cdot\vv)^2 $, $\eta>0$ and use that the last term is bounded by \eqref{v_est}. We have
\begin{align}
\frac{1}{3}d_t \|\rho\|_3^3 \leq \int_{\Omega} \eta \frac{2}{3} \rho^4 d\xb -\int_{\Omega}(2\rho \nabla \rho) \cdot (\nu \, \nabla \rho)\,d\xb+\C,\nonumber
\end{align}
and using (\ref{nu}),
\begin{align}
  \frac{1}{3}d_t \|\rho\|_3^3 \leq \int_{\Omega} \eta \frac{2}{3} \rho^4 d\xb  -\int_{\Omega}(2\rho \nabla \rho) \cdot ((\mu_0\rho^{-1}+\mu_1\rho)) \, \nabla \rho)\,d\xb+\C.\nonumber 
\end{align}
Now we intend to use (\ref{rho_bound1}) and the generalised Poincare inequality (\ref{poincare}). Hence, we add and subtract the bounded term $\|\rho\|_2^2$ to the right-hand side. 
\begin{align}
  \frac{1}{3}d_t \|\rho\|_3^3 - \int_{\Omega} \eta \frac{2}{3} \rho^4 d\xb  +\int_{\Omega} (\mu_1 |\nabla \rho^2|^2) +2\mu_0 \, |\nabla \rho|^2\,d\xb+\|\rho\|_2^2-\|\rho\|_2^2\leq \C.\nonumber
\end{align}
Next, we use the generalised Poincare inequality to conclude that, $\epsilon_1\|\rho^2\|_2^2\leq \|\rho\|_2^2+\int_{\Omega}\frac{1}{2} (\mu_1 |\nabla \rho^2|^2)d\xb$ for some $\epsilon_1>0$, such that
\begin{align}
  \frac{1}{3}d_t \|\rho\|_3^3 - \int_{\Omega} \eta \frac{2}{3} \rho^4 d\xb  +\epsilon_1\|\rho^2\|_2^2+\int_{\Omega}\frac{1}{2} (\mu_1 |\nabla \rho^2|^2) +2\mu_0 \, |\nabla \rho|^2\,d\xb\leq \|\rho\|_2^2+\C. \nonumber
\end{align}

By choosing $0<\eta<\epsilon_1$, we obtain the bounds,
\begin{align}
  \rho\in L^{\infty}_3,\quad\quad
  \rho^2\in L^{2}(H^1),\quad\quad
  \rho\in L^4_{12},\quad\quad
  \rho\in L^{6}_{6}.\label{strong_rho} 
\end{align}
(The last bound by (\ref{interp4}).)

Given these bounds, we can get even better bounds. Multiply the continuity equation by $\rho^{3}$ to obtain,
\begin{align}
  \int_{0,\Omega}^\T\rho^{3}\rho_t\,d\xb\,dt+
  \int_{0,\Omega}^\T\rho^{3}(\nabla \cdot \rho \vv)\,d\xb\,dt=
  \int_{0,\Omega}^\T\rho^{3}\nabla\cdot (\nu \nabla  \rho)\,d\xb\,dt.\label{rhoTenThird}
\end{align}
We focus on the second integral. By repeatedly integrating by parts and using the boundary conditions (\ref{no-slip}) and (\ref{adiabatic_bc}), we have
\begin{align}
  \int_{0,\Omega}^\T\rho^{3}(\nabla \cdot \rho \vv)\,d\xb\,dt=
  -\int_{0,\Omega}^\T (\nabla\rho^{3})\rho \cdot \vv\,d\xb\,dt&= \nonumber \\
  -\int_{0,\Omega}^\T\frac{3}{4}(\nabla \rho^{4} \cdot \vv)\,d\xb\,dt=
  \int_{0,\Omega}^\T\frac{3}{4}( \rho^{3} \cdot \rho\nabla\vv)\,d\xb\,dt&\leq\nonumber\\
  \int_{0,\Omega}^\T(\frac{3}{4}\rho^{3})^2 + (\rho\nabla \cdot \vv)^2\,d\xb\,dt. \nonumber
\end{align}
These integrals are controlled by \eqref{v_est} and \eqref{strong_rho}. From \eqref{rhoTenThird}, we thus obtain the following estimates: 
\begin{align}
  \rho\in L^{\infty}_{4},\quad\quad
  \rho^{3/2}\nabla \rho \in L^2_2,\quad\quad
  \rho\in L^5_{15},\quad\quad
  \rho\in L^{23/3}_{23/3}.\label{very_strong_rho} 
\end{align}

\begin{remark}
It is possible to bootstrap further and obtain a little better integrability of density but the estimates above will serve the current purposes.
\end{remark}

\subsection{Renormalised internal energy}

By integrating (\ref{RIE}) in space and employing the boundary conditions, we obtain,
\begin{align}
 \int_\Omega c_v(\rho H)_t \,d\xb = -\int_\Omega c_v\rho\nu H''|\nabla T|^2\,d\xb
 +\int_\Omega \left(H'p \nabla \cdot \vv\ +H'\nu\rho |\nabla\vv|^2\right)\,d\xb\nonumber\\
 - \int_{\Omega} H''4\kappa_r T^3|\nabla T|^2\,d\xb,\nonumber
\end{align}
which we integrate in time to obtain
\begin{align}
- \int_{\Omega,0}^\T H''4\kappa_r T^3|\nabla T|^2\,d\xb\,dt
-\int_{\Omega,0}^\T c_v\rho\nu H''|\nabla T|^2\,d\xb\,dt&=\nonumber \\
- \int_\Omega c_v(\rho H)|_{t=\T}\,d\xb+ \int_\Omega c_v(\rho H)|_{t=0}
\,d\xb 
  +\int_{\Omega,0}^{\T} \left(H'p \nabla \cdot \vv\ +H'\nu\rho |\nabla \vv|^2\right)\,d\xb\,dt.\label{RIE_est}
\end{align}
The particular choice, $H(T)=(1+T)^{1-\omega}$, $0<\omega<1$, implies:
\begin{itemize}
  \item $H(T)< T+1$.
  \item $0<H'(T)=(1-\omega)(1+T)^{-\omega}<(1-\omega)$ for all $T\geq 0$.
  \item $H''(T)=-\omega(1-\omega)(1+T)^{-\omega-1}<0$ for all $T>0$.  
\end{itemize}
Hence, the left-hand side is positive. The second integral on the right-hand side of (\ref{RIE_est}) is bounded by initial data and the first by (\ref{con_bounds}). (That is, use the $L^\infty_1$ estimates of $p$ and $\rho$ and note that $\rho H \lesssim \rho+p$.) Finally, the third integral on the right-hand side can be bounded as in \eqref{pdivv} and by \eqref{v_est}. 

We conclude that the left-hand side is bounded and the strongest bound is obtained from the first integral,
\begin{align}
\int_{\Omega,0}^{\T} (1+T)^{-\omega-1}T^3|\nabla T|^2\,d\xb\,dt\leq \C.\nonumber
\end{align}
From this bound, we have, $T^{1-\epsilon}\nabla T\in L^2_2$, or $\nabla(T^{2-\epsilon}) \in L^2_2$. By \eqref{temp_bound}, the generalised Poincare inequality and Sobolev embedding, we have $T^{2-\epsilon}\in L^2_6$  for any $\epsilon>0$. Here, we recast the bound as,
\begin{align}
\sqrt{T}\in L^{8-\epsilon}_{24-\epsilon}\quad\textrm{ for any}\quad \epsilon>0.\label{strong_T}
\end{align}
Next, we need the following interpolation inequality (see Appendix \ref{app1}),
\begin{align}
\|u\|^r_r\leq \|u\|^{r-r\theta}_p\|u\|^{r\theta}_q,\nonumber 
\end{align}
where $1=(r-r\theta)/p+r\theta/q$. We intend to use \eqref{strong_T} and \eqref{T_infty}. Therefore, we choose $r\theta=8$, $q=24$ and $p=1$ and obtain $r=8+2/3$. We conclude that $\sqrt{T}\in L^{8+2/3-\epsilon}_{8+2/3-\epsilon}$. By choosing $\omega$ sufficiently small, we get $\delta=2/3-\epsilon>0$ such that,
\begin{align}
T^4\in L^{1+\delta}_{1+\delta} \quad \textrm{for} \quad 1/12>\delta>0.\label{T4_L1}
\end{align}

\subsection{Improved estimates}\label{sec:improved}

By \eqref{temp_bound}, we have $T\in L^3_3$, which we use in
\begin{align}
  \int_0^\T \|p\|_2^2 \, dt = R^2 \int_{0,\Omega}^\T (\rho T)^2\,d\xb\,dt \leq
\C \int_0^\T \|\rho^2\|_{3}\|T^2\|_{3/2} \, dt \leq \C\int_0^\T \|\rho^2\|_3^3+\|T^2\|_{3/2}^{3/2}\,dt.\nonumber
\end{align}
Since $L^6_6$ is embedded in $L^{23/3}_{23/3}$, the density term is bounded by \eqref{very_strong_rho}. Hence,
\begin{align}
p\in L^2_2.\label{p_L2}
\end{align}
However, using the full integrability of $\rho$ and $T(\in L^4_4)$ we have that
\begin{align}
p\in L^{2+\delta}_{2+\delta}.\label{p_L2plus}
\end{align}
for some $\delta>0$.

Next, we turn to velocity,
\begin{align}
\|\vv\|_1\leq \|\rho^{-1/2}\|_2^2 + \|\sqrt{\rho}|\vv|\|_2^2.\nonumber
\end{align}
The right-hand side is bounded by \eqref{svolume_bounds} and \eqref{con_bounds}, and we have,
\begin{align}
\vv\in L^\infty_1.\label{v_infty}
\end{align}
Next, we use Nash' inequality,
\begin{align}
\|u\|_2^{1+2/n}\leq \C \|u\|_1^{2/n}\|Du\|_2,\label{nash}
\end{align}
where $u$ is a function, $n$ is the number of spatial dimensions, and $D$ the differential operator. Applying \eqref{nash} with $n=3$ to velocity gives,
\begin{align}
\int_0^\T( \|\vv\|_2^{5/3})^2\,dt \leq \C\int_0^\T  \|\vv\|_1^{4/3}\|\nabla\vv\|^2_2\,dt\leq \C\sup_{t\in [0,\T]} \|\vv\|_1^{4/3}\int_0^\T\|\nabla\vv\|^2_2\,dt, \label{nash_vel}
\end{align}
which is bounded by \eqref{v_infty} and \eqref{v_H1}. 

Next, we apply the interpolation inequality $\|u\|_3\leq \|u\|_1^{1/5}\|u\|_6^{4/5}$ to obtain,
\begin{align}
\int_0^\T \|\vv\|_3^{10/4}\,dt\leq \int_0^\T\|\vv\|_1^{10/4\cdot1/5}\|\vv\|_6^2\,dt,\label{v_interp}
\end{align}
which is bounded by \eqref{v_infty} and \eqref{v_H1}. We summarise \eqref{v_interp} and \eqref{nash_vel},
\begin{align}
  \vv\in L^{10/3}_2,\quad\quad
  \vv\in L^{10/4}_3.\label{strong_vel}
\end{align}

Turning to momentum, we consider
$\|\rho^2 |\vv|^2\|_1 \leq \|\rho^2\|_{3/2}\|\vv\|_6^2$,
and by \eqref{v_H1} and \eqref{strong_rho}, we conclude that
\begin{align}
  \rho \vv \in L^2_2\label{mom_l2}.
\end{align}

Next, we will improve (\ref{mom_l2}) slightly:
\begin{align}
  (\rho^2u^2)^{1+\delta}=\rho^{2(1+\delta)}u^{2\epsilon(1+\delta)}u^{2(1-\epsilon)(1+\delta)}.\nonumber  
\end{align}
This yields
\begin{align}
\|(\rho^2u^2)^{1+\delta}\|_1\leq \|\rho^{2(1+\delta)}u^{2\epsilon(1+\delta)}\|_{3/2}\|u^{2(1-\epsilon)(1+\delta)}\|_3.\nonumber
\end{align}
To bound the last factor, we require that $2(1-\epsilon)(1+\delta)\leq 2$. This is satisfied for $\delta\leq \epsilon$. Set $r=2(1+\delta)$ and consider the remaining factor
\begin{align}
  \|\rho^{r}u^{r\epsilon}\|_{3/2}= \|\rho^{r-\epsilon r/2}\rho^{\epsilon r/2}u^{r\epsilon}\|_{3/2}= \|\rho^{r-\epsilon r/2}(\sqrt{\rho} u)^{\epsilon r}\|_{3/2}&\leq\nonumber \\
\|\rho^{\frac{3}{2}(r-\epsilon r/2)}\|_q\|(\sqrt{\rho} u)^{\frac{3}{2}\epsilon r}\|_{p}&.\nonumber
  \end{align}
Set $\frac{3}{2}\epsilon r p= 2$, i.e. $p=\frac{2}{3\epsilon(1+\delta)}$. It is clear that we can choose $\epsilon>0$ so small, while keeping $\delta<\epsilon$, such that  $\frac{3}{2}(r-\epsilon r/2)q\leq 4$. (We use \eqref{con_bounds},\eqref{very_strong_rho} and \eqref{v_H1} to control all factors.) Hence,
\begin{align}
\rho u\leq L^{2+\delta}_{2+\delta}, \quad \textrm{for some}\quad \delta>0.\label{mom_l2plus}
\end{align}

Furthermore, by $\|(\rho |\vv|)^{8/5}\|\leq \|\sqrt{\rho}^{8/5}\|_{5}\|(\sqrt{\rho}|\vv|)^{8/5}\|_{5/4}=\|\sqrt{\rho}\|^{8/5}_{8}\|(\sqrt{\rho}|\vv|)\|^{8/5}_{2}\leq \|\sqrt{\rho}\|^{16/5}_{8}+\|(\sqrt{\rho}|\vv|)\|^{16/5}_{2}$. Using (\ref{very_strong_rho}) and (\ref{con_bounds}), we have
\begin{align}
\rho \vv\in L^{\infty}_{8/5}.\label{mom_linf}
\end{align}
Moreover, using (\ref{mom_linf}) together with (\eqref{strong_vel}) and the generalised H{\"o}lder inequality imply,
\begin{align}
\|\rho (\vv\otimes \vv) \|_{24/23}\leq \C\|\rho \vv\|_{8/5}\|\vv\|_{3},\nonumber
\end{align}
for some bounded constant $\C$, such that
\begin{align}
\rho (\vv\otimes \vv) \in L^{10/4}_{24/23}.\label{mom_flux}
\end{align}

Next, we turn to the terms appearing in the energy flux. We will need integrability a little better than $L^1_1$ and therefore we take $\delta,\epsilon>0$ and consider
\begin{align}
  \|\rho|\vv| |\vv|^2\|_{1+\delta}^{1+\delta}=\|(\rho|\vv| |\vv|^2)^{1+\delta}\|&=\nonumber\\
\|((\sqrt{\rho})^{1-\epsilon}(\sqrt{\rho})^{1+\epsilon}|\vv|^{1+\epsilon} |\vv|^{2-\epsilon})^{1+\delta}\|&=\nonumber\\
\|(\sqrt{\rho})^{(1+\delta)(1-\epsilon)} (\sqrt{\rho}|\vv|)^{(1+\delta)(1+\epsilon)} |\vv|^{(2-\epsilon)(1+\delta)}\|&\leq \nonumber\\
\|(\sqrt{\rho})^{(1+\delta)(1-\epsilon)} (\sqrt{\rho}|\vv|)^{(1+\delta)(1+\epsilon)}\|_{3/2}\| |\vv|^{(2-\epsilon)(1+\delta)}\|_3&. \nonumber
\end{align}
To bound the last factor in $L^2_6$ using (\ref{v_H1}), we demand that $(2-\epsilon)(1+\delta)\leq 2$. We proceed with the choice $\epsilon=1/100$ and $\delta=1/200$.
\begin{align}
  \|\rho|\vv| |\vv|^2\|_{1+\delta}^{1+\delta}&\leq
\|(\sqrt{\rho})^{197/198} (\sqrt{\rho}|\vv|)^{607/598}\|_{3/2}\| |\vv|\|_r^s, \nonumber
\end{align}
where $0<r<6$ and $0<s<2$. Next, we estimate
\begin{align}
  \|(\sqrt{\rho})^{(197/198)\frac{3}{2}} (\sqrt{\rho}|\vv|)^{(607/598)\frac{3}{2}}\|_{1}& \leq\nonumber\\
    \|(\sqrt{\rho})^{(197/198)\frac{3}{2}}\|_p\| (\sqrt{\rho}|\vv|)^{(607/598)\frac{3}{2}}\|_{q}&, \label{enflux1}
\end{align}
and choose $\frac{607}{598}\frac{3}{2}q=2$, i.e., $q=2392/1821$ and $p=2392/571$. The second factor in (\ref{enflux1}) is controlled by \eqref{con_bounds}. The first is,
\begin{align}
\|(\sqrt{\rho})^{(197/198)\frac{3}{2}}\|_{2392/571}\leq \|(\sqrt{\rho})\|^{197/132}_{52973/8473}\leq \|\sqrt{\rho}\|_8^\infty= \|\rho\|^{\infty}_4,\nonumber
\end{align}
where the right-hand side is controlled by (\ref{very_strong_rho}).  Hence,
\begin{align}
\rho |\vv|^2\vv \in L^{1+1/200}_{1+1/200}. \label{KE_flux_est}
\end{align}

We also need the following standard result:
$\int_{\Omega,0}^\T |\log \rho|^p\,d\xb \,dt  \leq C_p\left(\|\rho\|_2^2+\|\rho^{-1}\|_2^2\right),$
where $C_p$ is a bounded constant for any $p<\infty$.  Hence,
\begin{align}
  \log\rho \in L^p_p,\label{logr_Lpp}
\end{align}
for any arbitrary $p$, since $\rho,\rho^{-1}\in L^2_2$ by  (\eqref{svolume_b2} and \eqref{rho_bound2})

\subsection{Compactness of $\rho,\vv,T$}\label{sec:compact}

Assuming that we have a sequence of approximate solutions satisfying the a priori estimates, we will demonstrate that the primary variables converge a.e. (Such a sequence will be generated by means of a numerical scheme below.) 

{\bf Density:} We intend to apply Aubin-Lions Lemma (see Lemma \ref{lemma:aubin} in Appendix \ref{app1}) to the continuity equation. From \eqref{rho_bound2}, we know that $\rho\in L^2(H^1)$, and $H^1$ is compactly embedded i $L^2$.

We need a that $\rho_t$ is bounded in a space $X$, such that $L^2$ is continuously embedded in $X$. We try the negative Sobolev space $W^{-3,2}$ which is the set of all distributions, $u$, such that $\int_\Omega u\phi\,d\xb$ is bounded for all $\phi\in W^{3,2}$.

\begin{remark}
  $L^2$ is continuously embedded in $W^{-3,2}$ if $\|u\|_{W^{-3,2}}\leq \C \|u\|_{L^2}$, $u\in L^2$, Here, $\|u\|_{W^{-3,2}}=\sup_{\phi\in W^{3,2}(\Omega),\|\phi\|_{3,2}=1}\frac{|<\phi,u>|}{\|\phi\|_{3,2}}$. Since $|<\phi,u>|$ is bounded for all $u\in L^2$, we conclude that $L^2$ functions, are  bounded in $W^{-3,2}$.
\end{remark}

We proceed with the test functions $\phi\in W^{3,2}$, for which $\nabla \phi\in L^{\infty}$, and use them to test if  $\rho_t\in L^2(W^{-3,2})$.
\begin{align}
|<\partial_t \phi,\rho>|=|-<\nabla \phi,\nu \nabla \rho >-<\rho \vv,\nabla \phi>|\leq \nonumber\\
\|\nu \nabla \rho\|_1 \| \nabla \phi\|_{\infty}+\|\rho \vv\|_1\|\nabla \phi\|_{\infty}.\nonumber
\end{align}
Furthermore, we obtain
\begin{align}
\int_0^\T \sup_{\phi\in W^{3,2}, \|\phi\|= 1} |<\partial_t \phi,\rho>\,dt|^2& \leq \nonumber \\
\int_0^\T\sup_{\phi\in W^{3,2}, \|\phi\|= 1} \left(\|\nu \nabla \rho\|_1 \| \nabla \phi\|_{\infty}+\|\rho\vv\|_1\|\nabla \phi\|_{\infty}\right)^2 dt&.  \nonumber
\end{align}
First, we note that $\|\nabla \phi\|_{\infty}\leq\C \|\phi\|_{3,2}(\leq \C)$. Furthermore, \eqref{mom_l2}, \eqref{very_strong_rho}/\eqref{strong_rho}  and \eqref{ent_diff} bounds the fluxes  on the right-hand side in $L^2_2$ and we can use Aubin-Lions Lemma, to conclude (by Sobolev embedding) that $\rho$ is compact in $L^{10/3-\epsilon}_{10/3-\epsilon}$. Hence, a subsequence converges a.e. and since $\rho\in L^{23/3}_{23/3}$ by \eqref{very_strong_rho}, we have \emph{strong convergence of $\rho$ in $L^r_r$ for $r<23/3$.}

\vspace{0.5cm}

{\bf Temperature:} We assume that we have two sequences of approximate solutions $\rho^n,T^n$ satisfying the a priori estimates obtained from the equation: $\rho^n\in L^{23/2}_{23/2}; \quad T^n\in L^{4+\epsilon}_{4+\epsilon}$ (see \eqref{very_strong_rho} and \eqref{T4_L1}).

We have strong convergence of $\rho_n$ in $L^p_p$, $p<23/3$ and (subsequential) weak convergence of $T_n$ in $L^4_4$. 

From this, we have weak convergence of $\rho^nT^n$ in $L^2_2$ since for any $\phi\in L^2_2$
\begin{align}
<\phi, \rho^nT^n-\rho T> = <\phi T^n,\rho^n-\rho>+< \phi \rho, T^n-T> . 
\end{align}
For the first term to vanish $\phi T^n $ must be bounded in $L^{23/20+\epsilon}_{23/20+\epsilon}$, which it is. For the other term, we need $\rho \phi \in L^{4/3}_{4/3}$, which also holds.

Next, consider convergence of the $L^2$-norm:
\begin{align}
\int_0^T \int_{\Omega}(\rho^nT^n)^2-(\rho T)^2\,dx\,dt=\int_0^T \int_{\Omega}((\rho^n)^2-\rho^2)(T^n)^2+\rho^2((T^n)^2-T^2)\,dx\,dt,\nonumber
\end{align}
which approaches zero thanks to the strongly and weakly converging sequences.

Next, we consider convergence in norm:
\begin{align}
  \int_0^\T\|\rho^nT^n-\rho T\|_2^2\,dt &=\nonumber\\
  \int_0^\T<\rho^nT^n-\rho T,\rho^nT^n-\rho T>  \,dt& =\nonumber\\
  \int_0^\T(\rho^n)^2(T^n)^2-\rho^2T^2-2(\rho^nT^n-\rho T)\rho T  \,dt&. \nonumber
\end{align}
Convergence of the norm and the weak convergence in $L^2_2$ ensures that the last integral vanishes. Hence, we have strong convergence of $\rho T$ in $L^2_2$ and a.e. convergence of a subsequence. Since $\rho$ is a.e. convergent, there is a subsequence of $T_n$ that converges a.e.

\vspace{0.5cm}

{\bf Velocity:} Here we use that the sequence $\sqrt{\rho^n}u^n$ is bounded in $L^2_2$ \eqref{con_bounds}; we have $u^n$ bounded and weakly convergent in $L^{2.5}_{2.5}$, \eqref{strong_vel}; and $\sqrt{\rho^n}$ is strongly convergent up to $L^p_p$, $p<46/3$, thanks to \eqref{very_strong_rho}. 

Weak convergence of $\sqrt{\rho^n}u^n \in L^2_2$ requires that for $\phi\in L^2_2$, we have $\sqrt{\rho^n}\phi \in L^{5/3}_{5/3}$ and $\phi u \in L^{46/43+\epsilon}_{46/43+\epsilon}$. Both conditions can be confirmed with the available estimates.
 
Next, we check convergence of the $L^2_2$ norm:
\begin{align}
  \int_0^\T \int_{\Omega} \rho^n(u^n)^2-\rho u^2\,dx\,dt=
\int_0^T \int_{\Omega}(\rho^n-\rho)(u^n)^2+\rho((u^n)^2-u^2)\,dx\,dt. 
\end{align}
Since $(u^n)^2\in L^{5/4}_{5/4}$, we need $\rho$ to be strongly convergent in at least $L^5_5$, which indeed we have.

Hence, one can in the same way as for temperature deduce strong convergence of $\sqrt{\rho^n} u^n$ in $L^2_2$ and draw an a.e. converging subsequence and get a.e. convergence of $u^n$.

\begin{remark}
It is also possible to prove that $\rho u$ is compact in $L^{3/2}_{3/2}$ via Aubin-Lions Lemma. The fluxes are bounded in $L^{1+\epsilon}_{1+\epsilon}$ which is enough to show that $\partial_t (\rho u) \in L^{1+\epsilon}(W^{-3/2})$. Furthermore, $\rho \vv$ and $\nabla\cdot (\rho \vv)$ are both in $L^{1+\epsilon}_{1+\epsilon}$ and $L^{3/2}$ is compactly embedded in $W^{1,1}$. Hence, we have a.e. convergence of $\rho \vv$ and therefore also of $\vv$.
\end{remark}

\vspace{0.5cm}

{\bf Magnitude of velocity gradients:}

We apply the analogous arguments as above to the sequence $\sqrt{\rho^n}\sqrt{|u^n_x|}$. We have strong convergence of $\sqrt{\rho^n}\in L^{p<46/3}_{p<46/3}$ and weak convergence of $\sqrt{|u^n_x|} \in L^4_4$. This is enough to get \\ $\lim_{n\rightarrow \infty}\int_0^\T\left(\|\sqrt{\rho^n}\sqrt{|u^n_x|}\|_2^2-\|\sqrt{\rho}\sqrt{|u_x|}\|_2^2\right)\,dt\rightarrow 0$.

Weak convergence of $\sqrt{\rho^n}\sqrt{|u^n_x|} \in L^2_2$ can be verified with the help of the estimates. Hence, we can infer strong convergence of $\sqrt{\rho^n}\sqrt{|u^n_x|}\in L^2_2$ and draw an a.e. converging subsequence. We obtain that $|u^n_x|$ converges a.e. and strongly up to $L^{2-\epsilon}_{2-\epsilon}$ for $\epsilon>0$.

\vspace{0.5cm}

{\bf Magnitude of temperature and density gradients:} The argument can be mimicked for the sequences $\sqrt{\rho^n}\sqrt{|T^n_x|}$ and $\sqrt{\rho^n}\sqrt{|\rho^n_x|}$ that satisfy equally strong estimates. 

\vspace{0.5cm}

{\bf Summary:} $\rho,\vv,T$ and $|\nabla \vv|, |\nabla T|, |\nabla \rho|$ are all converging a.e.

\subsection{Weak sequential convergence }\label{sec:weak_seq}

Assuming that we have a sequence of approximate solutions that satisfies the a priori estimates, we will demonstrate that the estimates are sufficient to infer convergence to  an $L^1_1$ weak solution. To this end, we need that all inviscid fluxes are bounded in $L^{1+\epsilon}_{1+\epsilon}$, $\epsilon>0$ and a.e. convergence of the principal variables. The convergence of the diffusive fluxes are discussed separately.

Throughout, we assume that initial data is appropriately bounded. For instance, all variables bounded in $L_{\infty}$ and $\rho,T>0$ and bounded away from zero, would do.

\vspace{0.5cm}

{\bf Continuity equation:} We begin by testing the continuity equation against $\phi\in C_0^{1}(\bar\Omega\times [0,\T))$,
\begin{align}
\int_0^\T\left(<\phi_n,\rho_t>+<\phi,\nabla (\rho_n \vv_n)>\right)\,dt=\int_0^\T<\phi,\nabla\cdot(\nu_n\nabla\rho_n)>\,dt.\nonumber
\end{align}
Integration by parts results in
\begin{align}
  -<\phi,\rho_n>|_{t=0}-\int_0^\T\left(<\partial_t\phi,\rho_n>-<\nabla\cdot \phi,\rho_n \vv_n>\right)\,dt&=\nonumber \\
  -\int_0^\T<\nabla\phi,(\nu_n\nabla\rho_n)>\,dt&,\nonumber
\end{align}
which we require to hold for any $\phi\in C_0^{1}(\bar \Omega\times [0,\T))$.

  The first integral is well-defined thanks to the strong estimates of $\rho_n$. In the second, the momentum is bounded in $L^2_2$ by \eqref{mom_l2plus} and both $\rho_n$ and $\vv_n$ converge a.e. Hence, we have strong convergence of the inviscid flux term in $L^2_2$. The last integral is controlled by \eqref{ent_diff} and \eqref{strong_rho}. Since $\nu_n$ is a function of $\rho_n$, it converges a.e.. The weak convergence in $L^2_2$ of the density gradient, allows us to conclude that the product converges weakly to the correct limit.

\vspace{0.5cm}

{\bf Momentum equations:} The weak form is given by,
\begin{align}
  -<\phi,\rho_n\vv_n>|_{t=0}-\int_0^\T\left(<\partial_t\phi,\rho_n\vv_n>-<\nabla\cdot \phi,(\rho_n (\vv_n\otimes \vv_n)+p_n)\right)\,dt&=\nonumber \\
  -\int_0^\T<\nabla\phi,(\nu_n\nabla(\rho_n \vv_n))\,dt,\nonumber
\end{align}
where $\phi\in C_0^{1}(\Omega\times [0,\T); \Rbb^3)$.

The first integral is controlled by \eqref{mom_l2}. The  second integral is controlled by \eqref{mom_flux} which, along with a.e. convergence of $\rho_n,\vv_n$, implies weak convergence in $L^{1+\epsilon}_{1+\epsilon}$. 
Convergence of the pressure term follows from \eqref{p_L2plus}.

The integrand in the last term contains,
\begin{align}
  \rho (\nabla \rho) \cdot \vv,\quad
  \rho^2 \nabla\cdot \vv,\quad
  \rho^{-1} (\nabla \rho) \cdot \vv,\quad
   \nabla\cdot \vv.\nonumber
\end{align}
In the first term, we have $\rho\nabla\rho\in L^2_2$ by \eqref{strong_rho}. Hence, a subsequence converges weakly in $L^2_2$. By a.e. convergence of $\vv$ and $\vv\in L^{10/4}_{3}$ by \eqref{strong_vel}, we have strong convergence of $\vv\in L^{p}_{p}$, $1\leq p<2.5$. Hence, the product converges weakly in $L^{1+\epsilon}_{1+\epsilon}$ for some $\epsilon>0$.

The term $\rho^2 \nabla\cdot \vv$ is handled in the same way. A subsequence of $\nabla\cdot \vv$ converges weakly by \eqref{v_H1}. By \eqref{very_strong_rho} and a.e. convergence, $\rho$ converges strongly in $\rho\in L^{23/3-\epsilon}_{23/3-\epsilon}$. Hence, we have weak convergence of the second term in $L^{1+\epsilon}_{1+\epsilon}$ for some $\epsilon>0$.

The third term is handled like the first noting that $\nabla \log \rho \in L^2_2$ by \eqref{ent_diff}. Weak convergence of the fourth term is immediate (for a subsequence) by \eqref{v_H1}.

\vspace{0.5cm}

{\bf The total energy equation:} Finally, we turn to the weak form of the energy equation:
\begin{align}
  -<\phi,E_n>|_{t=0}-\int_0^\T\left(<\partial_t\phi,E_n>-<\nabla \phi,(E_n \vv_n+p_n\vv_n)\right)\,dt&=\nonumber \\
  -\int_0^\T<\nabla\phi,(\nu_n\nabla E_n)\,dt -\int_0^\T<\nabla\phi,(\kappa_r \nabla T^4_n)\,dt&,\nonumber 
\end{align}
with test functions $\phi\in C_0^{2}(\Omega\times [0,\T))$.

The first integral is handled in the same way as the inviscid flux in the momentum equation. Weak convergence of $p_n\vv_n$ appearing in the second integral follows from \eqref{p_L2plus} and \eqref{strong_vel}. The $\rho |\vv|^2\vv$ term is bounded in $L^{1+\epsilon}_{1+\epsilon}$ according to \eqref{KE_flux_est}. Furthermore, it converges a.e. Hence, it converges weakly in $L^1_1$ to the correct limit.

The radiation term is recast as $<\nabla \phi,\nabla T^4>=-<\Delta \phi,T^4>$. Strong convergence (of a subsequence) is ensured by \eqref{T4_L1} and a.e. convergence.

The pressure term in the diffusive flux is handled as follows,
\begin{align}
\nu\nabla p \sim (\rho+\rho^{-1}) \nabla p\sim\rho T \nabla \rho +\rho^2\nabla T+T\nabla\log\rho+\nabla T.\nonumber
\end{align}
These are easily handled by using \eqref{ent_diff}, \eqref{rho_bound2}, \eqref{p_L2plus} and  \eqref{strong_rho}, along with a.e. convergence of $\rho$ and $T$, and noting that by combining \eqref{rad_bound} and \eqref{ent_diff}, we get $\nabla T\in L^2_2$.

Finally, we consider
\begin{align}
  <\nabla \phi,\nu_n\nabla (\frac{1}{2}\rho_n|\vv_n|^2)>=
  -<\nabla\cdot((\nabla \phi)\nu_n), (\frac{1}{2}\rho_n|\vv_n|^2)>&=\nonumber\\
  -<\Delta\phi,\nu_n\frac{1}{2}\rho_n|\vv_n|^2>-<\nabla\phi\cdot\nabla \nu_n, (\frac{1}{2}\rho_n|\vv_n|^2)>&.\nonumber
\end{align}
Dropping the subscript we have terms of the form:
\begin{align}
  |\vv|^2,\quad \rho^2|\vv|^2,\quad  \rho|\vv|^2\nabla(\rho + \rho^{-1})\label{diff_en}.
\end{align}
Strong convergence of the first follows directly from  $\vv \in L^{2.5}_{2.5}$ (see (\ref{strong_vel})) and a.e convergence. In the second, we have that $\rho \vv$ is converging strongly in $L^{2+\delta}_{2+\delta}$ (see \eqref{mom_l2plus}).

The fourth is,
\begin{align}
-\nabla(\rho^{-1})\rho|\vv|^2=(\nabla \log \rho)|\vv|^2= \nabla (\log\rho|\vv|^2)-2\log\rho(\vv\cdot\nabla \vv).\nonumber
\end{align}
In the first term of the last expression, we can move the derivatives onto the test function. The bound $\log\rho |\vv|^2$ follows since $\log\rho \in L^p_p$ for any arbitrary $p$ (see (\ref{logr_Lpp})) and $ |\vv|^2\in L^{5/4}_{5/4}$.

We are left to bound $\log\rho(\vv\cdot\nabla \vv)$. We have that $\nabla \vv$ is converging weakly in $L^2_2$. Hence, we need that $\log(\rho)\vv$ is converging strongly in $L^2_2$. Since both $\rho$ and $\vv$ are a.e. convergent, we only need a bound on the product in a little better space than $L^2_2$. This follows since we have $\vv\in L^{10/4}_{10/4}$ and $\log(\rho)\in L^{11}_{11}$. 

Finally, we must deduce convergence of the third term of (\ref{diff_en}):
\begin{align}
  (\nabla\rho)\rho|\vv|^2=\frac{1}{2}(\nabla\rho^2)|\vv|^2=
  \frac{1}{2}\nabla(\rho^2|\vv|^2) - \rho^2\vv \nabla\vv.\label{odd_weak_form}
  \end{align}
In the first term, we can move the derivative to the test function and are left with $\rho^2|\vv|^2$ which is bounded in $L^{1+\epsilon}_{1+\epsilon}$ and a.e. convergent. The last step is to prove convergence of
\begin{align}
\rho^2\vv \nabla\vv= (\rho\vv)(\rho \nabla\vv).\label{difficult}
\end{align}
We have strong convergence of $(\rho\vv)\in L^{2+\delta}_{2+\delta}$ for some fixed $\delta>0$. Hence, if we have weak convergence of $(\rho \nabla\vv)\in L^{2-\epsilon}_{2-\epsilon}$ for $\epsilon>0$, (\ref{difficult}) is converging weakly in $L^1_1$.

By (\ref{v_est}), we know that $(\rho^n \nabla\vv^n)\in L^2_2$ and hence weakly converging to some limit $(\rho\nabla\vv)^*$. Using the strong and weak convergence of $\rho^n$ and $\nabla \vv^n$, we get that  $\rho^n \nabla\vv^n\rightharpoonup \rho\nabla\vv\in L^1_1$. That is when we test it against an $L^{\infty}_{\infty}$ function. Since the sequence converges to $(\rho\vv)^*$ for all test functions in $L^2_2$ including those that are also bounded in $L^\infty_{\infty}$, we conclude that  $(\rho\nabla\vv)^*=\rho\nabla\vv$ by uniqueness of the weak limit.

\vspace{0.5cm}

{\bf The entropy inequality:} $\rho S$ is convergent if $\rho\log(T)$ and $\rho\log(\rho)$ are both convergent in $L^{1+\epsilon}_{1+\epsilon}$. This follows from $\log\rho$ and $\log T$ being bounded in $L^2(H^1)$, which together with the strong estimates of $\rho$ along with a.e. convergence, ensures weak convergence of $\rho S \in L^{1+\epsilon}_{1+\epsilon} $

Convergence of the entropy flux follows from the $\log\rho$ and $\log T$ bounds and (\ref{mom_l2plus}), i.e., $\rho \vv \in L^{2+\epsilon}_{2+\epsilon}$ and a.e. convergence of $\rho,T,\vv$. ($\rho \vv$ converges a.e. and hence strongly in $L^{2+\delta}_{2+\delta}$ for $\delta<\epsilon$.)

Finally, the right-hand side of the entropy equation \eqref{ent_balance} is summable and non-positive ensuring that the entropy inequality, $(-\rho S)_t +\nabla\cdot (-\rho \vv S)\leq 0$, is satisfied in a distributional sense in $L^{1}_{1}$.

However, it is challenging to verify these relations for the numerical scheme below, and in this treatise we will demonstrate that entropy is diffused globally. That is,
\begin{align}
\int_\Omega -\rho S\,d\xb \leq \int_\Omega -\rho S\,d\xb|_{t=0}.\nonumber
\end{align}


\subsection{Weak entropy solution}\label{sec:weak_sol}

The weak sequential compactness above, demonstrates that it should be possible to construct weak entropy solutions to (\ref{eulerian}). Herein, we will utilise a numerical scheme to produce a sequence of approximate solutions. We will prove that the sequence converges to a weak solution in the following sense:
\begin{definition}\label{def:weak}
  We say that the triple $\{\rho,\vv,T\}$ is a \emph{weak entropy solution} to (\ref{eulerian}) with boundary conditions \eqref{no-slip} and \eqref{adiabatic_bc}, if   $\rho,T>0$ a.e. in $\Omega\times [0,\T]$ and the triple satisfies:
  \begin{itemize}
  \item For all $\phi\in C^{1}_0([0,\T)\times \bar \Omega)$, the continuity equation is satisfied weakly,
    \begin{align}
  -<\phi,\rho>|_{t=0}-\int_0^\T\left(<\partial_t\phi,\rho>+<\nabla \phi,\rho \vv>\right)\,dt&=
  -\int_0^\T<\nabla\phi,(\nu\nabla\rho)>\,dt.\nonumber
    \end{align}
\item For all $\phi\in C^{1}_0([0,\T)\times  \Omega)$, the momentum equation is satisfied weakly,
    \begin{align}
  -<\phi,\rho\vv>|_{t=0}-\int_0^\T\left(<\partial_t\phi,\rho\vv>+<\nabla\phi,(\rho (\vv\otimes \vv)+p)\right)\,dt&=\nonumber \\
  -\int_0^\T<\nabla\phi,\nu\nabla\rho\vv)\,dt.\nonumber
    \end{align}
  \item For all $\phi\in C^{2}_0([0,\T)\times \Omega)$, the internal energy equation is satisfied weakly,  
    \begin{align}
      -<\phi,E>|_{t=0}-\int_0^\T\left(<\partial_t\phi,E>+<\nabla\phi,(E+p)\vv>\right)\,dt&=\nonumber\\
      -\D       +\int_0^\T<\Delta\phi,(\kappa_r  T^4)\,dt, \nonumber
    \end{align}
where the diffusive term is interpreted as:
\begin{align}
 \D= <\nabla\phi,\nu\nabla E>=  <\nabla \phi,\nu_n\nabla (\frac{1}{2}\rho_n|\vv_n|^2+p)>,\nonumber
\end{align}
  and
  \begin{align}
 <\nabla \phi,\nu_n\nabla (\frac{1}{2}\rho_n|\vv_n|^2)>&=\nonumber\\
 -<\Delta\phi,\nu_n\frac{1}{2}\rho_n|\vv_n|^2>-<\nabla \phi, (\nabla\nu_n)\frac{1}{2}\rho_n|\vv_n|^2>&,\nonumber
\end{align}
and
\begin{align}
  \mu_1<\nabla \phi, (\nabla\rho_n)\frac{1}{2}\rho_n|\vv_n|^2> =
  -\mu_1<\Delta \phi, \frac{1}{4}\rho_n^2|\vv_n|^2> 
  -\mu_1<\nabla \phi, \frac{1}{2}\rho_n^2\vv_n\cdot\nabla\vv>,\nonumber
\end{align}

\begin{align}
  \mu_0<\nabla \phi, (\nabla\rho_n^{-1})\frac{1}{2}\rho_n|\vv_n|^2> =
  \frac{\mu_0}{2}<\Delta \phi, \log\rho_n|\vv_n|^2> 
  +\mu_0<\nabla \phi, (\log \rho_n)\vv_n\cdot\nabla\vv>.\nonumber
\end{align}

\item
  The entropy inequality
\begin{align}
\int_\Omega -\rho S\,d\xb \leq \int_\Omega -\rho S\,d\xb|_{t=0},\nonumber
\end{align}
is satisfied weakly.

  \end{itemize}
 
\qed
\end{definition}

\section{The numerical approximation scheme}\label{sec:scheme}

To generate a sequence of approximate solutions to (\ref{eulerian}), we use a finite volume scheme and consider the domain  $\Omega=(0,1)^3$. We discretise $\Omega$ with $N+1$ grid points in each direction: $x_i$, $y_j$, $z_k$, $i,j,k=0...N$. The discrete domain is denoted $\Omega_h$. Furthermore, let $x_{i+1/2}=\frac{x_{i+1}+x_i}{2}$ and similarly for $y_{j+1/2},z_{k+1/2}$.

Each grid point is associated with a control volume. In the interior of the domain, the control volumes are centred around the grid points: $\Vv_{ijk}=(x_{i-1/2},x_{i+1/2}]\times(y_{j-1/2},y_{j+1/2}]\times(z_{k-1/2},z_{k+1/2}]$ and, $V_{ijk}=|\Vv_{ijk}|$ denotes  the measure of the control volume. At the boundaries, the control volumes are defined in the same way (towards other grid points) and are closed with the domain boundary. (We use the standard node-centred finite volume scheme on a Cartesian grid.)

  The variables are constant within each control volume such that they form \emph{piecewise-constant functions in space}. We adopt the same notation for the approximations that we used for the continuous equation. For instance, $p_{ijk}$ is the value of the piecewise constant pressure and $\vv_{ijk}=(u_{ijk},v_{ijk},w_{ijk})$ is the velocity vector in  $\Vv_{ijk}$. 

Furthermore, we need the control-volume-boundary areas
\begin{align}
  \Ss^x_{jk}&=(y_{j+1/2}-y_{j-1/2})(z_{k+1/2}-z_{k-1/2}),\nonumber \\
  \Ss^y_{ik}&=(x_{i+1/2}-x_{i-1/2})(z_{k+1/2}-z_{k-1/2}),\nonumber \\
  \Ss^z_{ij}&=(x_{i+1/2}-x_{i-1/2})(y_{j+1/2}-y_{j-1/2}),\nonumber 
\end{align}

\begin{figure}[h]
\begin{center}
\includegraphics[height=5cm]{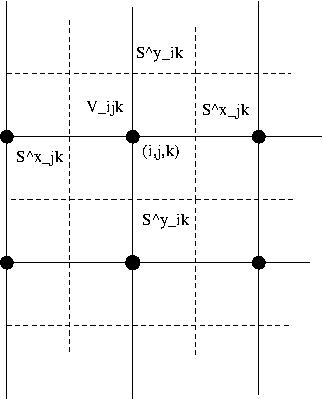}
\caption{A two-dimensional slice of a grid. $(i,j,k)$ denotes a point in the primal grid (solid lines) and the dashed lines indicate the dual grid (control volumes). $V_{ijk}$, $\Ss^x_{jk}, \Ss^y_{ik}$ are the control volume and control-volume-boundary areas associated with $(i,j,k)$.} 
\label{fig:cv}
\end{center}
\end{figure}

In Fig \ref{fig:cv}, a two-dimensional slice of control volumes is depicted (dashed lines). $\Ss^x,\Ss^y$ that make up the boundary of the control volume are also indicated for the point $(x_i,y_j,z_k)$.  

To simplify notation, we use the indices $-1/2$ and $N+1/2$ to denote the values at index $0$ and $N$, respectively. For instance, $x_{-1/2}=x_0$ and $p_{N+1/2,j,k}=p_{N,j,k}$.

We need the following notation for a scalar $a$:
\begin{align}
  \Delta_-^x a_{I+1jk}&=\Delta_+^x a_{Ijk}=a_{I+1jk}-a_{Ijk},\nonumber\\
  \Delta_-^y a_{iJ+1k}&=\Delta_+^y a_{iJk}=a_{iJ+1k}-a_{iJk},\nonumber\\
  \Delta_-^z a_{ijK+1}&=\Delta_+^z a_{ijK}=a_{ijK+1}-a_{ijK},\nonumber
\end{align}
where $I,J,K$ signify either a ``whole'' or ``half'' index and the boundary-index convention $-1/2\rightarrow 0$ etc, is applied. The difference operators are given as 
\begin{align}
  D^x_-a_{I+1jk}&=  \frac{\Delta^x_-a_{I+1jk}}{\Delta^x_- x_{I+1jk}},\nonumber \\
  D^y_-a_{iJ+1k}&=  \frac{\Delta^y_-a_{iJ+1k}}{\Delta^y_- y_{iJ+1k}},\label{diff_op} \\
  D^z_-a_{ijK+1}&=  \frac{\Delta^z_-a_{ijK+1}}{\Delta^z_- z_{ijK+1}},\nonumber
\end{align}
and similarly for $D_+^{x,y,z}$.

Furthermore, we need the two-point averages
\begin{align}
  \bar a_{i+1/2jk}&=\frac{a_{ijk}+a_{i+1jk}}{2},\nonumber\\
  \hat a_{i+1/2jk}& = \left\{\begin{array}{cc}\frac{\Delta^x_-  a_{i+1jk}}{\Delta_-^x\log a_{i+1jk}} & a_{ijk}\neq a_{i+1jk}\\
  a_{ijk} & a_{ijk}=a_{i+1jk}
  \end{array}
  \right.\quad\textrm{(log mean)},\label{means}\\
\check a_{i+1/2jk}&=\sqrt{a_{i+1jk}a_{ijk}}\quad\textrm{(geometric mean)},\nonumber
\end{align}
and for $a_i\geq 0$, the following elementary relations hold:
\begin{align}
\check a_{i+1/2}\leq \hat a_{i+1/2}\leq \bar a_{i+1/2}.\label{mean_rels} 
\end{align}

\begin{remark}
In (\ref{means}), the half index indicates in which direction the average is taken. See \cite{IsmailRoe09} for a stable numerical implementation of the log mean. 
\end{remark}

\begin{remark}
  The approximate derivative operators (\ref{diff_op}) are the same as the standard node-centred summation-by-parts (SBP) finite-volume scheme. (See \cite{NordstromForsberg03}.) In fact, since the grid is Cartesian, they coincide with the standard second-order SBP finite difference scheme. (See e.g. \cite{SvardNordstrom14}.)
\end{remark}

The following relation between arithmetic averages is frequently needed:
\begin{align}
  \overline {ab}_{i+1/2}&=\bar a_{i+1/2}\bar b_{i+1/2}+\frac{\Delta_+ a_i\Delta_+b_i}{4}.\label{split_av}
\end{align}
    The measure of the control volumes satisfies
\begin{align}
V_{ijk}=\Delta_-^xx_{i+1/2}\Ss^x_{jk}=\Delta_-^yy_{j+1/2}\Ss^y_{ik}=\Delta_-^zz_{k+1/2}\Ss^z_{ik},\label{V}
\end{align}
and we define the $L^2$-equivalent norm,
\begin{align}
\|a\|^2=\sum_{ijk=0}^NV_{ijk}a_{ijk}^2.\nonumber 
\end{align}
Throughout, we assume a uniform refinement of the grid such that $\Delta^x_- x_{i+1/2},\Delta^y_- y_{j+1/2},\Delta^z_- z_{k+1/2}$ are less than $h$ for all $i,j,k$. In this notation, all $\Ss\sim \mathcal{O}(h^2)$ and all $V\sim \mathcal{O}(h^3)$.

The discrete $L^2$-norms of the derivatives are given by
\begin{align}
  \|D^x_+a\|^2=\sum_{ijk=0}^{N-1,N,N}V_{ijk}(D_+^xa_{ijk})^2,\nonumber \\
  \|D^y_+a\|^2=\sum_{ijk=0}^{N,N-1,N}V_{ijk}(D_+^ya_{ijk})^2,\label{Dnorm} \\
  \|D^z_+a\|^2=\sum_{ijk=0}^{N,N,N-1}V_{ijk}(D_+^za_{ijk})^2.\nonumber 
\end{align}

\subsection{The numerical scheme}\label{sec:scheme}

The scheme is given by,
\begin{align}
  V_{ijk}(\uu_{ijk})_t+\Ss^x_{jk}\Delta^x_-\fb_{i+1/2jk}+\Ss^y_{ik}\Delta^y_-\gb_{ij+1/2k} +\Ss^z_{ij}\Delta^z_-\hb_{ijk+1/2}&=0\label{scheme}\\
  0\leq i,j,k,\leq N,\nonumber
\end{align}
where $\uu_{ijk}=(\rho_{ijk},\rho_{ijk}\vv_{ijk},E_{ijk})^T$. The fluxes consist of a convective and a diffusive part, $\fb_{i+1/2jk}=\fb^c_{i+1/2jk}-\fb^d_{i+1/2jk}$, and similarly for $\gb,\hb$. (These fluxes will be defined below.)

\begin{remark}
With the operators defined above, we can equivalently state the scheme on finite difference form as,
\begin{align}
  (\uu_{ijk})_t+D_-^z\fb_{i+1/2jk}+D_-^y\gb_{ij+1/2k}+D_-^z\hb_{ijk+1/2}&=0\label{fd_scheme}\\
0\leq i,j,k,\leq N\nonumber
\end{align}
\end{remark}

To define the fluxes, we use the notation:
\begin{align}
  \overline{|\vv|^2}_{i+1/2jk} &= (\overline{u^2})_{i+1/2jk}+(\overline{v^2})_{i+1/2jk}+(\overline{w^2})_{i+1/2jk},\nonumber\\
   \overline{|\vv|}^2_{i+1/2jk}&=(\bar u_{i+1/2jk})^2+(\bar v_{i+1/2jk})^2+(\bar w_{i+1/2jk})^2\quad \textrm{and analogously for y,z-means},\nonumber 
\end{align}

The convective fluxes (defined below) are  inspired by the entropy conservative and kinetic energy preserving flux proposed by Chandrashekar in \cite{Chandrashekar13} and we adopt his notation $\beta=\frac{1}{2RT}=\frac{\rho}{2p}$.

For $i=0...N-1$ and $j,k=0...N$:
\begin{align}
  p_{i+1/2jk}=&\frac{\bar \rho_{i+1/2jk}}{2\bar \beta_{i+1/2jk}}\nonumber\\
  \fb^{c,1}_{i+1/2jk}=&\overline {\rho u}_{i+1/2jk}\nonumber\\
  \fb^{c,2}_{i+1/2jk}=&
  \bar u_{i+1/2jk}\overline {\rho u}_{i+1/2jk}+p_{i+1/2jk}\nonumber\\
  \fb^{c,3}_{i+1/2jk}=&\bar v_{i+1/2jk}\overline {\rho u}_{i+1/2jk}\label{fc}\\
  \fb^{c,4}_{i+1/2jk}=&\bar w_{i+1/2jk}\overline {\rho u}_{i+1/2jk}\nonumber\\
  \fb^{c,5}_{i+1/2jk}=&\frac{1}{2(\gamma-1)\hat \beta_{i+1/2jk}}\overline {\rho u}_{i+1/2jk}-\frac{\overline{|\vv|^2}_{i+1/2jk}}{2}\overline {\rho u}_{i+1/2jk}\nonumber \\
  &+ \overline{|\vv|}^2_{i+1/2jk}\overline {\rho u}_{i+1/2jk}+p_{i+1/2jk}\bar u_{i+1/2jk} ,   \nonumber
\end{align}

\vspace{0.2cm}

For $j=0...N-1$ and $i,k=0...N$:
\begin{align}
p_{ij+1/2k}=&\frac{\bar \rho_{ij+1/2k}}{2\bar \beta_{ij+1/2k}}\nonumber\\
  \gb^{c,1}_{ij+1/2k}=&\overline{\rho v}_{ij+1/2k}\nonumber\\
  \gb^{c,2}_{ij+1/2k}=&\bar u_{ij+1/2k}\overline{\rho v}_{ij+1/2k}\nonumber\\
  \gb^{c,3}_{ij+1/2k}=
  &\bar v_{ij+1/2k}\overline{\rho v}_{ij+1/2k}+p_{ij+1/2k}\label{gc}\\
  \gb^{c,4}_{ij+1/2k}=&\bar w_{ij+1/2k}\overline{\rho v}_{ij+1/2k}\nonumber\\
  \gb^{c,5}_{ij+1/2k}=&\frac{1}{2(\gamma-1)\hat \beta_{ij+1/2k}}\overline{\rho v}_{ij+1/2k}-\frac{\overline{|\vv|^2}_{ij+1/2k}}{2}\overline{\rho v}_{ij+1/2k}\nonumber \\
&+\overline{|\vv|}^2_{ij+1/2k}\overline{\rho v}_{ij+1/2k}+p_{i+1/2jk}\bar v_{ij+1/2k}.
    \nonumber
\end{align}

\vspace{0.2cm}

For $k=0...N-1$ and $i,j=0...N$:
\begin{align}
  p_{ijk+1/2}=&\frac{\bar \rho_{ijk+1/2}}{2\bar \beta_{ijk+1/2}}\nonumber\\
  \hb^{c,1}_{ijk+1/2}=&\overline {\rho w}_{ijk+1/2}\nonumber\\
  \hb^{c,2}_{ijk+1/2}=&\bar u_{ijk+1/2}\overline {\rho w}_{ijk+1/2}\nonumber\\
  \hb^{c,3}_{ijk+1/2}=&\bar v_{ijk+1/2}\overline {\rho w}_{ijk+1/2}\label{hc}\\
  \hb^{c,4}_{ijk+1/2}=
  &\bar w_{ijk+1/2}\overline {\rho w}_{ijk+1/2}+p_{ijk+1/2}\nonumber\\
  \hb^{c,5}_{ijk+1/2}=&\frac{1}{2(\gamma-1)\hat \beta_{ijk+1/2}}\overline {\rho w}_{ijk+1/2}-\frac{\overline{|\vv|^2}_{ijk+1/2}}{2}\overline {\rho w}_{ijk+1/2}\nonumber \\
&+\overline{|\vv|}^2_{ijk+1/2}\overline {\rho w}_{ijk+1/2}+p_{ijk+1/2}\bar w_{ijk+1/2}.
    \nonumber
\end{align}
\begin{remark}
Above, we have used the convention that a numeral superscript signifies a component of a flux. 
\end{remark}

The diffusion coefficient is approximated as,
\begin{align}
  \nu_{i+1/2jk}&= \frac{\mu_0}{\hat \rho_{i+1/2jk}} + \mu_1\bar \rho_{i+1/2jk},\label{nu_approx}
\end{align}
and similarly in the y- and z-directions. 

Furthermore, we augment the diffusion with an with artificial (vanishing) component:
\begin{align}
  \lambda_{i+1/2jk}&=|\bar u_{i+1/2jk}|R^*_{i+1/2jk} + \frac{|\Delta^x_+ u_{ijk}|}{4},\label{ad_coeff2}\\
R^*_{i+1/2jk}&=\max(\frac{1}{2}, |\Delta^x_+ \log \rho_{ijk}|).\nonumber
\end{align}
We need the following Lemma:
\begin{lemma}\label{lemma:R}
\begin{align}
R^*_{i+1/2jk}\geq\nonumber\\ 
  \max\left(\frac{1}{2}, \left|\frac{\Delta^x_+ \rho_i}{12\bar \rho_{i+1/2}}\right|,\left|\frac{1}{2}\frac{(\sqrt{\rho_{i+1}}-\sqrt{\rho_i})}{\sqrt{\rho_{i+1}}+\sqrt{\rho_i}}\right|,
    \frac{\bar \rho_{i+1/2}}{2}\left(\frac{|\rho_{i+1}-\rho_i|}{\rho_{i+1}^2+\rho_{i+1}\rho_i+\rho_i^2}\right),\frac{|\rho_{i+1}-\rho_i|}{\hat \rho_{i+1/2}}\right)_{jk}\nonumber.
\end{align}
\end{lemma}
(As usual the obvious generalisations to the other two dimensions holds.)

\begin{proof}
We suppress the common $jk$ indices 
and check one term at the time. We use that $\rho_i>0$ (for all $i$) and the mean inequalities (\ref{mean_rels}). 

1. Trivial.

2.
\begin{align}
  \frac{|\Delta_+^x\rho_{i}|}{\bar\rho_{i+1/2}}\leq \frac{|\Delta_+^x\rho_i|}{\hat\rho_{i+1/2}}= |\Delta_+^x \log\rho_i|.\nonumber
\end{align}
3.
\begin{align}
  \frac{(\sqrt{\rho_{i+1}}-\sqrt{\rho_i})}{\sqrt{\rho_{i+1}}+\sqrt{\rho_i}}=
  \frac{(\rho_{i+1}-\rho_i)}{(\sqrt{\rho_{i+1}}+\sqrt{\rho_i})^2}=
  \frac{(\rho_{i+1}-\rho_i)}{\rho_{i+1}+\rho_i+2\sqrt{\rho_{i+1}}\sqrt{\rho_i}}&\leq\nonumber\\
  \frac{(\rho_{i+1}-\rho_i)}{\rho_{i+1}+\rho_i}&\leq \frac{1}{2}|\Delta_+^x\log\rho_i|\nonumber.
\end{align}
4.
  \begin{align}
  \frac{\bar \rho_{i+1/2}}{2}\left(\frac{|\rho_{i+1}-\rho_i|}{\rho_{i+1}^2+\rho_{i+1}\rho_i+\rho_i^2}\right)\leq\nonumber\\
  \frac{\bar \rho_{i+1/2}}{2}\left(\frac{|\rho_{i+1}-\rho_i|}{\frac{1}{2}(\rho_{i+1}+\rho_i)^2}\right)=\frac{\bar \rho_{i+1/2}}{2}\left(\frac{|\rho_{i+1}-\rho_i|}{2\bar\rho_{i+1/2}^2}\right)
  &\leq \frac{1}{4}|\Delta_+^x\log\rho_i |. \nonumber
\end{align}
5.
  \begin{align}
    \frac{|\rho_{i+1}-\rho_i|}{\hat \rho_{i+1/2}}&=|\Delta_+^x\log\rho_i|.\nonumber
\end{align}

\end{proof}

Next, we define
\begin{align}
  \tilde \nu_{i+1/2jk}&= \nu_{i+1/2jk}+(\Delta_+^x x_{i})\lambda_{i+1/2jk},\nonumber\\
  \tilde \nu_{ij+1/2k}&= \nu_{ij+1/2k}+(\Delta_+^y y_{j})\lambda_{ij+1/2k},\label{tilde_nu}\\
  \tilde \nu_{ijk+1/2}&= \nu_{ijk+1/2}+(\Delta_+^z z_{k})\lambda_{ijk+1/2}.\nonumber
\end{align}

The diffusive fluxes for $i,j,k=0...N-1$ are given by,
\begin{align}
  \fb^{d,1}_{i+1/2jk}&=\tilde \nu_{i+1/2jk}D^x_+\rho_{ijk},\nonumber \\
  \gb^{d,1}_{ij+1/2k}&=\tilde \nu_{ij+1/2k}D^y_+\rho_{ijk},\label{flux_d1}\\
  \hb^{d,1}_{ijk+1/2}&=\tilde \nu_{ijk+1/2}D^z_+\rho_{ijk},\nonumber 
\end{align}
and
\begin{align}
  \fb^{d,\{2,3,4\}}_{i+1/2jk}&=\tilde \nu_{i+1/2jk}D^x_+(\rho_{ijk}\{u_{ijk},v_{ijk},w_{ijk}\}),\nonumber \\
  \gb^{d,\{2,3,4\}}_{ij+1/2k}&=\tilde \nu_{ij+1/2k}D^y_+(\rho_{ijk}\{u_{ijk},v_{ijk},w_{ijk}\}),\label{flux_d} \\
  \hb^{d,\{2,3,4\}}_{ijk+1/2}&=\tilde \nu_{ijk+1/2}D^z_+(\rho_{ijk}\{u_{ijk},v_{ijk},w_{ijk}\}),\nonumber
\end{align}
and finally,
\begin{align}
  \fb^{d,5}_{i+1/2jk}=&\tilde \nu_{i+1/2jk}\left(\frac{(\pf_x)_{i+1/2jk}}{\gamma-1}+\frac{1}{2}D^x_+(\rho_{ijk}|\vv|^2_{ijk})+(\overline{|\vv|}^2_{i+1/2jk}-\overline {|\vv|^2}_{i+1/2jk})D^x_+\rho_{ijk}\right)\nonumber\\
  &+\kappa_rD^x_+T^4_{ijk}\nonumber\\
  \gb^{d,5}_{ij+1/2k}=&\tilde \nu_{ij+1/2k}\left(\frac{(\pf_y)_{ij+1/2k}}{\gamma-1}+\frac{1}{2}D^y_+(\rho_{ijk}|\vv|^2_{ijk})+(\overline{|\vv|}^2_{ij+1/2k}-\overline {|\vv|^2}_{ij+1/2k})D^y_+\rho_{ijk}\right)\nonumber \\
 & +\kappa_rD^y_+T^4_{ijk}\label{flux_d5}\\ 
  \hb^{d,5}_{ijk+1/2}=&\tilde \nu_{ijk+1/2}\left(\frac{(\pf_z)_{ijk+1/2}}{\gamma-1}+\frac{1}{2}D^z_+(\rho_{ijk}|\vv|^2_{ijk})+(\overline{|\vv|}_{ijk+1/2}^2-\overline {|\vv|^2}_{ijk+1/2})D^z_+\rho_{ijk}\right)\nonumber \\
 & +\kappa_rD^z_+T^4_{ijk}\nonumber
\end{align}
where
\begin{align}
  (\pf_x)_{i+1/2jk} &= \frac{1}{2\hat \beta_{i+1/2jk}}D_+^x\rho_{ijk}+\frac{\bar\rho_{i+1/2jk} }{2}D_+^x\frac{1}{\beta_{ijk}},\nonumber\\
  (\pf_y)_{ij+1/2k}&=\frac{1}{2\hat \beta_{ij+1/2k}}D_+^y\rho_{ijk}+\frac{\bar\rho_{ij+1/2k} }{2}D_+^y\frac{1}{\beta_{ijk}},\nonumber\\
  (\pf_z)_{ijk+1/2}&=\frac{1}{2\hat \beta_{ijk+1/2}}D_+^z\rho_{ijk}+\frac{\bar\rho_{ijk+1/2} }{2}D_+^z\frac{1}{\beta_{ijk}}.\nonumber
\end{align}

Furthermore, we  isolate the artificial diffusion by splitting the diffusive fluxes into
\begin{align}
  \fb^{d}_{i+1/2jk}&=\fb^{\nu}_{i+1/2jk}+\fb^{\lambda}_{i+1/2jk},\nonumber\\
  \gb^{d}_{ij+1/2k}&=\gb^{\nu}_{ij+1/2k}+\gb^{\lambda}_{ij+1/2k},\nonumber\\
  \hb^{d}_{ijk+1/2}&=\hb^{\nu}_{ijk+1/2}+\hb^{\lambda}_{ijk+1/2},\nonumber
\end{align}
where the $\lambda$-fluxes are obtained by replacing $\tilde \nu$ with $\lambda$, and $\kappa_r=0$, in (\ref{flux_d1}),(\ref{flux_d}) and (\ref{flux_d5}). (The $\nu$-fluxes are thus given by the $\fb^d$ with $\nu$ replacing $\tilde \nu$ and include the $\kappa_r$ part.)

{\bf Boundary conditions:} We invoke the no-slip condition strongly. That is, $u,v,w$ are all zero at the boundaries. Therefore,
\begin{align}
\uu_{0jk}=(\rho_{0jk},0,0,0,\frac{p_{0jk}}{\gamma-1})^T,\quad \uu_{Njk}=(\rho_{Njk},0,0,0,\frac{p_{Njk}}{\gamma-1})^T, \label{num_noslip}
\end{align}
which are used in the computation of $\fb_{1/2jk},\fb_{N-1/2jk}$. (Similarly in the other two directions.)

Since the velocities are zero at the boundary, there should be no scheme updating momentum at the boundary points. To to keep the simple structure of (\ref{scheme}), we enforce this by stipulating that every partial flux, $\fb,\fb^c,\fb^d,\fb^\lambda,...$  satisfy
\begin{align}
  \fb^{2,3,4}_{-1/2jk}=\fb^{2,3,4}_{1/2jk},\quad \gb^{2,3,4}_{0j+1/2k}=0,\quad \hb^{2,3,4}_{0jk+1/2}=0,\label{rm_momentum}
\end{align}
and similarly for all other boundaries. This choice amounts to $(\rho\vv)_t=0$ along the boundaries. Hence, if the momentum/velocity is initially zero, it will remain zero for all time.

Furthermore, invoking the no-slip in the remaining convective boundary flux components (for continuity and total energy), implies that they are all identically zero,
\begin{align}
  \fb^{c,1}_{-1/2jk}=\gb^{c,1}_{-1/2jk}=\hb^{c,1}_{-1/2jk}=0,&\label{no_slip_flux}\\
  \fb^{c,5}_{-1/2jk}=\gb^{c,5}_{-1/2jk}=\hb^{c,5}_{-1/2jk}=0,&\nonumber
\end{align}
and similarly at the other boundaries.

At the boundaries, the diffusive momentum fluxes are given by (\ref{rm_momentum}) which ensures that the momentum terms remain zero (in accordance with the no-slip condition).  Using no-slip and  homogeneous Neumann approximations for the wall normal gradient of density and temperature, the remaining two diffusive flux components are:
\begin{align}
  \fb^{d,1}_{-1/2jk}=\fb^{d,1}_{N+1/2jk}=0,&\quad \fb^{d,5}_{-1/2jk}=\fb^{d,5}_{N+1/2jk}=0,\nonumber \\
  \gb^{d,1}_{i,-1/2k}=\gb^{d,1}_{i,N+1/2k}=0,& \quad \gb^{d,5}_{i,-1/2k}=\gb^{d,5}_{i,N+1/2k}=0,\label{neumann_bc}\\
   \hb^{d,1}_{ij,-1/2}=\hb^{d,1}_{ij,N+1/2}=0,&   \quad  \hb^{d,5}_{ij,-1/2}=\hb^{d,5}_{ij,N+1/2}=0\nonumber
  \end{align}
Note that it is only the normal boundary flux that is affected by the Neumann boundary condition. The diffusive boundary fluxes (for the 1st and 5th components) tangent to the boundaries are not subject to any boundary conditions. 

{\bf Initial conditions:} We make the following assumptions on the initial data (t=0):
\begin{itemize}
\item All variables are uniformly bounded.
\item Density and temperature are positive and bounded away from zero.
\item At the boundary points, we set $\vv=0$ at $t=0$ (no-slip).
\item The initial data is projected onto the grid by e.g. computing the averages in each control volume. This also determines $\rho$ and $T$ at the boundary points. (We remark that injecting, instead of averaging, initial values to the grid points, does not change the overall formal accuracy of the method.)
\end{itemize}
It is possible to relax these assumptions but for practical purposes they are sufficient. \emph{Note that there is no smallness assumption on the initial data.}

The following is the main result of this paper:
\begin{theorem}\label{main_theo}
  Assume that at $t=0$ the initial data $\rho^0,\vv^0,T^0$ are bounded in $L^{\infty}(\Omega)$ and $\rho^0,T^0\geq constant>0$. Let $\{\rho_h,\vv_h,T_h\}$ be the family of numerical solutions generated from the initial data by the  scheme (\ref{scheme}) (see Section \ref{sec:scheme} for all approximations) as $h\rightarrow 0$. Then a subsequence converges strongly in the spaces,
  \begin{align}
    \rho_h\rightarrow \rho& \quad \in L^{23-\epsilon}_{23-\epsilon}\nonumber\\
    \vv_h\rightarrow \vv& \quad \in L^{10/4-\epsilon}_{10/4-\epsilon},\quad\quad (\epsilon>0),\nonumber \\
    T_h\rightarrow T &\quad \in L^{4-\epsilon}_{4-\epsilon},\nonumber 
  \end{align}
and $\{\rho, \vv, T\}$ is a weak entropy solution of the problem (\ref{eulerian}) in the sense of Def. \ref{def:weak}.
 
\end{theorem}

We summarise some of the estimates that resulting solution satisfies:
\begin{align}
      \rho \in L^{\infty}(L^4(\Omega))\cap L^{23/3}_{23/3}\cap L^5_{15},&\quad
      \rho^{5/2}\in L^2(H^1), \nonumber \\
      \rho^{-1}\in L^{\infty}_1\cap L^2(H^1)\cap L^2(L^6),&\quad
      \nabla(\log(\rho))\in L^2_2,\nonumber \\
      \vv \in L^\infty_1\cap L^2(H^1)\cap L^2(L^6)\cap L^{10/4}_3,&\quad
      \rho \nabla\vv\in L^2_2, \nonumber\\
      T\in L^{4+\delta}_{4+\delta},\quad \textrm{for some}\quad \delta>0,\nonumber\\
      T^{2-\epsilon}\in L^2(H^1)\cap L^2(L^6),\quad \textrm{for}\quad \epsilon>0,\nonumber\\
      \sqrt{T}\in L^{\infty}_1,\quad
      p\in L^{\infty}_2,&\quad
      \rho \log(T)\in L^{\infty}_1.\nonumber
  \end{align}

The remainder of the paper contains a proof of the theorem.

\section{Additional balance laws}\label{sec:add_laws}

Before we derive the discrete versions of the kinetic and internal energy balances, we state some auxiliary relations. For positive entities the following relations (exemplified with the positive density) hold:
\begin{align}
  \frac{1}{\check \rho_{i+1/2jk}} \geq  \frac{1}{\hat \rho_{i+1/2jk}}\geq \frac{1}{\bar \rho_{i+1/2jk}}.\label{rhoinf_rels}
\end{align}
They follow directly from (\ref{mean_rels}).

We shall also frequently use (\ref{split_av}) to obtain the following alternative form of the first component of the internal inviscid fluxes and artificial diffusion (here exemplified with the x-flux).
\begin{align}
  \fb^{c,1}_{i+1/2jk}-\fb^{\lambda,1}_{i+1/2jk}=&\overline {\rho u}_{i+1/2jk}-\lambda_{i+1/2jk}\Delta_+^x\rho_{ijk} \nonumber\\
 =&\bar\rho_{i+1/2jk} \bar u_{i+1/2jk}-\lambda^a_{i+1/2jk}\Delta_+^x\rho_{ijk},\label{alt_flux}
\end{align}
where
\begin{align}
  \lambda^a_{i+1/2jk}&=\lambda_{i+1/2jk}-\frac{\Delta^x_+ u_{ijk}}{4}>0,\nonumber
\end{align}
and
\begin{align}
  \fb^{c,1}_{i+1/2jk}-\fb^{\lambda,1}_{i+1/2jk}=&\overline {\rho}_{i+1/2jk} \bar u_{i+1/2jk}-\lambda^a_{i+1/2jk}\Delta_+^x\rho_{ijk} \nonumber\\
  =&\check\rho_{i+1/2jk} \bar u_{i+1/2jk}-\lambda^c_{i+1/2jk}\Delta_+^x\rho_{ijk},\label{alt_flux2}
\end{align}
where
\begin{align}
  \lambda^c_{i+1/2jk}&=\lambda_{i+1/2jk}-\frac{\Delta^x_+ u_{ijk}}{4}+\bar u_{i+1/2jk}\frac{\bar \rho_{i+1/2}-\check\rho_{i+1/2jk}}{\Delta_+^x\rho_{ijk} }>0.\nonumber
\end{align}
The last inequality follows from
\begin{align}
\frac{\bar \rho_{i+1/2}-\check\rho_{i+1/2jk}}{\Delta_+^x\rho_{ijk} }=
\frac{\frac{\rho_{i+1}+\rho_i}{2}-\sqrt{\rho_{i+1}\rho_i}}{\Delta_+^x\rho_{ijk} }&=\nonumber\\
\frac{1}{2}\frac{(\sqrt{\rho_{i+1}}-\sqrt{\rho_i})^2}{(\sqrt{\rho_{i+1}-\sqrt{\rho_i}})(\sqrt{\rho_{i+1}+\sqrt{\rho_i}})}=
\frac{1}{2}\frac{(\sqrt{\rho_{i+1}}-\sqrt{\rho_i})}{\sqrt{\rho_{i+1}+\sqrt{\rho_i}}}\label{ent_ad_coeff}
\end{align}
and Lemma \ref{lemma:R}.

\subsection{Kinetic and internal energy balances}

We derive the kinetic energy balance that will subsequently be subtracted from the total energy equation to produce the internal energy balance.

\begin{align}
\frac{d}{dt}K_{ijk}=-\frac{1}{2}|\vv|^2_{ijk}\frac{d\rho_j}{dt} + u_{ijk}\frac{d(\rho u)_{ijk}}{dt} + v_{ijk}\frac{d(\rho v)_{ijk}}{dt}  + w_{ijk}\frac{d(\rho w)_{ijk}}{dt}\label{KE0}.
\end{align}
We need the Leibniz rule
\begin{align}
  a_iD_-f_{i+1/2}=D_-(\bar a_{i+1/2}f_{i+1/2})-\frac{1}{2}f_{i+1/2}D_-a_{i+1}-\frac{1}{2}f_{i-1/2}D_-a_{i}. \label{Leibniz}
\end{align}
At the boundaries, we have the corresponding rules 
\begin{align}
a_0D_-f_{1/2}=D_-(\bar a_{1/2}f_{1/2})-\frac{1}{2}f_{1/2}D_-a_{1}\nonumber
\end{align}
where, in accordance with our convention, $\bar a_{-1/2}=a_0$. By adopting the convention that $D_-a_0\equiv 0$, (\ref{Leibniz}) applies to the boundaries as well, and the derivations below need not be repeated for the boundary points.

First we rewrite,
\begin{align}
  u_{ijk}\frac{d(\rho u)_{ijk}}{dt}&= \nonumber \\
  -u_{ijk}\left(D_-^x\fb^2_{i+1/2jk}+ D_-^y\gb^2_{ij+1/2k} +D_-^z\hb^2_{ijk+1/2}\right)&= \nonumber \\
  -\left(D_-^x(\bar u_{i+1/2jk}\fb^2_{i+1/2jk}) -\frac{1}{2}(D_-^xu_{i+1jk})\fb^2_{i+1/2jk}-\frac{1}{2}(D_-^xu_{ijk})\fb^2_{i-1/2jk}  \right)\nonumber \\
 - \left(D_-^y(\bar u_{ij+1/2k}\gb^2_{ij+1/2k}) -\frac{1}{2}(D_-^yu_{ij+1k})\gb^2_{ij+1/2k}-\frac{1}{2}(D_-^yu_{ijk})\gb^2_{ij-1/2k}  \right)\label{KE1} \\
 - \left(D_-^z(\bar u_{ijk+1/2}\hb^2_{ijk+1/2}) -\frac{1}{2}(D_-^zu_{ijk+1})\hb^2_{ijk+1/2}-\frac{1}{2}(D_-^zu_{ijk})\hb^2_{ijk-1/2}  \right),\nonumber 
\end{align}
and similarly for the momentum in the y- and z-directions.

Turning to the contribution from the continuity equation, and again using (\ref{Leibniz}), we have
\begin{align}
-\frac{1}{2}|\vv|^2_{ijk}\frac{d\rho_j}{dt}&=\nonumber\\
\frac{1}{2}|\vv|^2_{ijk}\left( D_-^x\fb^1_{i+1/2jk}+ D_-^y\gb^1_{ij+1/2k} +D_-^z\hb^1_{ijk+1/2}  \right)&=\nonumber\\
\frac{1}{2}\left( D_-^x(\overline{|\vv|^2}_{i+1/2jk}\fb^1_{i+1/2jk})-\frac{1}{2}(D_-^x|\vv|^2_{i+1jk})\fb^1_{i+1/2jk}-\frac{1}{2}(D_-^x|\vv|^2_{ijk})\fb^1_{i-1/2jk}\right)\nonumber \\
+\frac{1}{2}\left( D_-^y(\overline{|\vv|^2}_{ij+1/2k}\gb^1_{ij+1/2k})-\frac{1}{2}(D_-^y|\vv|^2_{ij+1/2k})\gb^1_{ij+1/2k}-\frac{1}{2}(D_-^y|\vv|^2_{ijk})\gb^1_{ij-1/2k}\right)\label{KE2} \\
+\frac{1}{2}\left( D_-^z(\overline{|\vv|^2}_{ijk+1/2}\hb^1_{ijk+1/2})-\frac{1}{2}(D_-^z|\vv|^2_{ijk+1/2})\hb^1_{ijk+1/2}-\frac{1}{2}(D_-^z|\vv|^2_{ijk})\hb^1_{ijk-1/2}\right).\nonumber
\end{align}
Using the identity,
\begin{align}
  D_-^x|\vv|^2_{ijk} = &D_-^x (u_{ijk}^2+v^2_{ijk}+w_{ijk}^2)= \nonumber \\
  &2\bar u_{i-1/2}D_-^xu_{ijk}
+2\bar v_{i-1/2}D_-^xv_{ijk}+2\bar w_{i-1/2}D_-^xw_{ijk},\nonumber
\end{align}
and (\ref{KE1})-(\ref{KE2}) in (\ref{KE0}), we obtain an expression for the kinetic energy balance:
\begin{align}
  \frac{dK_{ijk}}{dt}+D_-^x\Ff_{i+1/2jk}+D_-^y\Gf_{ij+1/2k}+D_-^z\Hf_{ijk+1/2}
  -(pD\vv)_{ijk}&=-\Df_{ijk},\label{KE_disc}\\
  0\leq i,j,k\leq N,\nonumber
\end{align}
where
\begin{align}
  \Ff_{i+1/2jk}&=-\frac{1}{2}\overline{|\vv|^2}_{i+1/2jk}\fb^1_{i+1/2jk}+\bar u_{i+1/2jk}\fb^2_{i+1/2jk}+\bar v_{i+1/2jk}\fb^3_{i+1/2jk}+\bar w_{i+1/2jk}\fb^4_{i+1/2jk} \nonumber \\
  \Gf_{ij+1/2k}&=-\frac{1}{2}\overline{|\vv|^2}_{ij+1/2k}\gb^1_{ij+1/2k}+\bar u_{ij+1/2k}\gb^2_{ij+1/2k}+\bar v_{ij+1/2k}\gb^3_{ij+1/2k}+\bar w_{ij+1/2k}\gb^4_{ij+1/2k}\label{KE_flux} \\
\Hf_{ijk+1/2}&=-\frac{1}{2}\overline{|\vv|^2}_{ijk+1/2}\hb^1_{ijk+1/2}+\bar u_{ijk+1/2}\hb^2_{ijk+1/2}+\bar v_{ijk+1/2}\hb^3_{ijk+1/2}+\bar w_{ijk+1/2}\hb^4_{ijk+1/2},\nonumber
\end{align}
and
\begin{align}
  (p D\vv)_{ijk}=&
    \quad\frac{1}{2}\left( (D_-^xu_{i+1jk})p_{i+1/2jk}+(D_-^xu_{ijk})p_{i-1/2jk}  \right)\nonumber \\
 & + \frac{1}{2}\left( (D_-^yv_{ij+1k})p_{ij+1/2k}+(D_-^yv_{ijk})p_{ij-1/2k}  \right)\nonumber \\
 & + \frac{1}{2}\left( (D_-^zw_{ijk+1})p_{ijk+1/2}+(D_-^zw_{ijk})p_{ijk-1/2}  \right).\nonumber
\end{align}
Using the boundary convention introduced in (\ref{Leibniz}) and the no-slip condition (implying that velocity differences along boundaries are zero), we have,  $ \Ff_{-1/2jk}\equiv 0$, etc. and
\begin{align}
  (p D\vv)_{0jk}=&\quad
    \frac{1}{2}\left( (D_-^xu_{1jk})p_{i+1/2jk}  \right)\nonumber 
\end{align}
and similarly for the other boundaries. The diffusive terms are,
\begin{align}
  \Df_{ijk}&=
  \frac{1}{2}\left(\tilde \nu_{i+1/2jk}\bar \rho_{i+1/2jk}(|D_-^x\vv_{i+1jk}|)^2+\tilde\nu_{i-1/2jk}\bar\rho_{i-1/2jk}(|D_-^x\vv_{ijk}|)^2  \right)\nonumber \\
& +\frac{1}{2} \left(\tilde\nu_{ij+1/2k}\bar \rho_{ij+1/2k}(|D_-^y\vv_{ij+1k}|)^2+\tilde\nu_{ij-1/2k}\bar \rho_{ij-1/2k}(|D_-^y\vv_{ijk}|)^2  \right)\nonumber \\
& +\frac{1}{2}\left(\tilde\nu_{ijk+1/2}\bar \rho_{ijk+1/2}(|D_-^z\vv_{ijk+1}|)^2+\tilde\nu_{ijk-1/2}\bar \rho_{ijk-1/2}(|D_-^z\vv_{ijk|})^2  \right),\nonumber 
\end{align}
where,
\begin{align}
  (|D_-^x\vv_{ijk}|)^2=(|D_-^xu_{ijk}|)^2+(|D_-^xw_{ijk}|)^2+(|D_-^xw_{ijk}|)^2, etc.\nonumber
\end{align}
At the boundary, most of the terms of $\Df_{ijk}$ are zero due to the no-slip condition and the boundary convention in (\ref{Leibniz}). For instance,
\begin{align}
  \Df_{0jk}=
  \frac{1}{2}\tilde \nu_{1/2jk}\bar \rho_{1/2jk}(|D_-^x\vv_{1jk}|)^2.\nonumber 
\end{align}

Subtracting the kinetic energy equation (\ref{KE_disc}) from the total energy equation of (\ref{fd_scheme}) gives the internal energy equation:
\begin{align}
  D_t\frac{p_{ijk}}{\gamma-1}+(pD\vv)_{ijk}-\Df_{ijk}
\nonumber\\
  +D_-^x(\frac{\overline{\rho u}_{i+1/2jk}}{2(\gamma-1)\hat \beta_{i+1/2jk}})  
+D_-^y(\frac{\overline{\rho v}_{ij+1/2k}}{2(\gamma-1)\hat \beta_{ij+1/2k}})  
+  D_-^y(\frac{\overline{\rho w}_{ijk+1/2}}{2(\gamma-1)\hat \beta_{ijk+1/2}})&=  \label{disc_IE}\\
D_-^x(\tilde \nu_{i+1/2jk}\frac{(\pf_x)_{i+1/2jk}}{\gamma-1})+
D_-^y(\tilde \nu_{ij+1/2k}\frac{(\pf_y)_{ij+1/2k}}{\gamma-1})+\nonumber \\
D_-^z(\tilde \nu_{ijk+1/2}\frac{(\pf_z)_{ijk+1/2}}{\gamma-1})\nonumber \\
+\kappa_r\left(D_-^xD^x_+T^4_{ijk}+D_-^yD^y_+T^4_{ijk}+D_-^z\kappa_rD^z_+T^4_{ijk}\right)&.\nonumber
\end{align}

\section{Stability}\label{sec:apriori_scheme}

Our aim is to derive the same a priori estimates for the semi-discrete system as  we found for the PDE system. To carry out the calculations, we need the following summation-by-parts formula,
\begin{align}
  \sum_{i=0}^{N}a_i\Delta_-b_{i+1/2}&=-a_0b_{-1/2}+a_Nb_{N+1/2}-\sum_{i=0}^{N-1}(\Delta_+a_i)b_{i+1/2},\label{SBP}
\end{align}
and product rule
\begin{align}
  \Delta_+a_ib_i&=\bar a_{i+1/2}\Delta_+b_i + \bar b_{i+1/2}\Delta_+a_i.\label{prod_rule}
\end{align}

\subsection{Conservation}

As in the continuous case, we obtain
\begin{align}
  \rho_{ijk},\, E_{ijk},p_i\in L^{\infty}_1,\quad \sqrt{\rho_{ijk}}|\vv_{ijk}|\in L^{\infty}_2,\label{diff_con} 
\end{align}
from conservation, positivity and the boundary conditions. 
\begin{remark}
In the notation above, $\rho_{ijk}$ is not a point value but symbolises the piecewise constant (in space) and continuous in time grid function on the entire domain that is bounded in $L^{\infty}(0,\T;L^1(\Omega))$
\end{remark}

\subsection{Entropy}

Here, we use the same scaling of the entropy as in \cite{Chandrashekar13}. Dropping the $ijk$ indices, the entropy variables are:
\begin{align}
\ww = \left(\frac{\gamma}{\gamma-1}-\frac{s}{\gamma-1}-\beta |\vv|^2,2\beta u, 2\beta v,2\beta w, -2\beta \right),\nonumber
\end{align}
where $|\vv|^2=u^2+v^2+w^2$ and the specific entropy,
\begin{align}
s=\log \frac{p}{\rho^\gamma}=-\log\beta-(\gamma-1)\log\rho-\log 2.\nonumber
\end{align}
Furthermore, the entropy potentials are: $\Psi^{x,y,z}=m_{1,2,3}=\rho \{u,v,w\}$.

An entropy estimate is obtained by contracting the scheme with the entropy variables, summing in space and integrating in time. We begin by multiplying the diffusive fluxes (without artificial viscosity) by $ \ww_{ijk}^T$ and summing in space. 

That is,
\begin{align}
  \sum_{ijk=0}^N \ww_{ijk}^T\left((\Delta_-^x\fb^\nu_{i+1/2jk})\Ss^x_{jk}+(\Delta^y_-\gb_{ij+1/2k}^\nu)\Ss^y_{ik}+(D_-^z\hb^\nu_{ijk+1/2})\Ss^z_{ij}\right)&=\nonumber \\
  \sum_{jk=0}^N (\ww_{Njk}^T\fb^\nu_{N+1/2jk}-\ww_{0jk}^T\fb^\nu_{-1/2jk})\Ss^x_{jk}
  +\sum_{ik=0}^N (\ww_{iNk}^T\gb_{iN+1/2k}^\nu-\ww_{i0k}^T\gb^\nu_{i,-1/2k})\Ss^y_{ik}\nonumber\\
  +\sum_{jk=0}^N (\ww_{ijN}^T\hb^\nu_{ijN+1/2}-\ww_{ij0}^T\hb^\nu_{ij,-1/2})\Ss^z_{ij}&\nonumber \\
  -\sum_{ijk=0}^{N-1,N,N} (\Delta_+^x\ww_{ijk})^T\fb^\nu_{i+1/2jk} \Ss^x_{jk}
  -\sum_{ijk=0}^{N,N-1,N}(\Delta^y_+\ww_{ijk})^T\gb_{ij+1/2k}^\nu\Ss^y_{ik}\label{diff_ent} \\
  -\sum_{ijk=0}^{N,N,N-1}(\Delta^z_+\ww_{ijk})^T\hb^\nu_{ijk+1/2}\Ss^z_{ij}&.\nonumber
\end{align}
The boundary conditions (\ref{num_noslip}) and numerical boundary fluxes (\ref{neumann_bc}) make all boundary terms vanish.

To further, manipulate (\ref{diff_ent}), we need the entropy-variable differences,
\begin{align}
  \Delta_+ \ww^1_i=&\frac{\Delta_+\rho_i}{\hat \rho_{i+1/2}}+\left(\frac{1}{(\gamma-1)\hat\beta_{i+1/2}}-\overline{|\uu|^2}_{i+1/2}\right)\Delta_+\beta_i,\nonumber\\
  &-2\bar u_{i+1/2} \bar \beta_{i+1/2}\Delta_+ u_i-2\bar v_{i+1/2} \bar \beta_{i+1/2}\Delta_+ v_i-2\bar w_{i+1/2} \bar \beta_{i+1/2}\Delta_+ w_i,\nonumber \\
\Delta_+ \ww^2_i=&2\bar \beta_{i+1/2} \Delta_+ u_i + 2\bar u_{i+1/2}\Delta_+ \beta_i,\label{deltaw} \\
\Delta_+ \ww^3_i=&2\bar \beta_{i+1/2} \Delta_+ v_i + 2\bar v_{i+1/2}\Delta_+ \beta_i,\nonumber \\
\Delta_+ \ww^4_i=&2\bar \beta_{i+1/2} \Delta_+ w_i + 2\bar w_{i+1/2}\Delta_+ \beta_i,\nonumber \\
\Delta_+ \ww^5_i=&-2\Delta_+ \beta_i.\nonumber 
\end{align}
The relations are stated generically to any difference direction. That is, $i$ represents one of the three indices and the other two are held constant. 
\begin{remark}
  A derivation of the identities (\ref{deltaw}) is found in \cite{Chandrashekar13}. We have verified their correctness and made the trivial extension to 3-D in the current treatise.
\end{remark}
Furthermore, we use (\ref{prod_rule}) to obtain the two auxiliary relations,
\begin{align}
  D_+^x(\rho_{ijk}u_{ijk})&=\bar \rho_{i+1/2}D_+^x(u_{ijk})+\bar u_{i+1/2jk}D_+^x(\rho_{ijk})\nonumber \\
\frac{1}{2}  D_+^x(\rho_i u^2_{i}) &= \frac{1}{2}  (D_+^x\rho_i) \overline {u^2}_{i+1/2} + \bar \rho_{i+1/2}\bar u_{i+1/2}D_+^xu_i, \nonumber  
\end{align}
which along with (\ref{deltaw}) is used to calculate the remaining (volume) diffusive terms in  (\ref{diff_ent}). There are a number of cancellations and the three directions are handled independently of each other. Suppressing the common $jk$-indices,  the terms in the x-direction from (\ref{diff_ent}) become,
\begin{align}
  (\Delta_+^x\ww_{i})\fb^\nu_{i+1/2}=\nu_{i+1/2jk}\left(\frac{(\Delta_+^x\rho_i)(D_+^x\rho_i)}{\hat \rho_{i+1/2}}\right.\nonumber \\
  +2\bar \beta_{i+1/2} \bar \rho_{i+1/2} \left((\Delta^x_+ u_i)(D^x_+ u_i)+(\Delta^x_+ v_i)(D^x_+ v_i)+(\Delta^x_+ w_i)(D^x_+ w_i)\right)\nonumber \\
 \left. -\frac{\bar \rho_{i+1/2}}{\gamma-1} \Delta_+^x\beta_{i}D_+^x\left(\frac{1}{\beta_{ijk}}\right)\right)-\frac{\kappa_r}{4}2(\Delta_+^x\beta_{ijk} )D_+^xT^4_i\label{ent_x}.
\end{align}
Furthermore, by (\ref{V}) and (\ref{Dnorm}), we have
\begin{align}
\sum_{ijk=0}^{N-1,N,N}(\Delta_+^x\ww_{ijk})\fb^\nu_{i+1/2jk}\Ss_{jk}=\sum_{ijk=0}^{N-1,N,N}(D_+^x\ww_{ijk})\fb^\nu_{i+1/2jk}V_{ijk},\nonumber
\end{align}
and similarly in the other directions.

Turning to the inviscid fluxes, we must verify that the convective numerical fluxes, augmented with artificial diffusion, are entropy dissipative in the following sense:
\begin{align}
  \Delta_+^x \ww^T_{ijk} (\fb^c_{i+1/2jk}-\fb^\lambda_{i+1/2jk})&\leq\Delta \Psi^x_{ijk}= \Delta^x_+ (\rho u)_{ijk},\nonumber \\
\Delta_+^y \ww^T_{ijk} (\gb^c_{ij+1/2k}-\gb^\lambda_{ij+1/2k})&\leq \Delta \Psi^y_{ijk}= \Delta^y_+ (\rho v)_{ijk},\label{EC} \\
\Delta_+^z \ww^T_{ijk} (\hb^c_{ijk+1/2}-\hb^d_{ijk+1/2})&\leq \Delta \Psi^z_{ijk}= \Delta^z_+ (\rho w)_{ijk}.\nonumber 
\end{align}
(The relations (\ref{EC}) are sometimes called \emph{Tadmor's shuffle conditions}.)

To reduce notation, we demonstrate that (\ref{EC}) holds in the one-dimensional case. (It is straightforward to extend the analysis to the three-dimensional case.) To carry out the derivation, we temporarily relabel $\ww^5$ as $\ww^3$ and similarly for the fluxes. In view of (\ref{EC}), we begin by calculating
\begin{align}
  LHS^c=(\Delta_+ \ww^1_i)\fb^{1,c}_{i+1/2}+(\Delta_+ \ww^2_i)\fb^{2,c}_{i+1/2}
  +(\Delta_+ \ww^3_i)\fb^{3,c}_{i+1/2}\nonumber
\end{align}
Using \eqref{deltaw} and \eqref{fc} (reduced to 1-D), we have
\begin{align}
  LHS^c=\nonumber \\
(\frac{\Delta_+\rho_i}{\hat \rho_{i+1/2}}+\left(\frac{1}{(\gamma-1)\hat\beta_{i+1/2}}-\overline{|u|^2}_{i+1/2}\right)\Delta_+\beta_i-2\bar u_{i+1/2} \bar \beta_{i+1/2}\Delta_+ u_i)(\overline{\rho u}_{i+1/2})\nonumber \\
+(2\bar \beta_{i+1/2} \Delta_+ u_i + 2\bar u_{i+1/2}\Delta_+ \beta_i)(\bar u_{i+1/2}\overline{\rho u}_{i+1/2}+p_{i+1/2})\nonumber \\
-2\Delta_+ \beta_i(\frac{1}{2(\gamma-1)\hat \beta_{i+1/2}}\overline {\rho u}_{i+1/2}-\frac{\overline{|u|^2}_{i+1/2}}{2}\overline {\rho u}_{i+1/2}\nonumber \\
+ \overline{|u|}^2_{i+1/2}\overline {\rho u}_{i+1/2}+p_{i+1/2}\bar u_{i+1/2})&.\nonumber
\end{align}
After some further manipulations, we obtain
\begin{align}
LHS^c=(\frac{\Delta_+\rho_i}{\hat \rho_{i+1/2}}+\left(\frac{1}{(\gamma-1)\hat\beta_{i+1/2}}-\overline{|u|^2}_{i+1/2}\right)\Delta_+\beta_i-2\bar u_{i+1/2} \bar \beta_{i+1/2}\Delta_+ u_i)(\overline{\rho u}_{i+1/2})\nonumber \\
+(2\bar \beta_{i+1/2} \Delta_+ u_i + 2\bar u_{i+1/2}\Delta_+ \beta_i)\bar u_{i+1/2}\overline{\rho u}_{i+1/2}+(2\bar \beta_{i+1/2} \Delta_+ u_i + 2\bar u_{i+1/2}\Delta_+ \beta_i)p_{i+1/2}\nonumber \\
-2\Delta_+ \beta_i(\frac{1}{2(\gamma-1)\hat \beta_{i+1/2}}\overline {\rho u}_{i+1/2}-\frac{\overline{|u|^2}_{i+1/2}}{2}\overline {\rho u}_{i+1/2}
+ \overline{|u|}^2_{i+1/2}\overline {\rho u}_{i+1/2}+p_{i+1/2}\bar u_{i+1/2})&.\nonumber
\end{align}
Cancelling velocity terms, leads to
\begin{align}
LHS^c=(\frac{\Delta_+\rho_i}{\hat \rho_{i+1/2}}+\left(\frac{1}{(\gamma-1)\hat\beta_{i+1/2}}\right)\Delta_+\beta_i)(\overline{\rho u}_{i+1/2})\nonumber \\
+(2\bar \beta_{i+1/2} \Delta_+ u_i + 2\bar u_{i+1/2}\Delta_+ \beta_i)p_{i+1/2}\nonumber \\
-2\Delta_+ \beta_i(\frac{1}{2(\gamma-1)\hat \beta_{i+1/2}}\overline {\rho u}_{i+1/2}
+p_{i+1/2}\bar u_{i+1/2})&.\nonumber
\end{align}
Cancelling pressure terms, results in
\begin{align}
LHS^c=(\frac{\Delta_+\rho_i}{\hat \rho_{i+1/2}})(\overline{\rho u}_{i+1/2})
+2\bar \beta_{i+1/2} \Delta_+ u_i p_{i+1/2}.\nonumber
\end{align}
Similarly, the artificial diffusion part of the left-hand side of (\ref{EC}) leads to
\begin{align}
LHS^{AD}=\lambda_{i+1/2}\left(\frac{(\Delta_+ \rho_i)^2}{\hat \rho_{i+1/2}}  +2\bar\beta_{i+1/2}\bar \rho_{i+1/2}(\Delta_+ u_i)^2 + \frac{\bar \rho_{i+1/2}}{\gamma-1}\Delta_+\beta_i(D_+\frac{1}{\beta_i})\right),\nonumber
\end{align}
in the same way as (\ref{ent_x}) was obtained. The complete left-hand side of (\ref{EC}), is given by $LHS=LHS^c+LHS^{AD}$. By substituting
\begin{align}
  \overline{\rho u}_{i+1/2}= \bar \rho_{i+1/2} \bar u_{i+1/2}+\frac{(\Delta_+ u_i)(\Delta_+ \rho_i)}{4}=\hat \rho_{i+1/2} \bar u_{i+1/2}+(\bar \rho_{i+1/2}-\hat \rho_{i+1/2}) \bar u_{i+1/2}+\frac{(\Delta_+ u_i)(\Delta_+ \rho_i)}{4},\nonumber
\end{align}
  we obtain
\begin{align}
LHS=(\frac{\Delta_+\rho_i}{\hat \rho_{i+1/2}})\left(\hat \rho_{i+1/2} \bar u_{i+1/2}+(\bar \rho_{i+1/2}-\hat \rho_{i+1/2}) \bar u_{i+1/2}+\frac{(\Delta_+ u_i)(\Delta_+ \rho_i)}{4}\right)
+2\bar \beta_{i+1/2} \Delta_+ u_i p_{i+1/2}\nonumber \\
+\lambda_{i+1/2}(\frac{(\Delta_+ \rho_i)^2}{\hat \rho_{i+1/2}}  +2\bar\beta_{i+1/2}\bar \rho_{i+1/2}(\Delta_+ u_i)^2 + \frac{\bar \rho_{i+1/2}}{\gamma-1}\Delta_+\beta_i(D_+\frac{1}{\beta_i}))&.\nonumber
\end{align}
Use $p_{i+1/2}=\frac{\bar \rho_{i+1/2}}{2\bar \beta_{i+1/2}}$ (defined in (\ref{fc})), to obtain
\begin{align}
LHS=\bar u_{i+1/2}(\Delta_+\rho_i)+\bar \rho_{i+1/2}(\Delta_+u_i)+\frac{(\Delta_+\rho_i)^2}{\hat \rho_{i+1/2}}\left(\frac{\bar \rho_{i+1/2}-\hat \rho_{i+1/2}}{\Delta_+\rho_i} \bar u_{i+1/2}+\frac{\Delta_+ u_i}{4}\right)
\nonumber \\
+\lambda_{i+1/2}(\frac{(\Delta_+ \rho_i)^2}{\hat \rho_{i+1/2}}  +2\bar\beta_{i+1/2}\bar \rho_{i+1/2}(\Delta_+ u_i)^2 + \frac{\bar \rho_{i+1/2}}{\gamma-1}\Delta_+\beta_i(D_+\frac{1}{\beta_i}))&.\nonumber 
\end{align}
The first two terms combine to $\Delta (\rho u)_i=\Delta \Psi_i$, which is the right-hand side of \eqref{EC}. The last terms are positive and the $\lambda_{i+1/2}$ bounds $|\frac{\bar \rho-\hat \rho}{\Delta_+\rho_i} \bar u+\frac{\Delta u}{4}|$. (In this expression and the next, we have suppressed the indices since all averages and differences are taken between $i,i+1$.) The latter is realised by the following calculation using (\ref{mean_rels}),
\begin{align}
  \left|\frac{\bar \rho-\hat \rho}{\Delta \rho} \right|=
  \left|\frac{\bar \rho-\frac{\Delta\rho}{\Delta\log\rho}}{\Delta \rho} \right|=
  \left|\frac{\bar \rho\Delta\log\rho-\Delta\rho}{\Delta\log\rho\Delta \rho} \right|=\left|\frac{1}{\Delta\log\rho}\left(\frac{\bar\rho}{\hat\rho}-1\right)\right|&=\nonumber\\
  \left|\frac{1}{\Delta\log\rho}\frac{\bar\rho-\hat\rho}{\hat\rho}\right|\leq \left|\frac{1}{\Delta\log\rho}\frac{\Delta\rho}{2\hat\rho}\right|=\left|\frac{\hat\rho}{2\hat \rho}\right| &=\frac{1}{2},   \label{2nd_order}
\end{align}
Hence, the fluxes (\ref{fc}),(\ref{gc}) and (\ref{hc}), augmented with the artificial diffusion fluxes,  satisfy (\ref{EC}).

Next, we use the above results when contracting the convective terms of the scheme (\ref{scheme}) with $\ww_{ijk}^T$, and use (\ref{EC}), to obtain
\begin{align}
  \sum_{ijk=0}^N \ww_{ijk}^T\left((\Delta_-^x(\fb^c_{i+1/2jk}-\fb^\lambda_{i+1/2jk}))\Ss^x_{jk}+(\Delta_-^y(\gb_{ij+1/2k}^c-\gb_{ij+1/2k}^\lambda))\Ss^y_{ik}\right.\nonumber\\
\left.  +(\Delta_-^z(\hb^c_{ijk+1/2}-\hb^\lambda_{ijk+1/2}))\Ss^z_{ij}\right)&\leq\nonumber \\
  \sum_{jk=0}^N (\ww_{Njk}^T\fb^c_{N+1/2jk}-\ww_{0jk}^T\fb^c_{-1/2jk}+\Psi^x_{0jk}-\Psi^x_{Njk})\Ss^x_{jk}\nonumber \\
  +\sum_{ik=0}^N (\ww_{iNk}^T\gb_{iN+1/2k}^c-\ww_{i0k}^T\gb^c_{i,-1/2k}+\Psi^y_{i0k}-\Psi^y_{iNk})\Ss^y_{ik}\nonumber\\
  +\sum_{jk=0}^N (\ww_{ijN}^T\hb^c_{ijN+1/2}-\ww_{ij0}^T\hb^c_{ij,-1/2}+\Psi^z_{ij0}-\Psi^z_{ikN})\Ss^z_{ij}&=0.\label{conv_ent}
\end{align}
All the boundary terms vanish  thanks to the no-slip condition implying that $w^{2,3,4}=0$; the first and fifth components are zero since the flux components are zero (see (\ref{no_slip_flux})); the entropy potentials also vanish due to the no-slip condition.

Combining (\ref{conv_ent}) and (\ref{diff_ent}), we compute the entropy estimate
\begin{align} 
\sum_{ijk=0}^N V_{ijk}(U_{ijk})_t
  +\sum_{ijk=0}^{N-1,N,N} (\Delta_+^x\ww_{ijk})^T\fb^d_{i+1/2jk} \Ss^x_{jk}\label{ent_est1} \\
  +\sum_{ijk=0}^{N,N-1,N}(\Delta^y_+\ww_{ijk})^T\gb_{ij+1/2k}^d\Ss^y_{ik}
  +\sum_{ijk=0}^{N,N,N-1}(\Delta^z_+\ww_{ijk})^T\hb^d_{ijk+1/2}\Ss^z_{ij}&\leq 0,\nonumber
\end{align}
where we have used $\ww_{ijk}^T(\uu_{ijk})_t= U(\uu_{ijk})_t=(U_{ijk})_t\quad$.

Integrating  (\ref{ent_est1}) and using (\ref{ent_x}) and (\ref{nu_approx}), we obtain (among others) the following bounds
\begin{align}
  \C\geq& \int_0^\T\sum_{ijk=0}^{N-1,N,N}\nu_{i+1/2jk}\frac{(D_+^x\rho_{ijk})^2}{\hat \rho}V_{ijk}\,dt\nonumber \\
  =&\int_0^\T
  \sum_{ijk=0}^{N-1,N,N}(\mu_1\bar \rho_{i+1/2jk}+\frac{\mu_0}{\hat\rho_{i+1/2jk}})\frac{(D_+^x\rho_{ijk})^2}{\hat \rho_{i+1/2jk}}V_{ijk}\,dt\nonumber\\
=&\int_0^\T\sum_{ijk=0}^{N-1,N,N}(\mu_1\frac{\bar \rho_{i+1/2jk}}{\hat\rho_{i+1/2jk}}(D_+^x\rho_{ijk})^2 +\mu_0(D_+^x\log\rho_{ijk})^2)V_{ijk}\,dt,
  \label{Dlogrho}
\end{align}
where $\C$ is a constant obtained from the initial data.
Noting that, $\bar \rho/\hat \rho\geq 1$, we obtain bounds on $D_+^x\rho_{ijk}$ and $D_+^x\log\rho_{ijk}$. Along with (\ref{diff_con}), we have,
\begin{align}
  \rho_{ijk}&\in L^2(0,\T;H^1(\Omega))\label{disc_rho_H1},\\
  \nabla_+\log\rho_{ijk} &\in L^2(0,\T;L^2(\Omega)),\nonumber
\end{align}
where $\nabla_+$ is short for the finite differences $(D_+^x,D_+^y,D_+^z)$

In the same way as in the derivation following (\ref{U_est}), and again by using that $\frac{\bar \rho_{i+1/2jk}}{\hat\rho_{i+1/2jk}}\geq 1$, we obtain the discrete counterpart of (\ref{rlogr}) and (\ref{ent_diff}):
\begin{align}
  \rho_{ijk}\log(\rho_{ijk})&\in L^{\infty}(0,\T,L^1(\Omega)),\nonumber\\
  \rho_{ijk}\log(T_{ijk})&\in L^{\infty}(0,\T,L^1(\Omega)),\label{rlogr_d}\\
  \nabla_+\log(T_{ijk})&\in L^{2}(0,\T,L^2(\Omega)).\nonumber
\end{align} 
(We remark that we get all the counterparts of (\ref{ent_diff}) but limit ourselves to list the estimates necessary for the convergence proof below.)

Next, we turn to the temperature estimate from radiation in (\ref{ent_est1}). We have,
\begin{align}
  \C\geq \int_0^\T \sum_{ijk=0}^N-V_{ijk}\kappa_r2(D_+^x\beta_{ijk} ) D_+^xT^4_{ijk}\,dt,
  \nonumber
\end{align}
where
\begin{align}
D_+^x\beta_{ijk} =-\frac{D_+^xT_{ijk}}{2R\check T_{i+1/2jk}^2},\nonumber
\end{align}
and
\begin{align}
  D_+^x(T_{i+1/2jk})^4 &=4\widetilde{T^3}_{i+1/2jk}D_+^xT_i\nonumber \\
  4\widetilde{T^3}_{i+1/2jk}&=(T_{i+1}^3+T_i^3+T_i^2T_{i+1}+T_{i+1}^2T_i)_{jk}.\nonumber
\end{align}
Hence, we have a bound on
\begin{align}
\frac{\sqrt{\widetilde{T^3}_{i+1/2jk}}}{\sqrt{T_{i+1jk}T_{ijk}}}D_+^xT_i\in L^2_2\nonumber.
\end{align}
We wish to show that
\begin{align}
D_+^xT^{3/2}_{ijk}\in L^2_2.\label{disc_temp}
\end{align}
Hence, we need to demonstrate that
\begin{align}
\frac{\sqrt{\widetilde{T^3}_{i+1/2jk}}}{\sqrt{T_{i+1jk}T_{ijk}}}\gtrsim \frac{D_+^xT^{3/2}_{ijk}}{ D^x_+T_{ijk}}.\label{cond_T}
\end{align}
We note that $\frac{{\widetilde{T^3}_{i+1/2jk}}}{T_{i+1jk}T_{ijk}}\gtrsim T_{i+1jk}+T_{ijk}$. That is, the left-hand side majorises the arithmetic mean (with an appropriate scaling) and twice the arithmetic mean majorises both $T_{ijk}$ and $T_{i+1jk}$. Furthermore,  by the mean-value theorem
\begin{align}
  \frac{(D_+^xT^{3/2}_{ijk})}{ (D^x_+T_{ijk})}=\frac{3}{2}\sqrt{T^*}, \quad T^*\in (T_i,T_{i+1}),\nonumber
\end{align}
showing that it is a monotone average. Since it is a monotone average, it is majorised by either end points and therefore by the left-hand side of (\ref{cond_T}). We conclude that (\ref{disc_temp}), holds. (The same arguments apply in the y- and z-directions.)

By repeating the arguments in the continuous case, we have a $|\B|$ (as in (\ref{B})) of non-zero measure and  $T_i\in L^{\infty}(0,\T;L^1(\B))$. Hence,
\begin{align}
  T_{ijk}^{3/2}&\in L^2(H^1),\nonumber \\
  T_{ijk}&\in L^3(L^9).\label{disc_temp_bound}
\end{align}
Finally, the bound on $\frac{(D_+^xT_{ijk})^2}{\check T^2_{i+1/2jk}}\widetilde{T^3_{i+1/2jk}}$ implies
\begin{align}
  T_{ijk}^{3/2}\nabla_+\log T_{ijk}&\in L^2(L^2).\label{strong_DlogT}
\end{align}

\subsection{Kinetic energy}

Summing (\ref{KE_disc}) by parts over the spatial domain, results in,
\begin{align}
  \sum_{ijk=1}^nV_{ijk}\frac{d}{dt}K_{ijk}=
  \sum_{jk=1}^{N,N}\left(\Ff_{-1/2jk}-\Ff_{N+1/2jk} \right)\Ss^x_{jk}\nonumber \\
  + \sum_{ik=1}^{N,N} \left(\Gf_{i,-1/2k}-\Gf_{i,N+1/2k} \right)\Ss^y_{ik}
  + \sum_{ij=1}^{N,N} \left( \Hf_{ij,-1/2}-\Hf_{ij,N+1/2}\right)\Ss^z_{ij}\label{KE_est1} \\
 -\sum_{ijk=1}^{N}V_{ijk}\left(\Df_{ijk}-(pDv)_{ijk}\right).\nonumber
\end{align}
The fluxes were defined in (\ref{KE_flux}) and vanish at the boundary thanks to no-slip.

As in the continuous case, we split the pressure term, using
\begin{align}
-\sum_{ijk=1}^{N-1,N,N} V_{ijk}(D_+^xu_{ijk})\frac{\bar \rho_{i+1/2jk}}{2\bar \beta_{i+1/2jk}}\leq \|\rho_{i+1/2jk}D_+^xu_{ijk}\|^2_2+\|\frac{1}{2\bar \beta_{i+1/2jk}}\|_2^2,\nonumber
  \end{align}
and similarly for the other terms.
Since $\nu_{i+1/2jk}$ contains a $\bar \rho_{i+1/2jk}$ term (see (\ref{nu_approx})) the first term on the right-hand side is controlled by $\sum_{ijk=1}^{N}V_{ijk}\Df_{ijk}$. The second term, is controlled by (\ref{disc_temp_bound}) via
\begin{align}
  \frac{1}{\bar \beta_{i+1/2}} = \frac{1}{\frac{1}{2RT_{i+1}}+\frac{1}{2RT_i}}=
\frac{1}{\frac{T_i+T_{i+1}}{2RT_{i+1}T_i}}=\frac{RT_{i+1}T_i}{\frac{T_i+T_{i+1}}{2}}\leq \frac{R\check T_{i+1/2}^2}{\bar T_{i+1/2}}\leq  R \bar  T_{i+1/2},\nonumber
\end{align}
where the last inequality follows from (\ref{mean_rels}) and positivity of temperature. 

Using that $\bar \rho/\hat\rho\geq 1$, we conclude that
\begin{align}
  \bar \rho_{i+1/2jk}\left(|D_+^xu_{ijk}|+|D_+^xv_{ijk}|  +|D_+^xw_{ijk}|\right)&\in L^2(0,\T;L^2),\nonumber\\
  |D_+^xu_{ijk}|+|D_+^xv_{ijk}|  +|D_+^xw_{ijk}|&\in L^2(0,\T;L^2),\quad \textrm{and}\label{disc_v_est}\\
  \sqrt{\frac{\bar\rho_{i+1/2jk}}{\hat\rho_{i+1/2jk}}}\left(|D_+^xu_{ijk}|+|D_+^xv_{ijk}|  +|D_+^xw_{ijk}|\right)&\in L^2(0,\T;L^2),\nonumber
\end{align}
and similarly in the other two directions. These bounds are the discrete equivalents of (\ref{v_est}) and also imply that
\begin{align}
\vv_{ijk}\in L^2(H^1) \cap L^2(L^6).\label{disc_v_H1}
\end{align}
We will also need the complete estimate obtained from $\Df_{ijk}$:
\begin{align}
\tilde \nu_{i+1/2jk}\bar \rho_{i+1/2jk}\left((D_+^xu_{ijk})^2+(D_+^xv_{ijk})^2  +(D_+^xw_{ijk})^2\right)&\in L^1(0,\T;L^1),\label{ke_nu}
\end{align}
and the corresponding bounds  in the other two directions.

\subsection{Specific volume}

We need the following identities
\begin{align}
  D_+^x\frac{1}{\rho_{i}^2}=\frac{1}{\Delta_+^x x_{i}}\left(\frac{1}{\rho_{i+1}}+\frac{1}{\rho_{i}}\right)\left(\frac{1}{\rho_{i+1}}-\frac{1}{\rho_{i}}\right)&=\nonumber \\
-\frac{1}{\Delta_+^x x_{i}}\frac{\rho_{i+1}+\rho_{i}}{\rho_{i+1}\rho_{i}}\frac{\rho_{i+1}-\rho_{i}}{\rho_{i+1}\rho_{i}}=-\frac{2\bar \rho_{i+1/2}}{\check \rho_{i+1/2}^4}D_+^x\rho_{i},\nonumber
 \end{align}
and
\begin{align}
D_+^x\frac{1}{\rho_i}=-\frac{1}{\check \rho_{i+1/2}^2}D_+^x\rho_i,\nonumber
\end{align}
where $\check\rho_{i+1/2}$ is given by (\ref{means}).
\begin{remark}
To reduce notation, we take the grid to be uniform in each direction, such that $\Delta_+^x x_{i}=\Delta_+^x x_{i+1/2}=h$. This is not a necessity and all estimates are derivable on a non-uniform grid. In fact, the scheme (and the on-going convergence proof) is readily generalised unstructured to Voronoi grids as well.
\end{remark}

Next, multiply the continuity equation of (\ref{scheme}) by $-1/\rho_{ijk}^2$,
\begin{align}
  V_{ijk}(\rho^{-1}_{ijk})_t &= \frac{1}{\rho_{ijk}^2}\left( \Delta_-^x\fb^{c\lambda,1})_{i+1/2jk}\Ss^x_{jk} + \Delta_-^y(\gb^{c\lambda,1})_{ij+1/2k}\Ss_{ik}^y+\Delta_-^z(\hb^{c\lambda,1})_{ijk+1/2}\Ss_{ij}^z \right)\nonumber \\
 & - \frac{1}{\rho_{ijk}^2}\left(\Delta_-^x(\nu_{i+1/2jk}D_+^x\rho_{ijk})\Ss_{jk}^x+
  \Delta_-^y(\nu_{ij+1/2k}D_+^y\rho_{ijk})\Ss_{ik}^y+
  \Delta_-^z(\nu_{ijk+1/2}D_+^z\rho_{ijk})\Ss_{ij}^z  \right)\label{spec_vol}
\end{align}
where $\fb^{c\lambda,1}=\fb^{c,1}-\fb^{\lambda,1}$ etc. We rewrite the flux as,
\begin{align}
  \fb^{c\lambda,1}_{i+1/2jk}=\overline{\rho u}_{i+1/2jk} -\lambda_{i+1/2jk} \Delta_+^x\rho_{ijk}=\rho^*_{i+1/2jk} \bar u_{i+1/2jk}-\lambda^*_{i+1/2jk}\Delta_+^x\rho_{ijk}\nonumber
\end{align}
where $\lambda^*=\lambda-\frac{ \rho^* -\bar \rho}{\Delta \rho}\bar u+\frac{\Delta u\Delta \rho}{4}\geq 0$. (Again, we momentarily suppress the indices $ijk$ and $i+1jk$.) Furthermore, we choose
\begin{align}
\rho^*=\frac{\check \rho^2}{\sqrt{\hat\rho}\sqrt{\bar \rho}}\nonumber,
\end{align}
and calculate
\begin{align}
|\rho^*-\bar \rho|&=\left|\frac{\check \rho^2}{\sqrt{\hat\rho}\sqrt{\bar \rho}}-\bar \rho\right|= \frac{|\check \rho^2-\bar\rho^{3/2}\sqrt{\hat \rho}|}{\sqrt{\hat\rho}\sqrt{\bar \rho}}.\nonumber
\end{align}
Since $\check \rho \leq \bar \rho$ and $\hat\rho\leq \bar \rho$, we have
\begin{align}
  \frac{|\rho^*-\bar \rho|}{\Delta\rho}\leq \frac{|\check \rho^2-\bar\rho^{2}|}{\sqrt{\hat\rho}\sqrt{\bar \rho}\Delta\rho}=\frac{|(\Delta\rho)^2|}{4\sqrt{\hat\rho}\sqrt{\bar \rho}\Delta\rho}\leq \frac{|(\Delta\rho)^2|}{4\hat\rho\Delta\rho}=\frac{1}{4}|\Delta\log\rho|.\nonumber
\end{align}
In view of (\ref{ad_coeff2}), the resulting artificial diffusion coefficient is
\begin{align}
\lambda^*=\lambda-(\frac{\rho^*-\bar \rho}{\Delta \rho}\bar u-\frac{\Delta u}{4})\geq 0.\nonumber
\end{align}
We introduce the notation, $\nu^*_{i+1/2jk}=\nu_{i+1/2jk}+(\Delta_+^xx_{i})\lambda^*_{i+1/2jk}$ in analogy with (\ref{tilde_nu}).

Next, we sum the x-terms of (\ref{spec_vol}) by parts and use the boundary conditions. Furthermore, we choose $\eta$ such that $\mu_0/2>\eta>0$, and obtain
\begin{align}
  \sum_{i=0}^N \Ss_{jk}^x \left(\frac{1}{\rho_{ijk}^2} \Delta_-^x( \rho^* \bar u)_{i+1/2jk} - \frac{1}{\rho_{ijk}^2}\Delta_-^x(\nu^*_{i+1/2jk}D_+^x\rho_{ijk})\right)&=\nonumber \\
  \sum_{i=0}^{N-1} \Ss_{jk}^x\Delta^x_+ x_{i}\left( \frac{2\bar \rho_{i+1/2jk}}{\check \rho_{i+1/2jk}^4}(D_+^x\rho_{ijk}) \rho^*_{i+1/2jk} \bar u_{i+1/2jk} -\nu^*_{i+1/2jk} \frac{2\bar \rho_{i+1/2jk}}{\check \rho_{i+1/2jk}^4}(D_+^x\rho_{ijk})^2\right)&\leq\nonumber \\
  \sum_{i=0}^{N-1} V_{ijk}\left( -\frac{2\bar \rho_{i+1/2jk}\rho^*_{i+1/2jk}}{\check \rho_{i+1/2jk}^2}(D_+^x\frac{1}{\rho_{ijk}}) \bar u_{i+1/2jk} -\nu^*_{i+1/2jk} 2\bar \rho_{i+1/2jk}(D_+^x\frac{1}{\rho_{ijk}})^2\right)&\leq\nonumber \\
   \sum_{i=0}^{N-1} V_{ijk}\left(\eta|\frac{2\bar \rho_{i+1/2}\rho^*_{i+1/2jk}}{\check \rho_{i+1/2jk}^2}D_+^x\frac{1}{\rho_{ijk}}|^2+ \eta^{-1}|\bar u_{i+1/2jk}|^2 - \frac{2\mu_0\bar \rho_{i+1/2jk}}{\hat \rho_{i+1/2jk}}|D_+^x\frac{1}{\rho_{ijk}}|^2\right)&,\nonumber
\end{align}
where we have used (\ref{nu_approx}) in the last row. Since $\frac{\bar \rho_{ijk}}{\check \rho_{ijk}}>1$, and by summing over $jk$, we obtain the bounds
\begin{align}
\sup_{0,\T}\|\rho^{-1}_{ijk}\|_1\leq \C,&\quad
\int_0^\T\|\nabla_+ \rho^{-1}_{ijk}\|_2^2\,dt,\leq \C,\label{svolume_disc}\\
\rho_{ijk}^{-1}\in L^2(H^1),&\quad
\rho_{ijk}^{-1}\in L^2(L^6).\nonumber
\end{align}
(Recall that $\nabla_+$ is short for the finite differences $(D_+^x,D_+^y,D_+^z)$.) These estimates correspond to (\ref{svolume_bounds})-(\ref{svolume_b2})). Furthermore, we get
\begin{align}
\sqrt{T_{ijk}}\in L^{\infty}_1,\label{T_infty_disc}
\end{align}
as in (\ref{T_infty}).

\subsection{Estimates from the continuity equation}

Multiply the continuity equation by $\rho_{ijk}$ and sum in space.
\begin{align}
\sum_{i,j,k=0}^{N}  \frac{1}{2} V_{ijk}(\rho^{2}_{ijk})_t =\nonumber\\
\sum_{i,j,k=0}^{N}   -\rho_{ijk}V_{ijk}\left( D_-^x(\fb^{c,1})_{i+1/2jk} + D_-^y(\gb^{c,1})_{ij+1/2k}+D_-^z(\hb^{c,1})_{ijk+1/2} \right)\nonumber \\
+\sum_{i,j,k=0}^{N}   \rho_{ijk}V_{ijk}\left(D_-^x\tilde \nu_{i+1/2jk}D_+^x\rho_{ijk}+
  D_-^y\tilde \nu_{ij+1/2k}D_+^y\rho_{ijk}+
  D_-^z\tilde \nu_{ijk+1/2}D_+^z\rho_{ijk}  \right).\nonumber
\end{align}
We sum by parts and thanks to the boundary conditions, the boundary terms vanish. (Again, all three directions are handled analogously and we focus on the x-direction and denote the remaining terms by $YZ$)
\begin{align}
\sum_{i,j,k=0}^{N}  \frac{1}{2} V_{ijk}(\rho^{2}_{ijk})_t &= \sum_{i,j,k=0}^{N-1,N,N} \left((D_+^x\rho_{ijk})\overline{\rho  u}_{i+1/2jk}  
 -(D_+^x \rho_{ijk})\tilde \nu_{i+1/2jk}D_+^x\rho_{ijk}\right)+YZ.\nonumber
\end{align}
Use $\overline{\rho  u}_{i+1/2}=\sqrt{\rho_{i}}\sqrt{\rho_i}u_i+\sqrt{\rho_{i+1}}\sqrt{\rho_{i+1}}u_{i+1}$ which allows us to proceed as in the continuous case. We have $L^{\infty}_2$ bounds of $\sqrt{\rho_i}u_i$. We need the diffusive term to provide a bound  $\sqrt{\rho_i}D_+^x\rho_i\in L^2_2$ and $\sqrt{\rho_{i+1}}D_+^x\rho_i$, which is what we get since $\tilde \nu_{i+1/2}\sim \bar \rho_{i+1/2}$ and $\sqrt{\bar \rho_{i+1/2}}\geq \frac{1}{2}(\sqrt{\rho_i}+\sqrt{\rho_{i+1}})$. We collect the estimates,
\begin{align}
  \rho_{ijk}\in L^\infty_2,\label{disc_rho_bound}\\
  \sqrt{\bar \rho_{i+1/2jk}}D_+^x\rho_{ijk}\in L^\infty_2.\nonumber
\end{align}

\vspace{0.5cm}

Next, we test the equation with $\rho_i^2$. The procedure follows the same route as above and we turn directly to the convective terms. In the x-direction, we have
\begin{align}
 \sum_{i=0}^N \rho_i^2D_-^x\fb^{c\lambda,1}_{i+1/2}=-\sum_{i=0}^{N-1} (D_+^x\rho_i^2)\fb^{c\lambda,1}_{i+1/2},\nonumber
\end{align}
where we have dropped the sum over $jk$ and suppressed those indices. As above, we denote $\fb^{c\lambda,1}=\fb^{c,1}-\fb^{\lambda,1}$. The boundary terms vanish thanks to the boundary conditions. Some further manipulations yield,
\begin{align}
  \sum_{i=0}^{N-1} (D_+^x\rho_i^2)\fb^{c\lambda,1}_{i+1/2}&=\nonumber\\
  \sum_{i=0}^{N-1} (D_+^x\rho_i^2)\left(\rho^*_{i+1/2}\bar u_{i+1/2}+\underbrace{\frac{\bar \rho_{i+1/2}-\rho^*_{i+1/2}}{\Delta_+^x\rho_i}\bar u_{i+1/2}\Delta_+^x\rho_i+\frac{\Delta_+^x u_i}{4}\Delta_+^x\rho_i-\lambda_{i+1/2}\Delta_+^x\rho_i}_{RT_{i+1/2}}\right)\nonumber.
\end{align}
Here, we choose the average:
\begin{align}
\rho^*_{i+1/2} = \frac{2}{3}\frac{D_+^x\rho_i^3}{D_+^x\rho_i^2}.\nonumber
\end{align}
Noting that  $(D_+^x\rho_i^2)=2\bar \rho_{i+1/2}D_+^x\rho_i$, the rest terms ($RT_{i+1/2}$) are bounded since
\begin{align}
\lambda_{i+1/2}\geq |\frac{\bar \rho_{i+1/2}-\rho^*_{i+1/2}}{\Delta_+^x\rho_i}\bar u_{i+1/2}+\frac{\Delta_+^x u_i}{4}|.\nonumber
\end{align}
The last inequality follows from the following calculation,
\begin{align}
  \frac{\bar \rho_{i+1/2}-\rho^*_{i+1/2}}{\Delta_+^x\rho_i}=
  \frac{\bar \rho_{i+1/2}-\frac{2}{3}\frac{\rho_{i+1}^3-\rho_i^3}{\rho_{i+1}^2-\rho_i^2}}{\Delta_+^x\rho_i}&=\nonumber\\
  \frac{\bar \rho_{i+1/2}-\frac{2(\rho_{i+1}^2+\rho_i^2+\rho_{i+1}\rho_i)}{3(\rho_{i+1}+\rho_i)}}{\Delta_+^x\rho_i}=
  \frac{\frac{3(\rho_{i+1}+\rho_i)^2-4(\rho_{i+1}^2+\rho_i^2+\rho_{i+1}\rho_i)}{6(\rho_{i+1}+\rho_i)}}{\Delta_+^x\rho_i}&=\nonumber\\
  \frac{-\rho_{i+1}^2-\rho_i^2+2\rho_{i+1}\rho_i}{6(\rho_{i+1}+\rho_i)\Delta_+^x\rho_i}=
  \frac{-(\rho_{i+1}-\rho_i)^2}{6(\rho_{i+1}+\rho_i)\Delta_+^x\rho_i}&=\nonumber\\
 - \frac{(\rho_{i+1}-\rho_i)}{6(\rho_{i+1}+\rho_i)}&,\nonumber
\end{align}
and (\ref{ad_coeff2}).

Using (\ref{SBP}), 
\begin{align}
  \sum_{i=0}^N \rho_i^2D_-^x\fb^{c\lambda,1}_{i+1/2}= -\sum_{i=0}^{N-1} (D_+^x\rho_i^2)\fb^{c\lambda,1}_{i+1/2}=-\frac{2}{3}\sum_{i=0}^{N-1} (D_+^x\rho_{i}^3) \bar u_{i+1/2}+(D_+^x\rho_{i}^2)RT_{i+1/2},\nonumber
  \end{align}
  where $RT$ denotes a damping term from the rest terms and artificial diffusion and we have used $u_{0jk}=u_{Njk}=0$ in the summation-by-parts step. We drop the rest terms, as they have the correct sign and recast the sums using the no-slip condition ($u_0=u_N=0$),
  \begin{align}
  \frac{2}{3}\sum_{i=0}^{N-1} (D_+^x\rho_{i}^3) \bar u_{i+1/2}=\frac{2}{3}\sum_{i=0}^{N-1} (D_+^x\rho_{i}^3) \frac{u_{i+1}+u_i}{2}=
  \frac{1}{3}\sum_{i=0}^{N-1} (D_+^x\rho_{i}^3) u_{i+1}
 +\frac{1}{3}\sum_{i=0}^{N-1} (D_+^x\rho_{i}^3) u_i\nonumber\\
   \frac{1}{3}\sum_{i=1}^{N-1} (D_+^x\rho_{i-1}^3) u_{i}
 +\frac{1}{3}\sum_{i=1}^{N-1} (D_+^x\rho_{i}^3) u_i=
   -\frac{1}{3}\sum_{i=1}^{N} \rho_{i-1}^3 D_-^xu_{i}
 -\frac{1}{3}\sum_{i=1}^{N} \rho_{i}^3 D_-^xu_i&=\nonumber\\
    -\frac{2}{3}\sum_{i=1}^{N} \overline{\rho^3}_{i-1/2}D_-^x u_{i}=
        -\frac{2}{3}\sum_{i=0}^{N-1} \overline{\rho^3}_{i+1/2}D_+^x u_{i}&.\label{spec_SBP}	
\end{align}
Since $\bar \rho_{i+1/2jk} D_+^xu_{ijk} \in L^2_2$ by (\ref{disc_v_est}), we can proceed as in the continuous derivation.  First,  
\begin{align}
\rho_{i+1}^3+\rho_i^3=(\rho_{i+1}+\rho_i)(\rho_{i+1}^2-\rho_i\rho_{i+1}+\rho_i^2).\nonumber
\end{align}
We can now split the term $|\overline{\rho^3}_{i+1/2}D_+^x u_{i}| \leq \eta (\rho_{i+1}^2-\rho_i\rho_{i+1}+\rho_i^2)^2 + \eta^{-1}((\rho_{i+1}+\rho_i)D_+^x u_{i})^2$  Hence, we need $\rho_{i+1}^2-\rho_i\rho_{i+1}+\rho_i^2\in L^2_2$. This follows if $\rho_i\in L^4_4$, which we will get from the diffusive term as in the continuous case.

To this end, we need to recast the diffusive term to bound the gradient of $\rho_{ijk}^2$. In the x-direction, we have
\begin{align}
  \sum_{i=0}^N \rho_i^2D_-^x(\nu_{i+1/2}D_+^x\rho_i)=-\sum_{i=0}^{N-1} (D_+^x\rho_i^2)\nu_{i+1/2}D_+^x \rho_i&=\nonumber \\
  -\sum_{i=0}^{N-1} 2\bar\rho_{i+1/2}(D_+^x\rho_i)\nu_{i+1/2}D_+\rho_i&.\nonumber
\end{align}
Since $\nu_{i+1/2}\sim\bar \rho_{i+1/2}$ and seeing that $\bar \rho_{i+1/2}D_+\rho_i=\frac{1}{2}D_+\rho_i^2$, we obtain the necessary estimate to bound $\rho_i\in L^4_4$. We obtain all the estimates of (\ref{strong_rho}):
\begin{align}
  \rho_{ijk}\in L^{\infty}_3, \quad\quad
  \rho_{ijk}^2\in L^{2}(H^1), \quad\quad
  \rho_{ijk}\in L^4_{12}, \quad\quad 
  \rho_{ijk}\in L^{6}_{6}. \label{disc_strong_rho} 
\end{align}

\vspace{0.5cm}

Next, we test against $\rho_i^{3}$. As above, we begin with the convective term (and drop the common $jk$ indices),
\begin{align}
 \sum_{i=0}^N \rho_i^3D_-^x\fb^{c\lambda,1}_{i+1/2}=-\sum_{i=0}^{N-1}  D_+^x\rho_i^3\fb^{c\lambda,1}_{i+1/2},\nonumber
\end{align}
and rewrite the flux as,
\begin{align}
\fb^{c\lambda,1}_{i+1/2}=\rho^*_{i+1/2}\bar u_{i+1/2} + (\bar \rho_{i+1/2}-\rho^*_{i+1/2})\bar u_{i+1/2} -\lambda^a_{i+1/2}\Delta^x_+\rho_i.\label{rewrite}
\end{align}
This time we choose,
\begin{align}
\rho^*_{i+1/2}=\frac{3}{4}\frac{\Delta^x_+\rho_i^4}{\Delta^x_+\rho_i^3}=\frac{3}{4}\frac{(\rho_{i+1}^2+\rho_i^2)(\rho_{i+1}+\rho_i)}{\rho_{i+1}^2+\rho_{i+1}\rho_i+\rho_i^2}\label{rho_mean4}.
\end{align}
Since, by the mean value theorem
\begin{align}
\frac{\Delta^x_+\rho_i^4}{\Delta^x_+\rho_i^3}=\frac{(\rho^3_{i+1})^{4/3}-(\rho_i)^{4/3}}{\rho_{i+1}^3-\rho_{i}^3}=\frac{4}{3}(\rho_\star^3)^{1/3},\nonumber
\end{align}
for some $\rho_\star\in (\rho_i,\rho_{i+1})$, we conclude that $\rho^*_{i+1/2}\in(\rho_i,\rho_{i+1})$. (That is, it is a monotone average.)

As before, we must verify that the artificial diffusion term dominates the error term in (\ref{rewrite}). That is, the following expression must be positive:
\begin{align}
 \lambda^a_{i+1/2}- \frac{\bar \rho_{i+1/2}-\rho^*_{i+1/2}}{\Delta^+_x\rho_i}\bar u_{i+1/2}&=\nonumber \\
\lambda_{i+1/2} -\frac{\Delta_+^x u_i}{4}- \frac{\bar \rho_{i+1/2}-\rho^*_{i+1/2}}{\Delta^x_+\rho_i}\bar u_{i+1/2}&.\label{rho4_ad}
\end{align}
Using (\ref{rho_mean4}), we calculate
\begin{align}
  \bar \rho_{i+1/2}-\rho^*_{i+1/2}=
  \bar \rho_{i+1/2}-\frac{3}{4}\frac{(\rho_{i+1}^2+\rho_i^2)(\rho_{i+1}+\rho_i)}{\rho_{i+1}^2+\rho_{i+1}\rho_i+\rho_i^2}&=\nonumber\\
  \bar \rho_{i+1/2}\left(1-\frac{3}{2}\frac{(\rho_{i+1}^2+\rho_i^2)}{\rho_{i+1}^2+\rho_{i+1}\rho_i+\rho_i^2}\right)=
  \bar \rho_{i+1/2}\left(\frac{\rho_{i+1}^2+\rho_{i+1}\rho_i+\rho_i^2-\frac{3}{2}(\rho_{i+1}^2+\rho_i^2)}{\rho_{i+1}^2+\rho_{i+1}\rho_i+\rho_i^2}\right)&=\nonumber\\
  -\frac{\bar \rho_{i+1/2}}{2}\left(\frac{(\rho_{i+1}-\rho_i)^2}{\rho_{i+1}^2+\rho_{i+1}\rho_i+\rho_i^2}\right).\nonumber
\end{align}
Using the last expression  and the definition of the artificial diffusion coefficient (\ref{ad_coeff2}), it is clear that (\ref{rho4_ad}) is positive.

Next, we recast convective term. Using (\ref{rho_mean4}) and the analogous derivation as (\ref{spec_SBP}), we get
\begin{align}
  -\sum_{i=0}^{N-1}  (D_+^x\rho_i^3)\rho^*_{i+1/2}\bar u_{i+1/2}=-\frac{3}{4}\sum_{i=0}^{N-1} (D_+^x\rho_i^4)\bar u_{i+1/2}=
  \frac{3}{4}\sum_{i=0}^{N-1} \overline{\rho^4}_{i+1/2} D_+^x u_{i}
  \nonumber.
\end{align}
We proceed as in the previous estimate and factorise $\overline{\rho^4}_{i+1/2}=\bar \rho_{i+1/2}\frac{\overline{\rho^4}_{i+1/2}}{\bar \rho_{i+1/2}}$ and use (\ref{disc_v_est}) to bound one part. The other is bounded since $\frac{\overline{\rho^4}_{i+1/2}}{\bar \rho_{i+1/2}}$ can be majorised by $\max_{j=i,i+1}|\rho_j^3|$, which in turn is bounded in $L^2_2$ by (\ref{disc_strong_rho}).

Turning to the diffusive terms, we have
\begin{align}
    \sum_{i=0}^N \rho_i^3D_+^x\fb^{d,1}_{i+1/2}=-\sum_{i=0}^{N-1}  (D_+^x\rho_i^3)\fb^{d,1}_{i+1/2}=
    -\sum_{i=0}^{N-1} D_+^x\rho_i^3\nu_{i+1/2}D_+\rho_{i}&=\nonumber\\
    -\sum_{i=0}^{N-1}  (\rho^2_{i+1}+\rho_i\rho_{i+1}+\rho_i^2)\nu_{i+1/2}(D_+\rho_{i})^2&.\label{diff_rho_high}
\end{align}
We use that $\nu_{i+1/2}\sim \rho_{i+1/2}$, which implies that $(\rho^2_{i+1}+\rho_i\rho_{i+1}+\rho_i^2)\nu_{i+1/2}$ is dominating $c\cdot (\max(\rho_i,\rho_{i+1}))$ for some $c>0$. Hence, (\ref{diff_rho_high}) dominates $D_+^x\rho_{i}^{5/2}$ since as, by the same argument as above
\begin{align}
\frac{D_+^x\rho_{i}^{5/2}}{D_+^x\rho_{i}}=\frac{5}{2}\rho_{\star}^{3/2},\nonumber
\end{align}
is an average and by the mean-value theorem, $\rho_\star\in(\rho_i,\rho_{i+1})$. Hence, we obtain 
\begin{align}
  \rho\in L^{\infty}_{4},\quad
  \nabla_+ \rho^{5/2} \in L^2_2\quad
  \rho\in L^5_{15},\quad
  \rho\in L^{23/3}_{23/3}.\label{disc_very_strong_rho} 
\end{align}

\subsection{Renormalised internal energy}

The internal energy balance was derived in (\ref{disc_IE}) and is repeated here for convenience,
\begin{align}
  D_t\frac{p_{ijk}}{\gamma-1}+(pD\vv)_{ijk}-\Df_{ijk}
\nonumber\\
  +D_-^x(\frac{\overline{\rho u}_{i+1/2jk}}{2(\gamma-1)\hat \beta_{i+1/2jk}})  
+D_-^y(\frac{\overline{\rho v}_{ij+1/2k}}{2(\gamma-1)\hat \beta_{ij+1/2k}})  
+  D_-^y(\frac{\overline{\rho w}_{ijk+1/2}}{2(\gamma-1)\hat \beta_{ijk+1/2}})&=  \label{disc_RIE}\\
D_-^x(\tilde \nu_{i+1/2jk}\frac{(\pf_x)_{i+1/2jk}}{\gamma-1})+
D_-^y(\tilde \nu_{ij+1/2k}\frac{(\pf_y)_{ij+1/2k}}{\gamma-1})+\nonumber \\
D_-^z(\tilde \nu_{ijk+1/2}\frac{(\pf_z)_{ijk+1/2}}{\gamma-1})\nonumber \\
+\kappa_r\left(D_-^xD^x_+T^4_{ijk}+D_-^yD^y_+T^4_{ijk}+D_-^z\kappa_rD^z_+T^4_{ijk}\right)&.\nonumber
\end{align}
As in the continuous case, we multiply by $H'(T_i)=H'_i$ and sum in space. Since $T\geq 0$, we have 
\begin{align}
H(T)=(1+T)^{1-\omega}<1+T, \,\,\omega>0,\nonumber\\
0<H'(T)<1,\quad
H''(T)<0,\quad
|H''|<c<1.\nonumber
\end{align}

\vspace{0.5cm}
We carry out the calculations for a few of the terms separately and begin with the temporal term of (\ref{disc_RIE}):
\begin{align}
H'(T_i)D_t\frac{p_{ijk}}{\gamma-1}= c_v(\rho_{ijk} H_{ijk})_t + c_v(H'_{ijk}T_{ijk}-H_{ijk})\partial_t\rho_{ijk}.\nonumber
\end{align}
Next, we multiply by $V_i$ and sum. Using the continuity equation and the boundary conditions/fluxes, we obtain
\begin{align}
 \sum_{ijk=0}^N V_{ijk}  H'(T_{ijk})D_t\frac{p_{ijk}}{\gamma-1}&=\nonumber\\
 \sum_{ijk=0}^N V_{ijk}  c_v(\rho_{ijk} H_{ijk})_t + c_v(H'_{ijk}T_{ijk}-H_{ijk})(-\nabla_-\cdot (\fb^{c,1}_{i+1/2jk},\gb^{c,1}_{ij+1/2jk},\hb^{c,1}_{ijk+1/2})\nonumber \\
 + \sum_{ijk=0}^N V_{ijk}c_v(H'_{ijk}T_{ijk}-H_{ijk})(\nabla_-\cdot (\fb^{d,1}_{i+1/2jk},\gb^{d,1}_{ij+1/2jk},\hb^{d,1}_{ijk+1/2})&,\label{rt1}
\end{align}
where $\nabla_{-/+}=(D_{-/+}^x,D_{-/+}^y,D_{-/+}^z)$.

\vspace{0.5cm}

We adopt the short-hand notation $\fb^T_{i+1/2jk}=\frac{\overline{\rho u}_{i+1/2jk}}{2(\gamma-1)\hat \beta_{i+1/2jk}}$ etc. for the discrete thermal fluxes and define the corresponding discrete renormalised thermal flux as, $\fb^{RT}_{i+1/2jk}=c_v\frac{H(T_{i+1jk})+H(T_{ijk})}{2}\overline{\rho u}_{i+1/2jk}$.
\begin{remark}
As above, we adopt the convention that values ``outside'' the domain  are identical with the boundary value. That is, $\fb^{RT}_{-1/2jk}=c_v\frac{H(T_{0jk})+H(T_{-1jk})}{2}\overline{\rho u}_{-1/2jk}= c_vH(T_0)(\rho u)_{0jk}=0$.
\end{remark}  
Multiplying the first convective thermal flux of (\ref{disc_RIE}) by $V_{ijk}H'(T_{ijk})$, yields
\begin{align}
 V_{ijk}H'(T_{ijk})D_-^x(\frac{\overline{\rho u}_{i+1/2jk}}{2(\gamma-1)\hat \beta_{i+1/2jk}})=V_{ijk}H'(T_{ijk})D_-^x(\fb^T_{i+1/2jk})&=\nonumber\\
V_{ijk}D_-^x \fb^{RT}_{i+1/2jk}\nonumber\\
 -V_{ijk}c_vD_-^x\frac{H(T_{i+1jk})+H(T_{ijk})}{2}\overline{\rho u}_{i+1/2jk}
 + V_{ijk}H'(T_{ijk})\frac{(\fb^{T}_{i+1/2jk}-\fb^T_{i-1/2jk})}{\Delta_-^xx_{i+1/2}}
  &= \nonumber\\
V_{ijk}D_-^x \fb^{RT}_{i+1/2jk}  +E^x_{ijk},\label{rt2}
\end{align}
where
\begin{align}
E^x_{ijk} = - \frac{V_{ijk}}{\Delta_-^xx_{i+1/2}} \left(c_v\frac{H(T_{i+1jk})+H(T_{ijk})}{2}-\frac{H'(T_{ijk})}{2(\gamma-1)\hat \beta_{i+1/2jk}}\right)\overline{\rho u}_{i+1/2jk}\label{Ex}\\
  +\frac{V_{ijk}}{\Delta_-^xx_{i+1/2}}\left(c_v\frac{H(T_{ijk})+H(T_{i-1jk})}{2}-\frac{H'(T_{ijk})}{2(\gamma-1)\hat \beta_{i-1/2jk}}\right)\overline{\rho u}_{i-1/2jk}.\nonumber
\end{align}
We carry out the same operations for the thermal fluxes in the $y,z$-directions resulting in the analogous definition of $E^{y,z}_{ijk}$.

Next, we multiply (\ref{disc_RIE}) by $V_{ijk}H'_{ijk}$, sum in space and use (\ref{rt1}) and (\ref{rt2}) (along with the similar expressions obtained in the other two directions) to obtain
\begin{align}
 \sum_{ijk=0}^N V_{ijk} \left( c_v(\rho_{ijk} H_{ijk})_t + D_-^x\fb^{RT}_{i+1/2jk}+D_-^x\gb^{RT}_{ij+1/2k}+D_-^z\hb^{RT}_{ijk+1/2}\right) +\sum_{ijk=0}^N(E^x_{ijk}+E^y_{ijk}+E^z_{ijk}) \nonumber \\
 + \sum_{ijk=0}^{N} V_{ijk}c_v(H'_{ijk}T_{ijk}-H_{ijk})\nabla_-\cdot (\fb^{d,1}_{i+1/2jk},\gb^{d,1}_{ij+1/2jk},\hb^{d,1}_{ijk+1/2})\nonumber\\
 - \sum_{ijk=0}^{N} V_{ijk}c_v(H'_{ijk}T_{ijk}-H_{ijk})\nabla_-\cdot (\fb^{c,1}_{i+1/2jk},\gb^{c,1}_{ij+1/2jk},\hb^{c,1}_{ijk+1/2})\nonumber\\
  +\sum_{ijk=0}^N V_{ijk}H'(T_{ijk})((pD\vv)_{ijk}-\Df_{ijk})
&=  \label{disc_RIE2}\\
\sum_{ijk=0}^N V_{ijk}H'(T_{ijk})\nabla_-\cdot \left(\tilde \nu_{i+1/2jk}\frac{(\pf_x)_{i+1/2jk}}{\gamma-1},\tilde \nu_{ij+1/2k}\frac{(\pf_y)_{ij+1/2k}}{\gamma-1},\tilde \nu_{ijk+1/2}\frac{(\pf_z)_{ijk+1/2}}{\gamma-1}\right)\nonumber \\
+\sum_{ijk=0}^N V_{ijk}H'(T_{ijk})\kappa_r\left(D_-^xD^x_+T^4_{ijk}+D_-^yD^y_+T^4_{ijk}+D_-^zD^z_+T^4_{ijk}\right).\nonumber
\end{align}
Combining the error terms $(E^{x,y,z}_{ijk})$ and the convective ``c,1'' terms lead to,
\begin{align}
 \sum_{ijk=0}^N V_{ijk} \left( c_v(\rho_{ijk} H_{ijk})_t + D_-^x\fb^{RT}_{i+1/2jk}+D_-^x\gb^{RT}_{ij+1/2k}+D_-^z\hb^{RT}_{ijk+1/2}\right) +\sum_{ijk=0}^N(\Ec^x_{ijk}+\Ec^y_{ijk}+\Ec^z_{ijk}) \nonumber \\
 + \sum_{ijk=0}^{N} V_{ijk}c_v(H'_{ijk}T_{ijk}-H_{ijk})\nabla_-\cdot (\fb^{d,1}_{i+1/2jk},\gb^{d,1}_{ij+1/2jk},\hb^{d,1}_{ijk+1/2})\nonumber\\
  +\sum_{ijk=0}^N V_{ijk}H'(T_{ijk})((pD\vv)_{ijk}-\Df_{ijk})
&=  \label{disc_RIE3}\\
\sum_{ijk=0}^N V_{ijk}H'(T_{ijk})\nabla_-\cdot \left(\tilde \nu_{i+1/2jk}\frac{(\pf_x)_{i+1/2jk}}{\gamma-1},\tilde \nu_{ij+1/2k}\frac{(\pf_y)_{ij+1/2k}}{\gamma-1}),\tilde \nu_{ijk+1/2}\frac{(\pf_z)_{ijk+1/2}}{\gamma-1}\right)\nonumber \\
+\sum_{ijk=0}^N V_{ijk}H'(T_{ijk})\kappa_r\left(D_-^xD^x_+T^4_{ijk}+D_-^yD^y_+T^4_{ijk}+D_-^zD^z_+T^4_{ijk}\right),\nonumber
\end{align}
where the error terms are given by
\begin{align}
\Ec^x_{ijk} = - \frac{c_vV_{ijk}}{\Delta_-^xx_{i+1/2}} \left(\frac{H(T_{i+1jk})+H(T_{ijk})}{2}-\frac{H'(T_{ijk})}{\widehat {(T^{-1})}_{i+1/2jk}}-(H(T_{ijk})-H'(T_{ijk})T_{ijk})\right)\overline{\rho u}_{i+1/2jk}\label{Ex}\\
  +\frac{c_vV_{ijk}}{\Delta_-^xx_{i+1/2}}\left(\frac{H(T_{ijk})+H(T_{i-1jk})}{2}-\frac{H'(T_{ijk})}{\widehat{ T^{-1}}_{i-1/2jk}}-(H(T_{ijk})-H'(T_{ijk})T_{ijk})\right)\overline{\rho u}_{i-1/2jk},\nonumber
\end{align}
and analogously for $\Ec^{y,z}_{ijk}$. Since $H'$ is bounded, we have that
\begin{align}
\frac{1}{\Delta_-^xx_{i+1/2}}  |\frac{H(T_{i+1jk})+H(T_{ijk})}{2}-H(T_{ijk})|\leq C |D^+_xT_{ijk}|.\nonumber
\end{align}
By (\ref{mean_rels}), $\frac{1}{\Delta_-^xx_{i+1/2}}|\frac{1}{\widehat {T^{-1}}_{i+1/2jk}}-T_{ijk}|\leq |D^+_x T_{ijk}|$. We conclude that,
\begin{align}
  |\sum_{ijk=0}^N\Ec^x_{ijk} | \leq \C\sum_{ijk=0}^N|(D_x^+T_{ijk})\overline{\rho u}_{i+1/jk}|\lesssim \|D_+T\|_2\|\overline{\rho u}\|_2. \nonumber
\end{align}
Since $D^+_x T_{ijk}\in L^2(H^1)$ by (\ref{rlogr_d}) and (\ref{disc_temp_bound}), and $\overline{\rho u}_{i+1/2jk}\in L^2_2$ by (\ref{disc_very_strong_rho}) and (\ref{disc_v_H1}), we conclude that the error term $|\sum \Ec^x_{ijk}|$ is bounded. By repeating the arguments in the $y,z$-directions, the remaining two error terms are also bounded.

Next, we consider the terms
\begin{align}
\sum_{ijk=0}^N V_{ijk}H'(T_{ijk})((pD\vv)_{ijk}-\Df_{ijk}),\nonumber
\end{align}
in (\ref{disc_RIE3}). Noting that $H'$ is bounded, the pressure term is controlled just like in (\ref{KE_est1}) and the viscous terms are bounded by the estimate obtained by (\ref{KE_est1}).

We must also control the diffusive terms in (\ref{disc_RIE3}). As usual, we present the arguments for the $x$-term. Namely, we use
\begin{align}
  (\pf_x)_{i+1/2jk} &= \frac{1}{2\hat \beta_{i+1/2jk}}D_+^x\rho_{ijk}+\frac{\bar\rho_{i+1/2jk} }{2}D_+^x\frac{1}{\beta_{ijk}},\nonumber
\end{align}
in the diffusive terms,
\begin{align}
  \sum_{ijk=0}^N V_{ijk}\left(H'(T_{ijk})D_-^x\left(\tilde \nu_{i+1/2jk}\frac{(\pf_x)_{i+1/2jk}}{\gamma-1}\right)-c_v(H'_{ijk}T_{ijk}-H_{ijk})D_-^x \fb^{d,1}_{i+1/2jk}\right)&=\nonumber\\
-  \sum_{ijk=0}^{N-1,N,N} V_{ijk}(D_+^xH'(T_{ijk}))\left(\tilde \nu_{i+1/2jk}\frac{1}{\gamma-1}\left(\frac{1}{2\hat \beta_{i+1/2jk}}D_+^x\rho_{ijk}+\frac{\bar\rho_{i+1/2jk} }{2}D_+^x\frac{1}{\beta_{ijk}}\right)\right)&\nonumber\\
  + \sum_{ijk=0}^{N-1,N,N} V_{ijk}c_v(D_+^x(H'_{ijk}T_{ijk}-H_{ijk}))\fb^{d,1}_{i+1/2jk} &=\nonumber\\
-  \sum_{ijk=0}^{N-1,N,N} V_{ijk}(D_+^xH'(T_{ijk}))\left(\tilde \nu_{i+1/2jk}\frac{1}{\gamma-1}\left(\frac{\bar\rho_{i+1/2jk} }{2}D_+^x\frac{1}{\beta_{ijk}}\right)\right)&\label{diff_rt}\\
   + \sum_{ijk=0}^{N-1,N,N} V_{ijk}c_v\left((D_+^x(H'_{ijk}T_{ijk}-H_{ijk}))-\frac{D_+^xH'(T_{ijk})}{c_v(\gamma-1)2\hat \beta_{i+1/2jk}} \right)\fb^{d,1}_{i+1/2jk} &,\nonumber
\end{align}
where we have summed by parts and used the boundary conditions/fluxes.

Using, 
\begin{align}
  \frac{1}{\hat \beta_{i+1/2}} = \frac{\Delta_+ \log\beta_i}{\Delta \beta_i}=2R\frac{-\Delta_+ \log T_i}{\Delta_+ T_{i}^{-1}}=2RT_{i}T_{i+1}\frac{\Delta_+ \log T_i}{\Delta_+ T_{i}}=2R\frac{\check T_{i+1/2}^2}{\hat T_{i+1/2}}\nonumber
\end{align}
and $c_v(\gamma-1)=c_v(\frac{c_p}{c_v}-1)=R$ the last term of (\ref{diff_rt}) can be rewritten as,
\begin{align}
  \frac{D_+^xH'(T_i)}{c_v(\gamma-1)2\hat \beta_{i+1/2jk}}\fb^{d,1}_{i+1/2jk}=
  \frac{D_+^xH'(T_i)}{2R\hat \beta_{i+1/2jk}}\fb^{d,1}_{i+1/2jk}&=\nonumber\\
  (D_+^xH'(T_i))\frac{\check T_{i+1/2}^2}{\hat T_{i+1/2}} \fb^{d,1}_{i+1/2jk}&.\nonumber
\end{align}
This allows us to rewrite the error term in (\ref{diff_rt}), i.e., the last sum (denoted $E$). We make the following manipulations.
\begin{align}
E\leq  c_v\left|\frac{\check T_{i+1/2}^2}{\hat T_{i+1/2}}D_+^xH'(T_{ijk})-D_+^x(H'_{ijk}T_{ijk}-H_{ijk})\right|| \fb^{d,1}_{i+1/2jk} |&=\nonumber\\
 c_v \left|\frac{\check T_{i+1/2}^2}{\hat T_{i+1/2}}D_+^xH'_{ijk}-\overline{H'}_{i+1/2jk}D_+^xT_{ijk}-\bar T_{i+1/2jk}D_+^xH'_{ijk}+D_+^xH_{ijk}\right|| \fb^{d,1}_{i+1/2jk} |&=\nonumber\\
 c_v \left|\left(\frac{\check T_{i+1/2}^2}{\hat T_{i+1/2}}-\bar T_{i+1/2jk}\right)D_+^xH'_{ijk}-(\overline{H'}_{i+1/2jk}-D^TH_{i+1/2jk})D_+^xT_{ijk}\right|| \fb^{d,1}_{i+1/2jk} |&,\nonumber 
\end{align}
where $D^TH_{i+1/2jk}=\frac{H_{i+1jk}-H_{ijk}}{T_{i+1jk}-T_{ijk}}$ denotes the finite difference w.r.t. temperature. By the mean value theorem, we have $D^TH_{i+1/2jk}=H'(\theta)$ for some $\theta\in [T_i,T_{i+1}]$ and $|H'(T)|\leq 1$ for all $T\geq 0$. Furthermore,
\begin{align}
(\frac{\check T_{i+1/2}^2}{\hat T_{i+1/2}}-\bar T_{i+1/2jk})D_+^xH'_{ijk}&=\nonumber
\frac{1}{\Delta_+x_{i-1/2}}(\frac{\check T_{i+1/2}^2}{\hat T_{i+1/2}}-\bar T_{i+1/2jk})\Delta_+^xH'_{ijk},\nonumber
\end{align}
and $|\frac{1}{\Delta_+x_{i-1/2}}(\frac{\check T_{i+1/2}^2}{\hat T_{i+1/2}}-\bar T_{i+1/2jk})|\lesssim |D_+^xT_{ijk}|$. Hence, we conclude that all temperature dependent factors appearing in $E$ scale as $H'|D_+^xT_{ijk}|$ where $H'$ is bounded.
Hence,
\begin{align}
E\lesssim |D_+^xT_{ijk}||\fb^{d,1}_{i+1/2jk}|,\nonumber
\end{align}
and to control the  last term in (\ref{diff_rt}), we need $\nabla_+T_{ijk}\in L^2_2$ and the fluxes $\fb^{d,1}, \gb^{d,1}, \hb^{d,1}$ all in $L^2_2$. The temperature bound is obtained from (\ref{disc_temp_bound}). The diffusive flux (in the x-direction) is given by
\begin{align}
\fb^{d,1}_{i+1/2jk}=(\frac{\mu_0}{\hat \rho_{i+1/2jk}}+\mu_1\bar \rho_{i+1/2jk}+(\Delta_+^xx_{ijk})\lambda_{i+1/2jk})D_+^x\rho_{ijk}\nonumber.
\end{align}
The first term is, $\mu_0D_+^x\log\rho_{ijk}$ which bounded in $L^2_2$ by (\ref{disc_rho_H1}). The second is in $L^2_2$ by (\ref{disc_strong_rho}). The third is proportional to $(R^*|\bar u|+\frac{\Delta u}{4})\Delta_+^x\rho_{ijk}$, where $R^*=\max(\frac{1}{2},\Delta \log\rho)$. Hence, we must control terms that are majorised by: $|u|\rho$, $\log\rho |u|\rho$, and $\rho\Delta u$. By noting that $\log \rho$ is controlled in any $L^p_p$ space by $L^2_2$ estimates of $\rho$ and $\rho^{-1}$ (c.f. (\ref{logr_Lpp})), the estimates (\ref{disc_v_H1}), (\ref{disc_very_strong_rho}) and (\ref{svolume_disc}) provide the necessary control. (As usual, the other two directions are handled analogously.)

We have now handled the last term of (\ref{diff_rt}) and we turn to the second last,
\begin{align}
-\sum_{ijk=0}^{N-1,N,N} V_{ijk}(D_+^xH'(T_{ijk}))\left(\tilde \nu_{i+1/2jk}\frac{1}{\gamma-1}\left(\frac{\bar\rho_{i+1/2jk} }{2}D_+^x\frac{1}{\beta_{ijk}}\right)\right)\label{fourier_diff}.
\end{align}
Using the mean-value theorem again, we have
\begin{align}
D_+^xH'(T_{ijk})=H''(\theta_{i+1/2jk})D_+^xT_{ijk}, \quad \theta_{i+1/2jk}\in [T_{ijk},T_{i+1jk}].\nonumber
\end{align}
Furthermore, $D_+^x\frac{1}{\beta_{ijk}}=2RD_+^xT_{ijk}$ such that (\ref{fourier_diff}) is a negative definite term. (Hence, it makes a negative contribution to the right-hand side of (\ref{disc_RIE3}), which is what we wish.)

We now have all estimates, to repeat the argument made for the continuous a priori estimate. Hence, we integrate (\ref{disc_RIE3}) in time and arrive at the desired bound
\begin{align}
  \C\geq -\int_0^\T\sum_{ijk=0}^{N-1,N,N} V_{ijk}\kappa_r D_+^xH'_{i+1/2jk} D_+^xT^4_{ijk}\,dt \label{disc_temp_est}.
\end{align}
Our goal is to mimic the continuous estimate. That is,  we wish to demonstrate that 
\begin{align}
  D_+T_{ijk}^{2-\omega}\in L^2_2.\nonumber
\end{align}
follows from (\ref{disc_temp_est}).

Dropping the recurring $jk$ indices, we use
\begin{align}
D_+T_i^4=(T_{i+1}^3+T_i^3+T_i^2T_{i+1}+T_{i+1}^2T_i)D_+T_i,\nonumber
\end{align}
to rewrite
\begin{align}
  -D_+H'_iD_iT_i^4=
  -\frac{\Delta_+H'_i}{\Delta T_i}(T_{i+1}^3+T_i^3+T_i^2T_{i+1}+T_{i+1}^2T_i)(D_+T_i)^2.
  \label{hdt4}
\end{align}
By the mean-value theorem, we have
\begin{align}
-\frac{\Delta_+H'_i}{\Delta T_i}=\omega(1-\omega)(1+\theta_{i+1/2})^{-\omega-1}>0,\quad \theta_{i+1/2}\in [T_i,T_{i+1}].\nonumber
\end{align}
We recast (\ref{hdt4}) as
\begin{align}
  -D_+H'_iD_iT_i^4=
  \omega(1-\omega)\frac{(T_{i+1}^3+T_i^3+T_i^2T_{i+1}+T_{i+1}^2T_i)}{(1+\theta_{i+1/2})^{\omega+1}}(D_+T_i)^2.
  \label{hdt5}
\end{align}
Without restriction, assume that $T_i\geq T_{i+1}\geq 0$. Then
\begin{align}
(1+\theta_{i+1/2})^{\omega+1}\leq (1+T_i)^{\omega+1}.\nonumber
\end{align}
Hence,
\begin{align}
 \frac{T_i^3}{(1+\theta_{i+1/2})^{\omega+1}}\geq \frac{T_i^3}{(1+T_i)^{\omega+1}}.\nonumber
\end{align}
Noting that all terms in (\ref{hdt5}) are positive, we can extract the bound
\begin{align}
T_i^{2-\omega}(D_+T_i)^2=(T_i^{1-\omega/2}D_+T_i)^2.\nonumber
\end{align}
By the mean-value theorem,
\begin{align}
D_+T_i^{2-\omega/2}=(2-\omega/2)\theta^{1-\omega/2}D_+T_i,
\end{align}
for some $\theta \in [T_i,T_{i+1}]$. Since $T_i\geq \theta$, $T_i^{1-\omega/2}D_+T_i$ bounds $D_+T_i^{2-\omega/2}$.
Finally, we note that the choice $T_{i+1}<T_i$ was arbitrary and we can reverse the argument. 
Hence, we have
\begin{align}
D_+T_i^{2-\omega/2}\in L^2_2\nonumber
\end{align}
and, as in the continuous case, Sobolev embeddings and an interpolation inequality ensures that
\begin{align}
T_i\in L^{4+\delta}_{4+\delta},\label{disc_strong_T}
\end{align}
 for some $\delta>0$.

\subsection{Improved estimates}\label{sec:disc_improved}

The derivations of the improved estimates are the same as in Section \ref{sec:disc_improved}. We summarise the ones that are necessary to prove convergence:
\begin{align}
  p_{ijk} \in L^{2+\delta}_{2+\delta},\quad \delta>0,\label{imp_p}\\
  \vv_{ijk} \in L^\infty_1, \quad \vv\in L^{10/3}_2\cap L^{10/4}_3,\label{imp_v}\\
  \quad \rho_{ijk}|\vv_{ijk}|\in L^{2+\delta}_{2+\delta},\quad \delta>0,\label{imp_mom}\\
  (\rho(\vv_l)\vv_m)_{ijk} \in L^{1+\delta}_{1+\delta}, \quad \delta>0,\quad l,m=\{1,2,3\},\label{imp_ke}\\
  (\rho|\vv|^2\vv)_{ijk} \in L^{1+\delta}_{1+\delta},\quad \delta>0. \label{imp_ke_flux}
\end{align}
We also add the following estimate that is a consequence of the bounds (\ref{svolume_disc}) and (\ref{disc_very_strong_rho}):
\begin{align}
\log \rho \in L^p_p,\quad \textrm{for any}\quad  0<p<\infty\label{logr_Lp}.
\end{align}

\section{Convergence to a weak solution}\label{sec:ex}

We begin by introducing some notation. First, we denote discrete grid functions with a subscript $h$, ie., $\rho_h$ is a field whose $ijk$th component is $(\rho_{ijk})$. A solution obtained from the numerical scheme is hence denoted $\uu(t)_h=(\rho(t)_h,\mm(t)_h,E(t)_h)$, where $(\mm_h)_{ijk}=\rho_{ijk}(u_{ijk},v_{ijk},w_{ijk})$. We will use the convention that the fields are multiplied component-wise such that e.g. $(m_1)_h=\rho_hu_h$ and $\mm_h=\rho_h\vv_h$. In this notation, $E_h=\frac{1}{2}\rho_h|\vv_h|^2+\frac{p_h}{\gamma-1}$ with components $(E_h)_{ijk}=\frac{1}{2}\rho_{ijk}|\vv_{ijk}|^2+\frac{p_{ijk}}{\gamma-1}$. Furthermore, differences will be collected into fields using the same the short-hand notation, e.g. $D_+^x\rho_h$ and $\lambda_hD_+^x\rho_h$. Finally, it is often necessary to look at the precise structure of terms and then we may say that e.g. $a_{i+1/2jk}b_{ijk}$ is bounded in some space, by that we mean that the field $a_hb_h$ is bounded in said space.

\subsection{Solvability of the semi-discrete scheme}

The first step towards demonstrating the existence of a weak solution, is to establish that the semi-discrete system (\ref{scheme}) can be solved in time. Given the strong a priori estimates and that $\rho,T>0$  (as have been demonstrated above), we can apply  Carath{\'e}odory's existence theorem, to conclude that there exists a unique and continuous solution to the semi-discrete system (\ref{scheme}), for $0\leq t\leq \T$ and for any $h>0$. Hence, we can generate a sequence of triplets $\{\rho_h,\vv_h,T_h\}$.

The next step is to show that this triplet converges (up to a subsequence) to a weak solution in the sense of Def. \ref{def:weak} as $h\rightarrow 0$. To that end, we proceed as in the preliminary discussion on compactness in Section \ref{sec:compact}. A difference in the semi-discrete case, is that we must also show that the approximate numerical fluxes are consistent with the mathematical fluxes, i.e., that the error terms vanish. (In the subsequent analysis, we will tacitly draw subsequences where necessary.)

\subsection{Convergence a.e.}

First, we apply Aubin-Lions lemma (Lemma \ref{lemma:aubin}) to the continuity equation. The arguments are initially the same as in Section \ref{sec:compact}, partial sums replacing partial integration. Hence, we begin to verify that the numerical fluxes of the continuity equation reside in $L^1_1$. 

The total flux in the x-direction is given by (\ref{fc}) and (\ref{flux_d1}):
\begin{align}
  \fb^{1}_{i+1/2jk}&=\frac{(\rho u)_{i+1jk}+(\rho u)_{ijk}}{2}-\tilde\nu_{i+1/2jk}D_+^x \rho_{ijk}\label{total_flux_rho}\\
  \tilde \nu_{i+1/2jk}&=\frac{\mu_0}{\hat\rho_{i+1/2jk}}+\mu_1\frac{\rho_{i+1jk}+\rho_{ijk}}{2}+(\Delta_+^xx)_{i+1/2}\lambda_{i+1/2jk}\nonumber\\
    \lambda_{i+1/2jk}&\sim |\Delta_+^x u_{ijk}| + |\bar u_{i+1/2jk}||\Delta_+^x\log\rho_{ijk}|\nonumber.
\end{align}
The boundedness of the first two flux terms follows (\ref{imp_mom}).

Next, we consider the artificial diffusion term. $\Delta^x_+ u_h$ is bounded (and vanishes) by (\ref{disc_v_H1}) in $L^2_2$. $|\Delta_+^x \log\rho_h||\bar u_h|$ is bounded in (at least) $L^2_2$ by (\ref{imp_v}) and by $\log\rho_h$ being bounded in an arbitrary $L^p_p$ space (see (\ref{logr_Lp})).  Using the bound (\ref{disc_rho_H1}) on the density gradient, we have
\begin{align}
  \lim_{h\rightarrow 0}\int_0^\T\|(\Delta_+^xx_{i+1/2})\lambda_{i+1/2jk}D_+^x\rho_{ijk}\|_1\,dt &\leq\nonumber\\
  \lim_{h\rightarrow 0} h\int_0^\T(\|\lambda_{i+1/2jk}\|_2^2+\|D_+^x\rho_{ijk}\|^2_2)\,dt&\rightarrow 0.\nonumber
\end{align}
Likewise, $\frac{1}{\hat\rho_{i+1/2jk}}D_+^x\rho_{ijk}=D_+^x\log\rho_{ijk}$, which is 
 bounded in $L^2_2$ by e.g. (\ref{disc_rho_H1}), and 
$\left(\mu_1\frac{\rho_{i+1}+\rho_{i}}{2}\right)D_+^x\rho_{ijk}$ by (\ref{disc_strong_rho}).
Hence, we may conclude, as in Section \ref{sec:compact} using Aubin-Lions Lemma, that $\rho_{ijk}$ converges a.e. in the limit $h\rightarrow 0$.

Using the analogous arguments as in Section \ref{sec:compact}, we obtain that $\rho_i,\vv_i,p_i,T_i$ all converges a.e. Likewise, we can use the arguments of Section \ref{sec:compact} to conclude that the magnitude of the gradients converge a.e. and in particular that,
\begin{align}
 \{\, |\Delta_+^{x,y,z}\rho_{ijk}|,\quad
  |\Delta_+^{x,y,z} u_{ijk}|,\quad
  |\Delta_+^{x,y,z} T_{ijk}|\,\}\rightarrow 0\quad \textrm{a.e.}.\label{diff_vanish}
\end{align}

\subsection{Weak convergence of fluxes}\label{sec:disc_compact}

Before considering the numerical fluxes, we briefly explain the techniques we use and why convergence is a little more delicate than in Section \ref{sec:weak_seq}. First  $(m_1)_h=\rho_hu_h\rightarrow \rho u$ strongly in $L^2_2$ by the same arguments as laid out in Section \ref{sec:weak_seq}. (We remind that $\rho_hu_h$ is the pointwise multiplication of the two fields, i.e. $\rho_{ijk}u_{ijk}$.)

This also ensures that e.g. $\overline{\rho u}_{i+1/2jk}=\frac{\rho_{i+1jk}u_{i+1jk}+\rho_{ijk}u_{ijk}}{2}$ converges to the correct limit in $L^2_2$. However, it does not ensure that e.g. $\bar \rho_{i+1/2jk}\bar u_{i+1/2jk}$ converges to the same limit since,
\begin{align}
\bar \rho_{i+1/2jk}\bar u_{i+1/2jk}=\overline{\rho u}_{i+1/2jk}-\frac{\Delta_+^x \rho_{ijk}\Delta_+^xu_{ijk}}{4}.\nonumber
\end{align}
Here, $\overline{\rho u}_{i+1/2jk}$ yields the desired limit and the last term is an error term that must vanish in $L^1_1$. We can see that it does, since $D_+^x\rho_{ijk}$  reside in $L^2_2$ it follows that $\Delta_+^x \rho_{ijk}\sim h D_+^x\rho_{ijk}$ vanishes in $L^2_2$ and the same is true for $D_+^xu_{ijk}$.

The error term can also be handled in another way: From the estimates of $\rho_h$ and $u_h$, we can immediately conclude that $\bar \rho_h\bar u_h$ is bounded in $L^{1+\epsilon}_{1+\epsilon}$, i.e., the sequence is equi-integrable. Since that is also true for $(\rho u)_h$, we conclude that the error term has to reside in $L^{1+\epsilon}_{1+\epsilon}$. Then by a.e. convergence of the magnitude of the gradients, using (\ref{diff_vanish}) it follows that that $(|\Delta_+^x \rho|)_h$ as well as $(|\Delta_+^x u|_h$ vanish a.e. (Note that the last approach does not rely on an explicit bound of the error term via gradient estimates.)

We will employ both methods to demonstrate consistency of the numerical fluxes.

\vspace{0.5cm}

{\bf Continuity equation:} We need to show that $\fb^1_h$ in (\ref{total_flux_rho}) converges weakly in $L^{1+\epsilon}_{1+\epsilon}$, $\epsilon>0$. (As usual, the $yz$-fluxes are handled in the same way and omitted.)

As discussed above,  $\overline{\rho u}_{i+1/2jk}$ converges strongly in (at least) $L^2_2$, by (\ref{imp_mom}) and a.e. convergence of $\rho_i$ and $u_i$. The remaining terms are handled as in Section \ref{sec:weak_seq}: First, we use a.e. convergence of $\rho$ and the estimates in (\ref{svolume_disc}), to deduce strong convergence of $\rho^{-1}_h$ in $L^{2+\delta}_{2+\delta}$ for some $\delta>0$. Then we consider the diffusive $x$-flux:
\begin{align}
  \tilde \nu_{i+1/2jk}D_+^x\rho_{ijk}&=\nonumber\\
  \left((\Delta_+^xx_{ijk})(R^*|\bar u_{i+1/2jk}|+\frac{|\Delta u|}{4})+\frac{\mu_0}{\hat\rho_{i+1/2jk}}+\mu_1\bar \rho_{i+1/2jk}\right)D_+^x\rho_{ijk}.\nonumber
\end{align}
These terms are proportional to: $|\Delta\log \rho||\bar u|\Delta\rho, |\Delta u|\Delta \rho,  D_+\rho^2 $ and $D_+\log\rho$, which are controlled in $L^{1+\epsilon}_{1+\epsilon}$ by (\ref{disc_very_strong_rho}), (\ref{disc_rho_H1}), (\ref{disc_v_H1}), (\ref{logr_Lp}) and (\ref{imp_v}). Weak convergence of the two diffusion terms in $L^{1+\epsilon}_{1+\epsilon}$ is immediate and the artificial diffusion terms vanish since $|\Delta \rho|\rightarrow 0$ a.e. by  (\ref{diff_vanish}) (and the remaining factors are functions of the primary variables that converge a.e.).

\vspace{0.5cm}

{\bf Momentum equations:} 
We consider a few typical terms and the remaining can be handled in the same way. We show that all terms are bounded in $L^{1+\epsilon}_{1+\epsilon}$, $\epsilon>0$, and convergent a.e. to the correct limit.

Convergence of the temporal term follows immediately from the strong convergence of $\rho_h\vv_h$ in $L^2_2$. (See \eqref{imp_mom}.)

The inviscid velocity terms are all of the form:
\begin{align}
  \bar u_{i+1/2jk}\overline{(\rho v)}_{i+1/2jk}= \overline{u\rho v}_{i+1/2jk}-\frac{1}{4}\Delta_+^x u_{ijk}\Delta_+^x(\rho v)_{ijk}.\nonumber
\end{align}
$\overline{u\rho v}_{i+1/2jk}$ is equi-integrable thanks to  $\rho_h \vv_h\in L^2_2$ (c.f. (\ref{imp_mom})) and $\vv_h\in L^{5/2}_{5/2}$. Pointwise convergence of $\rho$ and $\vv$ ensure the correct limit. Noting that $\bar u_{i+1/2jk}\overline{(\rho v)}_{i+1/2jk}$ is bounded in the same space by the same estimates, we conclude immediately that the error term $\frac{1}{4}\Delta_+^x u_{ijk}\Delta_+^x(\rho v)_{ijk}$ is bounded in $L^{1+\delta}_{1+\delta}$. Furthermore, by (\ref{diff_vanish}) both differences approach zero a.e. Hence, we conclude that these terms converge weakly in $L^{1+\epsilon}_{1+\epsilon}$ for $0<\epsilon<\delta$.

Next, we consider the pressure term (suppressing the $jk$ indices):
\begin{align}
  \frac{1}{R}\frac{\bar \rho_{i+1/2}}{2\bar \beta_{i+1/2}}=
  \frac{(\rho_{i+1}+\rho_i)}{2}\frac{2}{\frac{1}{T_{i+1}}+\frac{1}{T_i}}&=\nonumber\\
  \frac{\rho_{i+1}T_{i+1} + \rho_iT_i}{2}-
  \rho_i(\frac{T_i}{2}-\frac{1}{\frac{1}{T_i}+\frac{1}{T_{i+1}}})
   -\rho_{i+1}(\frac{T_{i+1}}{2}-\frac{1}{\frac{1}{T_i}+\frac{1}{T_{i+1}}}).\nonumber
\end{align}
The $\rho_iT_i$ terms converge to $p_i$ in $L^2_2$ thanks to a.e. convergence and (\ref{imp_p}). Moreover, $\bar \rho_{i+1/2}$ is controlled by (\ref{disc_very_strong_rho}) and $\frac{1}{2\bar \beta_{i+1/2}}\sim \frac{T_iT_{i+1}}{T_i+T_{i+1}}\lesssim T_{i+1}+T_i$ is controlled by (\ref{disc_strong_T}). Hence, $\frac{\bar \rho_{i+1/2}}{2\bar \beta_{i+1/2}}$ is bounded in (a better space than) $L^2_2$. Therefore, we conclude that the error terms are bounded in $L^2_2$.

Consider the first error term:
\begin{align}
  \rho_i(\frac{T_i}{2}-(\frac{T_{i+1}T_i}{T_i+T_{i+1}})=  \rho_iT_i(\frac{1}{2}-(\frac{T_{i+1}}{T_i+T_{i+1}})&=\nonumber \\
  \rho_iT_i(\frac{T_i+T_{i+1}-2T_{i+1}}{2(T_i+T_{i+1})})=  \rho_iT_i(\frac{T_i-T_{i+1}}{2(T_i+T_{i+1})})&\nonumber
\end{align}
In the last expression, we note that $\rho_i,T_i,T_{i+1}$ are a.e convergent and $\Delta T_i\rightarrow 0$ a.e. Hence, the error terms vanish.

\vspace{0.5cm}

Turning to the diffusive terms, we must show that
\begin{align}
\nu D_+^x \rho \vv \sim (\bar \rho + \frac{1}{\hat\rho})(\bar \rho D_+^x \vv+ \bar \vv D_+^x \rho),\nonumber
\end{align}
converge weakly to the correct limits. (Here, and frequently below, we have dropped the indices altogether, since the calculations only involve to neighbouring points.) That is to say that the following error terms converge to zero:
\begin{align}
  \left| \overline{ \rho^2} D_+^x \vv-(\bar \rho)^2 D_+^x \vv\right|,\nonumber\\
 \left| \frac{1}{2}\bar \vv D_+^x \rho^2- \bar \vv \bar \rho D_+^x \rho\right|,\nonumber\\
\left|\frac{\bar \rho}{\hat\rho} D_+^x \vv-D_+^x \vv \right|,\nonumber\\
 \left| \frac{1}{\hat \rho}\bar \vv D_+^x \rho-\bar \vv D_+^x \log\rho\right|.\nonumber
\end{align}
In the first, both terms are bounded in $L^{1+\delta}_{1+\delta}$ by the a priori estimates. (We have $D\vv\in L^2_2$ and $\rho_i$ in a better space than $L^4_4$.) Hence, also the difference is bounded in the same space. Furthermore, $|\overline{ \rho^2} -(\bar \rho)^2|=\frac{1}{2}|(\Delta_+^x\rho)^2|$ which vanishes a.e. such that it converges strongly in a space better than $L^4_4$. Hence, the error is a product of a weakly and a strongly convergent sequence, such that the product converges weakly to zero.

The second error term is identically zero, since $D_+^x \rho_i^2 =2\bar \rho_{i+1/2} D_+^x \rho_i$ and both terms are bounded in $L^{1+\delta}_{1+\delta}$ by (\ref{imp_v}) and (\ref{disc_strong_rho}).

In the third, we have a direct bound on $D_+^x \vv\in L^2_2$. Furthermore, $1/\hat \rho\in L^{2+\delta}_{2+\delta} $ thanks to the $L^{\infty}_1$ and $L^2(H^1)$ bounds of $\rho_i^{-1}$, and since $\hat \rho_{i+1/2}$ is a monotone average bounded below and above by $\rho_{i}$ and $\rho_{i+1}$. As before, we have $\bar \rho D_+^x \vv\in L^2_2$ and therefore $\frac{\bar \rho}{\hat\rho} D_+^x \vv \in L^{1+\delta}_{1+\delta}$. Hence, the difference is also bounded in $ L^{1+\delta}_{1+\delta}$. Seeing that $|\frac{\bar\rho}{\hat\rho}-1|= |\frac{\bar\rho-\hat\rho}{\hat\rho}|\leq \frac{|\Delta_+^x\rho|}{\hat\rho} $ approaches zero a.e., we have the desired convergence.

The last term is identically zero, since  $\frac{ D_+^x \rho}{\hat \rho}= D_+^x\log \rho$.

\vspace{0.1cm}

Turning to the artificial diffusion, we consider,
\begin{align}
\lambda_{i+1/2jk}\Delta_+^x(\rho u)_{ijk}.\label{mom_ad}
\end{align}
This term is bounded in $L^{1+\epsilon}_{1+\epsilon}$ since $\lambda_h\sim |\log\rho_h||u|_h\in L^{2+\delta}_{2+\delta}$ (by (\ref{logr_Lp}) and (\ref{imp_v})) and by using (\ref{imp_mom}) to bound $\rho_h u_h$. To see that (\ref{mom_ad}) vanishes a.e., we recast it as,
\begin{align}
|\lambda_{i+1/2jk}\Delta_+^x(\rho u)_{ijk}|\leq\lambda_{i+1/2jk}(|\bar \rho_{i+1/2jk}|\Delta_+^x u_{ijk}|+ |\bar u_{i+1/2jk}||\Delta_+^x\rho_{ijk}|).\nonumber
\end{align}
Since all factors converge a.e., and the differences to zero a.e., the artificial diffusion term is a.e. converging to zero. Hence, it vanishes in $L^{1+\epsilon}_{1+\epsilon}$, for some $\epsilon>0$.

\vspace{0.5cm}

{\bf Total energy:}

The equation for total energy is the most difficult one. We consider some representative terms and demonstrate convergence one by one. The remaining terms follow the same pattern.

\vspace{0.1cm}

\textsc{Time derivative} [$p_{ijk}/(\gamma-1)+\frac{1}{2}\rho_{ijk}|\vv_{ijk}|^2$]: The pressure part converges strongly (at least in $L^2_2$) thanks to (\ref{imp_p}) and a.e. convergence of $\rho_h,T_h$. The kinetic energy is bounded by \eqref{imp_ke} and converges a.e., and therefore also strongly in $L^{1+\epsilon}_{1+\epsilon}$ ($0<\epsilon <\delta$).

\vspace{0.1cm}

\textsc{ inviscid flux}: The x-flux is given by:
\begin{align}
    \fb^{c,5}_{i+1/2jk}=&\frac{1}{2(\gamma-1)\hat \beta_{i+1/2jk}}\overline {\rho u}_{i+1/2jk}-\frac{\overline{|\vv|^2}_{i+1/2jk}}{2}\overline {\rho u}_{i+1/2jk}\nonumber \\
&+ \overline{|\vv|}^2_{i+1/2jk}\overline {\rho u}_{i+1/2jk}+p_{i+1/2jk}\bar u_{i+1/2jk}    \label{fc5}
\end{align}
Consider the internal-energy flux, $IE_1=\frac{1}{2(\gamma-1)\hat\beta_{i+1/2}}\overline{\rho u}_{i+1/2}$ (dropping the recurring $jk$ indices.). Our aim is to prove consistency with $\frac{\overline{pu}_{i+1/2}}{\gamma-1}$. That is, we need to show that
\begin{align}
E_1=|IE_1-\frac{\overline{pu}_{i+1/2}}{\gamma-1}|\rightarrow 0.\nonumber
\end{align}
Or,
\begin{align}
  E_1\leq |IE_1-\frac{R\bar T_{i+1/2}\overline{\rho u}_{i+1/2}}{\gamma-1}|+
|\frac{R\bar T_{i+1/2}\overline{\rho u}_{i+1/2}}{\gamma-1}-  \frac{\overline{pu}_{i+1/2}}{\gamma-1}|\rightarrow 0.\nonumber
\end{align}
The available a priori estimates provide bounds on all terms (in at least $L^{1+\epsilon}_{1+\epsilon}$). Hence, both differences are bounded in the same space. The second one is easily handled by using (\ref{split_av}),
\begin{align}
|\frac{R\bar T_{i+1/2}\overline{\rho u}_{i+1/2}}{\gamma-1}-  \frac{\overline{pu}_{i+1/2}}{\gamma-1}|\lesssim |\Delta_+^x T_i||\Delta_+^x (\rho u)_i|.\nonumber
\end{align}
This approaches zero a.e since $|\Delta_+^x T_i|$ does and $\rho u$ converges a.e.

Turning to the first part of the error, we need an estimate of 
\begin{align}
A=\bar T-\frac{1}{\widehat{T^{-1}}}= \bar T -\frac{1}{\frac{\Delta 1/T}{\Delta \log T^{-1}}}&=\nonumber \\ 
\bar T -\frac{1}{\frac{\Delta T^{-1}}{-\Delta \log T}}=
\bar T +\frac{\Delta \log T}{\Delta T^{-1}}&=\nonumber \\
\bar T +\frac{\Delta \log T}{-\frac{\Delta T}{\check T^2}}=
\bar T -\frac{\check T^2}{\hat T}&=
\frac{\bar T\hat T-\check T^2}{\hat T}\label{T_diffA1},
\end{align}
where again we have suppressed the indices for brevity.
We introduce $T_*^2=\bar T\hat T\in [\hat T^2,\bar T^2]$, to get,
\begin{align}
A=  \frac{(T_*+\check T)(T_*-\check T)}{\hat T}\leq
  \frac{(T_*+\check T)|\Delta T|}{\frac{\Delta T}{\Delta \log T}}\leq
  (T_*+\check T)|\Delta \log T|.
  \label{T_diffA2}
\end{align}
Hence,
\begin{align}
|IE_1-\frac{R\bar T_{i+1/2}\overline{\rho u}_{i+1/2}}{\gamma-1}|\sim (T_*+\check T)_{i+1/2}|\Delta_+ \log T_i|(\overline{\rho u})_{i+1/2}.\nonumber
  \end{align}
We already know that this difference is equi-integrable and we conclude that it vanishes since all factors converge a.e, and in particular $|\Delta T|\rightarrow 0$. (The term $p_{i+1/2jk}\bar u_{i+1/2jk}$ in (\ref{fc5}) can be handled in the same way.)

Next, we turn to the convective part of the inviscid flux. The kinetic-energy fluxes appearing in (\ref{fc5}) are bounded in $L^{1+\epsilon}_{1+\epsilon}$ by the same arguments leading to (\ref{KE_flux_est}) (seeing that  both terms are bounded by products of $\rho \vv$ and $|\vv|^2$.).

The numerical kinetic-energy flux is compared with the simple average, $ \frac{1}{2}\overline{\rho |\vv|^2u}_{i+1/2jk}$, which is bounded (as in (\ref{KE_flux_est})) and converges to the correct limit in the weak formulation thanks to the a.e. convergence of $\rho$ and $\vv$.

The error term that should vanish is thus:
\begin{align}
  E_2= \frac{1}{2}\overline{\rho |\vv|^2u}_{i+1/2jk}-\left(-\frac{\overline{|\vv|^2}_{i+1/2jk}}{2}\overline {\rho u}_{i+1/2jk}
+ \overline{|\vv|}^2_{i+1/2jk}\overline {\rho u}_{i+1/2jk}\right).\nonumber
\end{align}
Since both the numerical flux and the target flux are bounded in $L^{1+\epsilon}_{1+\epsilon}$, the error $E_2$ is bounded in the same space.

Using (\ref{split_av}), we recast the error,
\begin{align}
 E_2= \frac{1}{2}\overline{\rho u}\overline{ |\vv|^2}+\frac{\Delta\rho u\Delta |\vv|^2}{8}-\left(-\frac{\overline{|\vv|^2}}{2}\overline {\rho u}
    + \overline{|\vv|}^2\overline {\rho u}\right).\nonumber
\end{align}
Noting that $\Delta |\vv|^2=2\bar \vv D\vv$ and $\overline{ |\vv|^2}- \overline{|\vv|}^2=\frac{1}{4} (|\Delta \vv|)^2$, the error terms take the form
\begin{align}
  \frac{\Delta\rho u\Delta |\vv|^2}{8}&=  \frac{\Delta(\rho u)\bar\vv\Delta\vv}{4},\nonumber\\
\overline{\rho u}(\overline{ |\vv|^2}- \overline{|\vv|}^2)&= \frac{1}{4}\overline{\rho u}(|\Delta \vv|)^2\nonumber.
\end{align}
Both of these terms can be bounded as products of $\rho \vv$ and $|\vv|^2$ and are thus equi-integrable. We see that they vanish a.e. since they are products of primitive variables and differences.

\vspace{0.1cm}

\textsc{ Diffusive flux:}

We begin by considering the kinetic-energy part of the flux in the x-direction:
\begin{align}
\tilde \nu_{i+1/2jk}\left(\frac{1}{2}(D_+^x\rho_{ijk} |\vv_{ijk}|^2)+(\overline{|\vv|}^2_{i+1/2jk}-\overline{|\vv|^2}_{i+1/2jk})D_+^x\rho_{ijk}\right)\label{copy_diff_ke}
\end{align}
The last two terms two terms approximate the same quantity and we begin by showing that their sum converges to zero (in $L^1_1$). 
Since $\overline{ |\vv|^2}- \overline{|\vv|}^2=\frac{1}{4} (|\Delta \vv|)^2$, we get,
\begin{align}
  \tilde \nu(\overline{|\vv|}^2 -\overline{| \vv|^2})D_+^x\rho  =  \frac{\tilde \nu}{4} (|\Delta_+^x \vv|)^2D_+^x\rho  =h\frac{\tilde \nu}{4} (|D_+^x \vv|)^2\Delta_+^x\rho  \nonumber
\end{align}
Using (\ref{ke_nu}), we conclude that the term vanishes in $L^1_1$. (Note that the $h$-factor also makes it equi-integrable.)

Next, we turn to the part of (\ref{copy_diff_ke}) that should converge to the correct average as given in Def. \ref{def:weak}. That is, the term
\begin{align}
  \tilde \nu_{i+1/2jk}\frac{1}{2}D_+^x\rho_{ijk} |\vv_{ijk}|^2.\label{KE_diff_term}
\end{align}
First, we handle the terms arising from the $\bar \rho$ dependence of $\tilde \nu$. We also limit the analysis to one velocity component as the others can be handled analogously. We have (when suppressing the indices),
\begin{align}
  \bar \rho D (\rho u^2) =\bar \rho (\bar \rho Du^2+ \overline{u^2}D\rho).\label{KE_diff1}
\end{align}
The first term of (\ref{KE_diff1}):
\begin{align}
\bar \rho \bar \rho Du^2 = \bar \rho\bar \rho \bar u Du = (\bar \rho \bar u)(\bar \rho Du)= (\overline{\rho u}-\frac{\Delta\rho\Delta u}{4})\bar \rho Du.\label{ke_err1}
\end{align}
We use that $\bar \rho Du\in L^2_2$ and $\overline{\rho u}\in L^{2+\delta}_{2+\delta}$. We need the error term to vanish. To this end we note that $\Delta\rho\Delta u = h(\rho_1D u -\rho_2D u)$ and $\rho_{1,2}Du$ is bounded in $L^2_2$. Furthermore, the $h$-factor makes (\ref{ke_err1}) both equi-integrable and ensures that it vanishes in $L^1_1$.

Furthermore, $\overline{\rho u}\bar \rho Du$ approximates $\overline{\rho^2 u} Du$ up to an error term,
\begin{align}
\Delta(\rho u)\Delta \rho Du = \bar u(\Delta \rho)^2 Du +\bar \rho(\Delta \rho) (\Delta u)Du.\nonumber
\end{align}
The last term is the same as in (\ref{ke_err1}) and handled in the same way. The second last is rewritten as:
\begin{align}
  \bar u(\Delta \rho)^2Du = \bar u\Delta \rho (\Delta \rho Du) = (\Delta (\rho u)- (\Delta u) \bar \rho) (\Delta \rho Du).&\nonumber
\end{align}
The $\Delta (\rho u)(\Delta \rho Du)$ is equi-integrable by  $\bar \rho Du\in L^2_2$ and $\rho u\in L^{2+\delta}_{2+\delta}$ (see \ref{imp_mom}) and (\ref{disc_v_est})). It vanishes a.e. since $|(\Delta \rho u)(\Delta \rho Du)|\lesssim |\rho u||\Delta \rho| |Du|$ where $|Du|$ and $\rho u$ converges a.e and $|\Delta \rho|\rightarrow 0$ a.e. Finally, the $(\bar \rho\Delta u) (\Delta \rho Du)$ is the same as a previously handled term.

\vspace{0.2cm}

The second term of (\ref{KE_diff1}) is $\overline{u^2}\bar \rho D\rho=\frac{1}{2}\overline{u^2} D\rho^2 = \frac{1}{2}D(\rho^2u^2)-\overline{\rho^2}\bar uDu$. (This is the numerical counterpart of (\ref{odd_weak_form}).)  
In view of (\ref{odd_weak_form}) and by using (\ref{split_av}), the $\overline{\rho^2}\bar uDu$-term results in an error term,
\begin{align}
-(\overline{\rho^2}\bar u-\overline{\rho^2 u})Du=\frac{\Delta \rho^2\Delta u}{4}Du= \frac{1}{2}(\Delta \rho\Delta u)(\bar\rho Du)\nonumber.
\end{align}
We have already handled a term of this form and know that it vanishes.

The above considerations show that the $\bar \rho$-dependent part of $\tilde \nu$ in the kinetic-energy flux can be recast as
\begin{align}
\frac{1}{4}D(\rho^2u^2)-\frac{1}{2}\overline{\rho^2u}Du+\textrm{error terms},\nonumber
\end{align}
where the error terms all vanish. The remaining part converges to the correct limit by the same arguments as in Section \ref{sec:weak_seq}. In short, in the term $\frac{1}{4}D(\rho^2u^2)$ we move the difference to the test function which results in a momentum term  that is equi-integrable and a.e. convergent. Moreover, recasting the second term as
\begin{align}
\overline{\rho^2u}Du = \frac{1}{2}\left(\rho^2_1u_1+\rho^2_2u_2\right)Du &=\nonumber \\  \frac{1}{2}\left( (\rho_1u_1)(\rho_1Du)+(\rho_2u_2)(\rho_2)\right)Du &,\nonumber 
\end{align}
($1,2$ denotes the two points appearing in the average and differences) and using the same estimates as above, shows that it is equi-integrable. To demonstrate convergence to the correct limit, we argue that $\bar \rho Du$ converges weakly to the correct limit as in Section \ref{sec:weak_seq}. Strong convergence of momentum, completes the argument.

Second, we handle the terms arising from the $(\hat \rho)^{-1}$ dependence of $\tilde \nu$ in (\ref{KE_diff_term}). Namely,
\begin{align}
  \frac{1}{\hat \rho} D (\rho u^2)= \frac{1}{\hat \rho} \left(\overline{u^2}D \rho + 2\bar\rho\bar u Du\right)=  \overline{u^2}D \log \rho + 2\frac{\bar\rho}{\hat\rho}\bar u Du&=\nonumber\\
  D(u^2 \log \rho) - 2\overline{\log\rho} \bar u Du+2\frac{\bar\rho}{\hat\rho}\bar u Du.\nonumber
\end{align}
The difference in the $D(u^2 \log \rho)$ is moved to the test function. It is bounded in $L^{1+\epsilon}_{1+\epsilon}$ and consistent with the corresponding term in Def. \ref{def:weak}. In the second, $\bar u Du$ is bounded in $L^{1+\epsilon}_{1+\epsilon}$ and weakly convergent (by the analogous argument as $\rho\nabla \vv$ is weakly convergent. See end of Section \ref{sec:compact}.) Since $\log\rho$ can be bounded in any $L^p_p$ space, and is a.e. convergent, it is also strongly convergent in a sufficiently good space. This ensures convergence to the right limit of $\overline{u\log\rho}  Du$. The error term $|\overline{\log\rho} \bar u Du-\overline{u\log\rho}  Du|$ vanishes, since it is majorised by $|\Delta u||\Delta\log\rho||Du|$. This is equi-integrable thanks to $\bar u Du\in L^{1+\epsilon}_{1+\epsilon} $ and the $L^p_p$ estimates of $\log \rho$. It approaches zero a.e, since all factors are a.e. convergent and $|\Delta\log\rho|\rightarrow 0$.

Finally, the last term should converge to $\bar u Du$, which is weakly convergent in $L^1_1$. Hence, we must show that the following error term vanishes:
\begin{align}
  |(1-\frac{\bar\rho}{\hat\rho})\bar u Du|=
  |\frac{\hat\rho-\bar\rho}{\hat\rho}\bar u Du|&\leq \nonumber \\
  |\frac{\Delta \rho}{\hat\rho}\bar u Du|&=
  |(\Delta \log \rho)\bar u Du|.&\nonumber
\end{align}
We have already shown that such a term is equi-integrable and vanishes.

Turning to the artificial diffusion part of (\ref{KE_diff_term}), we need that
\begin{align}
h\lambda D(\rho u^2) \sim h(|u|+|\Delta \log\rho||u|+|\Delta u|) D(\rho u^2),\nonumber
\end{align}
 vanishes in $L^1_1$. We begin with
\begin{align}
h|\Delta u|D(\rho u^2)  = |\Delta u|\Delta (\rho u^2).\label{ke_err2}
\end{align}
Let the subscripts $1,2$ denote the two points on which the flux term is based. Then,
\begin{align}
  (u_1-u_2)(\rho_1u_1^2-\rho_2u_2^2)=
  \rho_1u_1^3-\rho_2u_2^3-\rho_1u_1^2u_2-\rho_2u_2^2u_1&=\nonumber\\
  \rho_1u_1^3-\rho_2u_2^3-\rho_1u_1(u_1u_2)-\rho_2u_2(u_2u_1)&\leq \nonumber\\
  \rho_1u_1^3+\rho_2u_2^3-\rho_1u_1(u_1^2+u_2^2)-\rho_2u_2(u_2^2+u_1^2), \nonumber
\end{align}
shows that (\ref{ke_err2}) can be majorised by $\overline{\rho u}\overline{u^2}$, which can be split in the same way leading to the bound (\ref{imp_ke_flux}). Hence, it is equi-integrable. Moreover, $|\Delta u|\rightarrow 0$ a.e. ensures that (\ref{ke_err2}) vanishes a.e. 

The same argument applies directly to $h|u|D(\rho u^2)$.

This leaves us with the artificial diffusion term,
\begin{align}
|\Delta \log\rho||\bar u|\Delta (\rho u^2).\label{err_ad}
\end{align}
As previously shown, $|\bar u|\Delta (\rho u^2)\in L^{1+\epsilon}_{1+\epsilon}$ and $\Delta \log\rho$ can be bounded in a sufficiently high $L^p_p$ space to obtain equi-integrability of the error term (\ref{err_ad}). Moreover, (\ref{err_ad}) vanishes a.e. since $\Delta (\rho u^2) = (\Delta \rho) \overline{u^2}+2\bar u \bar \rho \Delta u$ and the differences vanish.

\vspace{0.2cm}

\textsc{Diffusive flux, internal energy  [$\tilde \nu_{i+1/2jk}\frac{(\pf_x)_{i+1/2jk}}{\gamma-1}$]:}

\begin{align}
 \tilde \nu_{i+1/2jk} (\pf_x)_{i+1/2jk} &= \tilde \nu_{i+1/2jk}\left(\frac{1}{2\hat \beta_{i+1/2jk}}D_+^x\rho_{ijk}+\frac{\bar\rho_{i+1/2jk} }{2}D_+^x\frac{1}{\beta_{ijk}}\right)\label{num_KE_diff_flux}
\end{align}
The ``target terms'' that, up to some constants, should be approximated are,
\begin{align}
  \rho T D\rho,\quad
  \rho^2  DT,\quad
  TD\log\rho,\quad
  DT\nonumber.
\end{align}
We begin with the terms associated with $\tilde \nu\sim \bar\rho$. They generate error terms of the form:
\begin{align}
  \left(\frac{\bar\rho}{2\hat\beta}-\overline{R\rho T} \right)D\rho&,\label{IE_T1}\\
  \left(\bar\rho^2-\overline{\rho^2} \right)DT&= -\frac{1}{4}\left(\rho_1-\rho_2 \right)^2DT.\label{IE_T2}
\end{align}
The second term, (\ref{IE_T2}), is bounded by $\rho \in L^4_4$ and $DT\in L^2_2$. Furthermore, $DT$ converges weakly and $(\Delta \rho)^2$ strongly since $|\Delta\rho|\rightarrow 0$ a.e. Hence, the error term (\ref{IE_T2}) vanishes in $L^1_1$. The term $\overline{\rho^2} DT$ converges to the correct limit thanks to the strong convergence of $\rho$ and weak convergence of $DT$.

Turning to (\ref{IE_T1}), we make the following manipulations,
\begin{align}
  \left(\frac{\bar\rho}{2\hat\beta}-\overline{R\rho T} \right)=
  \left(\frac{\bar\rho}{2\hat\beta}-R\bar\rho \bar T +R\frac{\Delta\rho\Delta T}{4} \right)=
    \left(\bar{\rho}\left(\frac{1}{2\hat\beta}-R\bar T\right) +R\frac{\Delta\rho\Delta T}{4} \right)\label{IE_err}
\end{align}
Inserting (\ref{IE_err}) in (\ref{IE_T1}), we see that we need to control $\Delta\rho\Delta T D\rho$, or equivalently, $(\Delta \rho)^2D T$. This is the same as term as (\ref{IE_T2}), which we have already handled. This leaves us with $\bar{\rho}\left(\frac{1}{2\hat\beta}-R\bar T\right)D\rho$ (from (\ref{IE_err}) and (\ref{IE_T1})). To handle this term, we use (\ref{T_diffA1}) and (\ref{T_diffA2}), to obtain
\begin{align}
  |\frac{1}{2R\hat \beta}-\bar T|\lesssim \bar T|\Delta \log T|.\nonumber
\end{align}
Hence, we must control $\bar \rho (D\rho)\bar T |\Delta \log T|$. $\bar \rho D\rho=\frac{1}{2}D\rho^2\in L^2_2$ which converges weakly. We need strong convergence (to zero) of $|\bar T\Delta \log T|$ in a better space than $L^2_2$. We note that it converges a.e. to zero and only need the bound. This follows by interpolating the bounds  $T^{3/2}D\log T\in L^2_2$ (see (\ref{strong_DlogT})) and $D\log T\in L^2_2$ (see (\ref{rlogr_d})) such that $TD\log T\in L^2_2$.  This implies that $|\bar T\Delta \log T| \in L^{2+\delta}_{2+\delta}$, $\delta>0$. Finally, the target entity $\overline{R\rho T} D\rho$ converges to the correct limit thanks to the strong convergence of $p\in L^{2+\delta}_{2+\delta}$ (following from (\ref{imp_p}) and a.e. convergence) and the weak convergence of $D\rho\in L^2_2$.

Next, we consider error terms of (\ref{num_KE_diff_flux}) for $\tilde \nu\sim 1/\hat\rho$:
\begin{align}
  \left|\left(\frac{1}{2\hat\beta}-\overline{T} \right)\frac{D\rho}{\hat\rho}\right|=\left|\left(\frac{1}{2\hat\beta}-\overline{T} \right)D\log\rho\right|&\lesssim \bar T |\Delta\log T||D\log\rho| \nonumber\\
  \left|\left(1-\frac{\bar \rho}{\hat\rho} \right)DT\right|\leq\left|\left(\frac{|\Delta \rho| }{\hat\rho} \right)DT\right|&\leq|\Delta\log\rho | | DT|\sim |D\log\rho| |\Delta T|.
  \nonumber
\end{align}
The first is handled in the same way as above, only now we are using that $D\log \rho \in L^2_2$. The second is handled by the strong bound on $T$ and a.e. convergence, implying strong convergence of $|\Delta T|$ (to zero), and $D\log\rho\in L^2_2$. 
 Here, the target is $\bar T D\log\rho$ and it is easy to verify convergence to the correct limit.

The error terms of (\ref{num_KE_diff_flux}) for $\tilde \nu\sim h((\frac{1}{2}+|\Delta \log\rho|)|\bar u|+\Delta u)$ (artificial diffusion):
\begin{align}
  h((\frac{1}{2}+|\Delta \log\rho|)|\bar u|+|\Delta u|)\left( \frac{D\rho}{\hat \beta}+\bar\rho DT\right)&=\nonumber\\
 ((\frac{1}{2}+ |\Delta \log\rho|)|\bar u|+|\Delta u|)\left( \frac{\Delta\rho}{\hat \beta}+\bar\rho \Delta T\right)\nonumber
\end{align}
The strong bounds imply that
\begin{align}
  \left( \frac{\Delta\rho}{\hat \beta}+\bar\rho \Delta T\right)&\in L^{2+\epsilon}_{2+\epsilon}\quad \textrm{and}\nonumber\\
  (\frac{1}{2}+|\Delta \log\rho|)|\bar u|+|\Delta u|&\in L^2_2.\nonumber
\end{align}
(Note that $\Delta T\in L^4_4$ and $1/\hat\beta \in L^4_4$.) This gives equi-integrability.  The a.e. vanishing differences ensures that the error terms disappear.

\vspace{0.2cm}

\textsc{ Radiation diffusion} [$\kappa_rD_-^xD_+^xT^4_{ijk}$]: This term is straightforwardly handled by moving the differences to the test function. This results in a $T^4_h$ term, which is strongly convergent thanks to (\ref{disc_strong_T}) and the a.e convergence of temperature.

\vspace{0.5cm}

{\bf Entropy inequality:}  The global entropy diffusivity condition in Def. \ref{def:weak} follows directly from (\ref{ent_est1}), the strong bounds on the solution and the positive definiteness of the diffusive terms appearing in the left-hand side of (\ref{ent_est1}).

\vspace{0.5cm}

\section{Concluding remarks}\label{sec:final}

The main result of this paper (Theorem \ref{main_theo}) is the existence of weak entropy solutions to the alternative Navier-Stokes system (\ref{eulerian}). These weak solutions ensure positivity of temperature and density, except possibly on a set of Lebesgue measure zero. To the best of our knowledge, such results are not available for the standard Navier-Stokes-Fourier system.

The weak solutions are obtained as the limits of solutions to a finite volume scheme. We make some further remarks on the scheme:
\begin{itemize}

\item As presented, the scheme is formally first-order accurate for smooth solutions (on regular grids). This is due to the first-order artificial-diffusion coefficent, namely the ``$1/2$'' appearing in  (\ref{ad_coeff2}), which results in an ``upwind''-type diffusion. As this is only needed in the entropy estimate (see (\ref{2nd_order})), we can replace $\lambda$ in (\ref{ad_coeff2}) by
  \begin{align}
    \lambda_{i+1/2jk}&=|\bar u_{i+1/2jk}|R^\#_{i+1/2jk}+ \frac{|\Delta_+^xu_{ijk}|}{4}\nonumber \\
    R^\#_{i+1/2jk}&=\max\left(\left|\frac{\bar \rho_{i+1/2jk}-\hat \rho_{i+1/2jk}}{\Delta _+^x\rho_{ijk}} \right|,|\Delta_+^x\log\rho_{ijk}|\right)\nonumber
  \end{align}
  The first entry in the maximum is $\mathcal{O}(|\Delta_+^x\rho_{ijk}|)$. (Formally, this implies second-order accuracy.) Since $R^\#$ is dominated by $R^*$ (c.f. (\ref{2nd_order})), and converges to zero a.e., it does not affect the convergence proof.

  However, in order to prevent that the finite arithmetic pollutes the numerical solutions, the first entry of $R^\#$ requires a well-conditioned numerical approximation. This should be obtainable in the same way that the approximation of the log mean was derived (see Appendix B of \cite{IsmailRoe09}), but we leave this as future work.

\item  Unfortunately, the generalisation to high-order accuracy (three or higher for smooth solutions), discovered in \cite{FisherCarpenter13} is effectively ruled out, since it requires the use of entropy-conservative fluxes, which in turn appears incompatible with the current approach.

\item The scheme is trivially generelisable to the unstructured finite-volume framework developed by Eymard et al. (\cite{Eymard_etal97}) for Voronoi control volumes and two-point fluxes. The convergence proof holds in this case as well, thanks to the summation-by-parts property of such schemes.

\item The extension of this theory to allow some form of far-field boundary condition to close the domain for external flows, appears to be non-trivial. The ideas used in \cite{Svard21} may be a possible route to this end.
\end{itemize}

Future work, also includes investigating if (\ref{eulerian}), like many other fluid systems, has some weak-strong uniqueness property.

\section{Acknowledgements}

I am very grateful to Prof. Josef M{\'a}lek whose help has enabled me to write this article.

\bibliographystyle{alpha}
\bibliography{ref}

\newcommand{\etalchar}[1]{$^{#1}$}
\begin{thebibliography}{RDO{\etalchar{+}}19}

\bibitem[Bre05a]{Brenner05_1}
H.~Brenner.
\newblock Kinematics of volume transport.
\newblock {\em Physica A}, 349:11, 2005.

\bibitem[Bre05b]{Brenner05_2}
H.~Brenner.
\newblock Navier-{S}tokes revisited.
\newblock {\em Physica A}, 349:60, 2005.

\bibitem[Cha13]{Chandrashekar13}
P.~Chandrashekar.
\newblock Kinetic energy preserving and entropy stable finite volume schemes
  for compressible {E}uler and {N}avier-{S}tokes equations.
\newblock {\em Comm. Comp. Phys.}, 14(5):1252--1256, 2013.

\bibitem[DS21]{DolejsiSvard21}
V{\'i}t Dolej{\v s}{\'i} and Magnus Sv{\"a}rd.
\newblock Numerical study of two models for viscous compressible fluid flows.
\newblock {\em J. Comp. Phys}, 427:110068, 2021.

\bibitem[EGH97]{Eymard_etal97}
R.~Eymard, T.~Gallou{\"e}t, and R.~Herbin.
\newblock Finite volume methods.
\newblock In P.G. Ciarlet and J.L. Lions, editors, {\em Handbook of Numerical
  Analysis}, volume~7, pages 713--1020. Elsevier, 1997.

\bibitem[Eva10]{Evans}
L.~C. Evans.
\newblock {\em Partial Differential Equations, 2nd ed.}
\newblock American Mathematical Society, 2010.

\bibitem[FC13]{FisherCarpenter13}
T.~C. Fisher and M.~H. Carpenter.
\newblock High-order entropy stable finite difference schemes for nonlinear
  conservation laws: Finite domains.
\newblock {\em J. Comput. Phys.}, 252:518--557, 2013.

\bibitem[FN17]{FeireislNovotny}
E.~Feireisl and A.~Novotn\'y.
\newblock {\em Singular limits in thermodynamics of viscous fluids, 2nd Ed.}
\newblock Birkh\"auser, 2017.

\bibitem[FV10]{FeireislVasseur10}
Eduard Feireisl and Alexis Vasseur.
\newblock {\em New Perspectives in Fluid Dynamics: Mathematical Analysis of a
  Model Proposed by Howard Brenner}, pages 153--179.
\newblock Birkh{\"a}user Basel, Basel, 2010.

\bibitem[GSH21]{GassnerSvard21}
G.J. Gassner, M.~Sv{\"a}rd, and F.J. Hindenlang.
\newblock Stability issues of entropy-stable and/or split-form high-order
  schemes.
\newblock {\em to appear in Journal of Scientific Computing}, 2021.

\bibitem[IR09]{IsmailRoe09}
F.~Ismail and P.L. Roe.
\newblock Affordable, entropy-consistent euler flux functions {II}: {E}ntropy
  production at shocks.
\newblock {\em J. Comp. Phys.}, pages 5410--5436, 2009.

\bibitem[NFAE03]{NordstromForsberg03}
J.~Nordstr\"om, K.~Forsberg, C.~Adamsson, and P.~Eliasson.
\newblock Finite volume methods, unstructured meshes and strict stability for
  hyperbolic problems.
\newblock {\em Applied Numerical Mathematics}, 45(4):453--473, June 2003.

\bibitem[RDO{\etalchar{+}}19]{Reddy_etal19}
M.~H.~L. Reddy, S.~K. Dadzie, R.~Ocone, M.~K. Borg, and J.~M. Reese.
\newblock Recasting {N}avier-{S}tokes equations.
\newblock {\em Journ. Phys. Comm.}, 3(10):1--25, 2019.

\bibitem[SDP21]{SayyariDalcin21}
M.~Sayyari, L.~Dalcin, and M.~Parsani.
\newblock Development and analysis of entropy stable no-slip wall boundary
  conditions for the {E}ulerian model for viscous and heat conducting
  compressible flows.
\newblock {\em SN PDEA}, 2, 2021.

\bibitem[Sim86]{Simon86}
J.~Simon.
\newblock Compact sets in the space {$L^p(O,T;B)$}.
\newblock {\em Annali di Matematica Pura ed Applicata}, 146:65--96, 1986.

\bibitem[SN14]{SvardNordstrom14}
M.~Sv{\"a}rd and J.~Nordstr{\"o}m.
\newblock Review of summation-by-parts schemes for initial-boundary-value
  problems.
\newblock {\em J. Comput. Phys.}, 268:17--38, 2014.

\bibitem[Sto45]{Stokes1845}
G.~G. Stokes.
\newblock On the theories of internal friction of fluids in motion, and of the
  equilibrium and motion of elastic solids.
\newblock {\em Transactions of the Cambridge Philosophical Society},
  8:287--305, 1845.

\bibitem[Sv{\"a}18]{Svard18}
M.~Sv{\"a}rd.
\newblock A new {E}ulerian model for viscous and heat conducting compressible
  flows.
\newblock {\em Physica A}, 506:350--375, 2018.

\bibitem[Sv{\"a}21]{Svard21}
M.~Sv{\"a}rd.
\newblock Entropy stable boundary conditions for the {E}uler equations.
\newblock {\em J. Comput. Phys.}, 426:1--17, 2021.

\end{thebibliography}

\appendix

\section{Prerequisites} \label{app1}

An elementary identity: $\|u^r\|_p^q=\|u\|_{rp}^{rq}$.

\vspace{0.5cm}

A list of standard inequalities:

\begin{itemize}

\item Cauchy-Schwarz inequality: If $u,v\in L^2$ then $<u,v>\leq \|u\|_2\|v\|_2$.

\item H\"older's inequality: $\|uv\|_1\leq \|u\|_p\|v\|_q$ where $(p,q)\in [1,\infty)$ are H\"older conjugates, $1/p+1/q=1$.

\item A generalisation of H\"older's inequality (see \cite{FeireislNovotny}): $\|uv\|_r\leq \|u\|_p\|v\|_q$ where $(p,q,r)\in [1,\infty)$ are H\"older conjugates, $1/p+1/q=1/r$.

\item Minkowski (triangle) inequality: $\|u+v\|_p\leq \|u\|_p+\|v\|_p$ for $p\in[1,\infty)$.

\item Young's inequality: For $a,b\in \RR$, $a\cdot b\leq \frac{a^p}{p}+\frac{b^q}{q}$, where $(p,q)$ is the H\"older conjugate.

\item Young's inequality is often combined with H\"older in the following way:
\begin{align}
\|uv\|_1\leq \epsilon\|u\|_p\epsilon^{-1}\|v\|_q\leq \frac{(\epsilon\|u\|)^p}{p}+\frac{(\epsilon^{-1}\|v\|)^q}{q}.
\end{align}
  
\item Ladyzhenskaya's inequality in 3-D: $\|u\|_4\leq \C \|u\|^{1/4}_2\|\nabla u\|_2^{3/4}$.

\item Nash' inequality in $\RR^n$: $\|u\|_2^{1+2/n}\leq \C \|u\|^{2/n}_1\|Du\|_2$.

\item Sobolev embedding (in 3-D): $\|u\|_6\leq  \C(\|u\|_2+\|\nabla u\|_2)$.

\item Riesz-Thorin (interpolation) inequality:
\begin{align}
  \|u\|_r\leq \|u\|_p^{1-\theta}\|u\|^\theta_q\nonumber \\
  \frac{1}{r}=\frac{1-\theta}{p}+\frac{\theta}{q},\quad 0<\theta<.1\nonumber
\end{align}
An important special case of this inequality is the ``10/3''-rule,
\begin{align}
\|u\|^{10/3}_{10/3}\leq \|u\|^{4/3}_2\|u\|^2_6\label{interp1}
\end{align}
which is often used in combination with Sobolev embedding. Some other special cases:
\begin{align}
\|u\|^{20/3}_{20/3}&\leq \|u\|^{2/3}_1\|u\|^6_{18}\label{interp2} \\
\|u\|^{3}_{3}&\leq \|u\|^{3/5}_1\|u\|^{12/5}_{6}\label{interp3}  \\
\|u\|^{6}_{6}&\leq \|u\|^{2/3}_3\|u\|^{4}_{12}\label{interp4} 
\end{align}

\item A special case of the Gagliardo-Nirenberg inequality:
\begin{align}
\|u\|_3\leq \C\|u\|_2^{1/2}\|\nabla u\|_2^{1/2}\nonumber
\end{align}

\end{itemize}

A generalised Poincare inequality: 
\begin{theorem}\label{theo:poincare}
  Let $1\leq p \leq \infty$, $0<\Gamma<\infty$, $V_0>0$ and let $\Omega\subset \RR^N$ be a bounded Lipschitz domain. Then there exists a positive constant $c=(p,\Gamma,V_0)$ such that,
  \begin{align}
    \|v\|_{W^{1,p}}(\Omega)\leq  c\left[ \|\nabla v\|_{L^p(\Omega;R^N)}+(\int_V|v|^{\Gamma}d\xb)^{\frac{1}{\Gamma}}\right]\label{poincare}
    \end{align}
\end{theorem}
for any measurable $V\subset \Omega$, $|V|>V_0$ and any $v\in W^{1,p}(\Omega)$.
\begin{proof}
This is Theorem 11.20 in \cite{FeireislNovotny}, where it is also proved.
\end{proof}

Aubin-Lions Lemma (see e.g. \cite{Simon86}):
\begin{lemma}\label{lemma:aubin}
For Banach spaces $X\subset B \subset Y$ where the embedding of $X$ in $B$ is compact and $B$ in $Y$ is continuous. Let $U=\{u|u\in L^p(0,\T;X)\,\textrm{and}\,u_t\in L^q(0,\T;Y)\}$ where $1\leq p,q <\infty$. Then $U$ is compactly embedded in $L^p(0,\T;B)$.
\end{lemma}

The following Lemma for weakly and strongly converging sequences is a standard result.
\begin{lemma}
  Let $1<p<\infty$ and $1/p+1/q =1$ and $\Omega \subset \RR^N$. If
  \begin{align}
    u_n \rightharpoonup u \quad \textrm{in} \quad L^p(\Omega),\nonumber \\
    v_n \rightarrow v \quad \textrm{in} \quad L^q(\Omega),\nonumber
  \end{align}
then 
    \begin{align}
    u_nv_n \rightharpoonup uv \quad \textrm{in} \quad L^1(\Omega),\nonumber
  \end{align}
  \end{lemma}

The following Lemma, that we state without proof, is a consequence of Vitali's Convergence Theorem:
\begin{lemma}
  Let $\Omega\subset \RR^N$ be bounded. If $u_n\rightarrow u$ a.e. in $\Omega$ and $u_n$ is bounded in $L^p(\Omega)$, $p>1$. Then, $u_n\rightarrow u$ in $L^r$ for all $1\leq r<p$.
  \end{lemma}

The next theorem is found in \cite{Evans} (Appendix D, Theorem 3).
\begin{theorem}
  Let $u_n\in L^p(\Omega)$, $1<p<\infty$ be a uniformly bounded sequence. Then, there is a subsequence (still denoted $\{u_n\}$) and a function $u\in L^p(\Omega)$, such that $u_n\rightharpoonup u$ in $L^p(\Omega)$.
\end{theorem}


\end{document}